%% file: DAG-rootfile.tex
\documentclass{article}
\usepackage{latexsym, amscd, amsfonts, eucal, mathrsfs, amsmath, amssymb, amsthm, xypic,xr, makeidx, stmaryrd}
\def\Z{\mathbf{Z}}

\DeclareMathOperator{\Spectra}{Sp}
\DeclareMathOperator{\PreSpectra}{PSp}

\newcommand{\Andre}{Andr\'{e}}

\DeclareMathOperator{\im}{im}
\DeclareMathOperator{\Gap}{Gap}

\newcommand{\LLPres}[1]{\mathcal{P}{\text r}^{\text L}_{#1}}

\newcommand{\RRPres}[1]{\mathcal{P}{\text r}^{\text R}_{#1}}
\newcommand{\dplus}{-}
\newcommand{\dminus}{+}
\DeclareMathOperator{\Exc}{Exc}
\DeclareMathOperator{\EM}{\mathcal{E}\mathcal{M}}

\DeclareMathOperator{\Rexx}{Rex}

\DeclareMathOperator{\connSpectra}{\Spectra^{conn}}
\DeclareMathOperator{\RPres}{\mathcal{P}r^{R}}
\DeclareMathOperator{\Group}{\mathcal{G}rp}

\DeclareMathOperator{\rght}{R}
\DeclareMathOperator{\lft}{L}

\newcommand{\toposref}[1]{T.\ref{HTT-#1}}
\newcommand{\stableref}[1]{S.\ref{STA-#1}}

\newcommand{\degree}{\text{o}}
\newcommand{\bfA}{{\mathbf A}}

\DeclareMathOperator{\Kan}{\mathcal{K}an}
\DeclareMathOperator{\cDelta}{{\bf \Delta}}

\DeclareMathOperator{\Set}{\mathcal{S}et}
\DeclareMathOperator{\sSet}{\mathcal{S}et_{\Delta}}

\DeclareMathOperator{\Nerve}{N}
\DeclareMathOperator{\Cat}{\mathcal{C}at}
\newcommand{\h}[1]{\rm{h} \! #1}
\newcommand{\heart}[1]{{#1}^{\heartsuit}}
\newcommand{\Adjoint}[4]{\xymatrix@1{#2 \ar@<.4ex>[r]^-{#1} & #3 \ar@<.4ex>[l]^-{#4}}}

\renewcommand{\boxtimes}{\odot}

\DeclareMathOperator{\coker}{coker}
\DeclareMathOperator{\Chain}{Ch}

\newcommand{\et}{\'{e}t}

\newcommand{\bigdot}{\bullet}

\DeclareMathOperator{\Stab}{Stab}

\newcommand{\Ab}{ \mathcal{A}b}
\newcommand{\fin}{\text{fin}}
\DeclareMathOperator{\bHom}{Map}

\newcommand{\Cech}{\v{C}ech\,}

\DeclareMathOperator{\Mod}{\mathcal{M}od}

\DeclareMathOperator{\calM}{\mathcal{M}}

 \DeclareMathOperator{\G}{\mathbf{G}}

\DeclareMathOperator{\Ext}{Ext} 
 
\DeclareMathOperator{\DerRing}{\mathcal{SCR}}

\DeclareMathOperator{\Ex}{Ex}
\DeclareMathOperator{\Q}{\mathbf{Q}}
 
\DeclareMathOperator{\colim}{colim}

\DeclareMathOperator{\calA}{\mathcal{A}}
\DeclareMathOperator{\bd}{\partial}

\newcommand{\Sphere}{S}
\DeclareMathOperator{\LPres}{\mathcal{P}r^{L}}

\DeclareMathOperator{\calE}{\mathcal{E}}

\DeclareMathOperator{\calO}{\mathcal{O}}

\DeclareMathOperator{\FinSpace}{\mathcal{S}_{\ast}^{fin}}
\DeclareMathOperator{\FinnSpace}{\mathcal{S}^{fin}}

\DeclareMathOperator{\FinSpectra}{\mathcal{S}_{\infty}^{fin}}
\DeclareMathOperator{\calB}{\mathcal{B}}

\DeclareMathOperator{\Spec}{{\bf Spec}}

\DeclareMathOperator{\Hom}{Hom} 
\DeclareMathOperator{\HH}{H} 
\DeclareMathOperator{\id}{id} \DeclareMathOperator{\Fun}{Fun}
\DeclareMathOperator{\calC}{\mathcal{C}}
\DeclareMathOperator{\calI}{\mathcal{I}}

\DeclareMathOperator{\calJ}{\mathcal{J}}

\DeclareMathOperator{\SSet}{\mathcal{S}}

\DeclareMathOperator{\calD}{\mathcal{D}}
\DeclareMathOperator{\Ind}{Ind} 
\DeclareMathOperator{\calP}{\mathcal{P}} \topmargin=0in

\newtheorem{theorem}{Theorem}[section]
\newtheorem{lemma}[theorem]{Lemma}
\newtheorem{proposition}[theorem]{Proposition}

\newtheorem{corollary}[theorem]{Corollary}

\theoremstyle{definition}
\newtheorem{definition}[theorem]{Definition}
\newtheorem{construction}[theorem]{Construction}
\newtheorem{example}[theorem]{Example}

\newtheorem{notation}[theorem]{Notation}

\newtheorem{warning}[theorem]{Warning}
\newtheorem{remark}[theorem]{Remark}

\makeindex

\begin{document}

\title{Derived Algebraic Geometry I: Stable $\infty$-Categories}

\maketitle
\tableofcontents

\include{DAG-I-stable}

\end{document}

%% file: DAG-I-stable.tex
 \section{Introduction}\label{stable1}
 
There is a very useful analogy between topological spaces and 
chain complexes with values in an abelian category. For example, it is customary to speak of
{\em homotopies} between chain maps, {\em contractible} complexes, and so forth. The analogue of the homotopy category of topological spaces is the {\em derived category} of
an abelian category $\calA$, a triangulated category which provides a good
setting for many constructions in homological algebra. However, it has long been recognized that for many purposes the derived category is too crude: it identifies homotopic morphisms of chain complexes without remembering {\em why} they are homotopic. It is possible to correct this
defect by viewing the derived category as the homotopy category of an underlying $\infty$-category $\calD(\calA)$. The $\infty$-categories which arise in this way have special features that reflect their ``additive'' origins: they are {\it stable}. 

The goal of this paper is to provide an introduction to the theory of stable $\infty$-categories. We will begin in \S \ref{stable2} by introducing the definition of stability and some other basic terminology. In many ways, an arbitrary stable $\infty$-category $\calC$ behaves like the derived category of an abelian category: in particular, we will see in \S \ref{stable3} that for every stable $\infty$-category $\calC$, the homotopy category $\h{\calC}$ is triangulated (Theorem \ref{surmite}). In \S \ref{stable4} we will establish some other simple consequences of stability; for example, stable $\infty$-categories admit finite limits and colimits (Proposition \ref{surose}). 

The appropriate notion of {\em functor} between stable $\infty$-categories is an {\em exact} functor: that is, a functor which preserves finite colimits (or equivalently, finite limits: see Proposition \ref{funrose}). The collection of stable $\infty$-categories and exact functors between them can be organized into an $\infty$-category, which we will denote by $\Cat_{\infty}^{\Ex}$. In \S \ref{stable5}, we will study the $\infty$-category $\Cat_{\infty}^{\Ex}$; in particular, we will show that it is stable under limits and filtered colimits in $\Cat_{\infty}$. The formation of limits in $\Cat_{\infty}^{\Ex}$ provides a tool for addressing the classical problem of ``gluing in the derived category''.

In \S \ref{stable6}, we will review the theory of t-structures on triangulated categories. We will see that, if $\calC$ is a stable $\infty$-category, there is a close relationship between t-structures on the homotopy category $\h{\calC}$ and localizations of $\calC$. We will revisit this subject in \S \ref{stable16}, where we show that, under suitable set-theoretic hypotheses (to be described in \S \ref{stable15}), we can construct a t-structure ``generated'' by an arbitrary collection of objects of $\calC$. 

The most important example of a stable $\infty$-category is the $\infty$-category $\Spectra$ of {\em spectra}. The homotopy category of $\Spectra$ can be identified with the classical {\it stable homotopy category}. There are many approaches to the construction of $\Spectra$. In \S \ref{stable8} we will adopt the most classical perspective: we begin by constructing an $\infty$-category $\FinSpectra$ of {\em finite spectra}, obtained from the $\infty$-category of finite pointed spaces by formally inverting the suspension functor. The stability of $\FinSpectra$ follows from the classical {\em homotopy excision theorem}.
We can then define the $\infty$-category $\Spectra$ as the $\infty$-category of $\Ind$-objects of $\FinSpectra$. The stability of $\Spectra$ follows from a general result on $\Ind$-objects (Proposition \ref{kappstable}). 

There is another description of the $\infty$-category $\Spectra$ which is perhaps more familiar: it can be viewed as the $\infty$-category of {\em infinite loop spaces}, obtained from the $\infty$-category $\SSet_{\ast}$ of pointed spaces by formally inverting the loop functor.
More generally, one can begin with an arbitrary $\infty$-category $\calC$, and construct a new $\infty$-category $\Stab(\calC)$ of {\it infinite loop objects of $\calC$}. The $\infty$-category $\Stab(\calC)$ can be regarded as universal among stable $\infty$-categories which admits a left exact functor to $\calC$ (Proposition \ref{urtusk21}). This leads to a characterization of $\Spectra$ by a mapping property: namely, $\Spectra$ is freely generated under colimits (as a stable $\infty$-category) by a single object, the {\it sphere spectrum} (Corollary \ref{choccrok}).

A classical result of Dold and Kan asserts that, if $\calA$ is an abelian category, then the category of simplicial objects in $\calA$ is equivalent to the category of nonnegatively graded chain complexes in $\calA$. In \S \ref{doldkan}, we will formulate and prove an $\infty$-categorical version of this result, where the abelian category $\calA$ is replaced by a stable $\infty$-category. Here we must replace the notion of ``chain complex'' by the related notion ``filtered object''. If 
$\calC$ is a stable $\infty$-category equipped with a t-structure, then every filtered object of $\calC$ determines a spectral sequence; we will give the details of this construction in \S \ref{filttt}.

In \S \ref{stable10}, we will return to the subject of homological algebra. We will explain how to pass from a suitable abelian category $\calA$ to a stable $\infty$-category $\calD^{\dplus}(\calA)$, which we will call the {\it derived $\infty$-category of $\calA$}. The homotopy category of $\calD^{\dplus}(\calA)$ can be identified with the classical derived category of $\calA$. 

Our final goal in this paper is to characterize $\calD^{\dplus}(\calA)$ by a universal mapping property. In \S \ref{stable14}, we will show that $\calD^{\dplus}(\calA)$ is universal among stable $\infty$-categories equipped with a suitable embedding of the ordinary category $\calA$ (Corollary \ref{sutty}).


The theory of stable $\infty$-categories is not really new: most of the results presented here are well-known to experts. There exists a sizable literature on the subject in the setting of {\it stable model categories} (see, for example, \cite{stablemodel}). The theory of stable model categories is essentially equivalent to the notion of a {\em presentable} stable $\infty$-category, which we discuss in \S \ref{stable15}. For a brief account in the more flexible setting of Segal categories, we refer the reader to \cite{toenK}.

In this paper, we will use the language of {\it $\infty$-categories} (also called quasicategories or weak Kan complexes), as described in \cite{topoi}. We will use the letter T to indicate references to \cite{topoi}. For example, Theorem \toposref{mainchar} refers to Theorem \ref{HTT-mainchar} of \cite{topoi}.

\section{Stable $\infty$-Categories}\label{stable2}  

In this section, we will introduce our main object of study: stable $\infty$-categories. We begin with a brief review of some ideas from \S \toposref{chmdim}.

\begin{definition}\index{object!zero}\index{zero!object}\index{$\infty$-category!pointed}\index{pointed $\infty$-category}
Let $\calC$ be an $\infty$-category. A {\it zero object} of $\calC$ is an object which is
both initial and final. We will say that $\calC$ is {\it pointed} if it contains a zero object.
\end{definition}

In other words, an object $0 \in \calC$ is zero if the spaces
$\bHom_{\calC}(X,0)$ and $\bHom_{\calC}(0,X)$ are both contractible
for every object $X \in \calC$. Note that if $\calC$ contains a zero object, then
that object is determined up to equivalence. More precisely, the 
full subcategory of $\calC$ spanned by the zero objects is a contractible Kan complex
(Proposition \toposref{initunique}).

\begin{remark}\label{tustpoint}
Let $\calC$ be an $\infty$-category. Then $\calC$ is pointed if and only if the following conditions are satisfied:
\begin{itemize}
\item[$(1)$] The $\infty$-category $\calC$ has an initial object $\emptyset$.
\item[$(2)$] The $\infty$-category $\calC$ has a final object $1$.
\item[$(3)$] There exists a morphism $f: 1 \rightarrow \emptyset$ in $\calC$.
\end{itemize}
The ``only if'' direction is obvious. For the converse, let us suppose that $(1)$, $(2)$, and $(3)$ are satisfied. We invoke the assumption that $\emptyset$ is initial to deduce the existence of a morphism $g: \emptyset \rightarrow 1$. Because $\emptyset$ is initial, $f \circ g \simeq \id_{\emptyset}$, and because $1$ is final, $g \circ f \simeq \id_{1}$. Thus $g$ is a homotopy inverse to $f$, so that $f$ is an equivalence. It follows that $\emptyset$ is also a final object of $\calC$, so that $\calC$ is pointed.
\end{remark}

\begin{remark}\index{zero!morphism}\index{morphism!zero}
Let $\calC$ be an $\infty$-category with a zero object $0$. For
any $X,Y \in \calC$, the natural map $$\bHom_{\calC}(X,0) \times
\bHom_{\calC}(0,Y) \rightarrow \bHom_{\calC}(X,Y)$$ has contractible
source. We therefore obtain a well defined morphism $X \rightarrow Y$
in the homotopy category $\h{\calC}$, which we will refer to as the {\it zero morphism}
and also denote by $0$. 
\end{remark}


\begin{definition}\index{triangle}\index{triangle!exact}\index{triangle!coexact}\index{exact!triangle}
Let $\calC$ be a pointed $\infty$-category. A {\it triangle} in $\calC$ is a diagram
$\Delta^1 \times \Delta^1 \rightarrow \calC$, depicted as
$$ \xymatrix{ X \ar[r]^{f} \ar[d] & Y \ar[d]^{g} \\
0 \ar[r] & Z }$$
where $0$ is a zero object of $\calC$. We will say that a triangle in $\calC$ is {\it exact} if
it is a pullback square, and {\it coexact} if it is a pushout square.
\end{definition}

\begin{remark}
Let $\calC$ be a pointed $\infty$-category. A triangle in $\calC$ consists of the following data:
\begin{itemize}
\item[$(1)$] A pair of morphisms $f: X \rightarrow Y$ and $g: Y \rightarrow Z$ in $\calC$.
\item[$(2)$] A $2$-simplex in $\calC$ corresponding to a diagram
$$ \xymatrix{ & Y \ar[dr]^{g} & \\
X \ar[ur]^{f} \ar[rr]^{h} & & Z }$$
in $\calC$, which identifies $h$ with the composition $g \circ f$.
\item[$(3)$] A $2$-simplex
$$ \xymatrix{ & 0 \ar[dr] & \\
X \ar[ur] \ar[rr]^{h} & & Z }$$
in $\calC$, which we may view as a {\em nullhomotopy} of $h$.
\end{itemize}
We will sometimes indicate a triangle by specifying only the pair of maps
$$ X \stackrel{f}{\rightarrow} Y \stackrel{g}{\rightarrow} Z,$$
with the data of $(2)$ and $(3)$ being implicitly assumed.
\end{remark}\index{nullhomotopy}

\begin{definition}\index{kernel}\index{cokernel}
Let $\calC$ be a pointed $\infty$-category containing a morphism $g: X \rightarrow Y$.
A {\it kernel} of $g$ is an exact triangle
$$ \xymatrix{ W \ar[r] \ar[d] & X \ar[d]^{g} \\
0 \ar[r] & Y. }$$
Dually, a {\it cokernel} for $g$ is a coexact triangle
$$ \xymatrix{ X \ar[r]^{g} \ar[d] & Y \ar[d] \\
0 \ar[r] & Z. }$$
We will sometimes abuse terminology by simply referring to $W$ and $Z$ as the kernel and cokernel of $g$. We will also write $W = \ker(g)$ and $Z = \coker(g)$.\index{ZZZker@$\ker(f)$}\index{ZZZcoker@$\coker(f)$}
\end{definition}

\begin{remark}\label{cokerdef}
Let $\calC$ be a pointed $\infty$-category containing a morphism $f: X \rightarrow Y$.
A cokernel of $f$, if it exists, is uniquely determined up to equivalence. More precisely,
consider the full subcategory $\calE \subseteq \Fun( \Delta^1 \times \Delta^1, \calC)$ spanned by the coexact triangles. Let $\theta: \calE \rightarrow \Fun(\Delta^1,\calC)$ be the forgetful functor, which associates to a diagram
$$ \xymatrix{ X \ar[r]^{g} \ar[d] & Y \ar[d] \\
0 \ar[r] & Z }$$
the morphism $g: X \rightarrow Y$. Applying Proposition \toposref{lklk} twice, we deduce that
$\theta$ is a Kan fibration, whose fibers are either empty or contractible (depending on whether or not a morphism $g: X \rightarrow Y$ in $\calC$ admits a cokernel). In particular, if every
morphism in $\calC$ admits a cokernel, then $\theta$ is a trivial Kan fibration, and therefore admits a section $\coker: \Fun(\Delta^1, \calC) \rightarrow \Fun(\Delta^1 \times \Delta^1, \calC)$, which is well defined up to a contractible space of choices. We will often abuse notation by also letting
$\coker: \Fun(\Delta^1, \calC) \rightarrow \calC$ denote the composition
$$ \Fun(\Delta^1, \calC) \rightarrow \Fun(\Delta^1 \times \Delta^1, \calC) \rightarrow \calC,$$where
the second map is given by evaluation at the final object of $\Delta^1 \times \Delta^1$.
\end{remark}

\begin{remark}\label{coprodunt}
The functor $\coker: \Fun(\Delta^1, \calC) \rightarrow \calC$ can be identified with
a left adjoint to the left Kan extension functor $\calC \simeq \Fun( \{1\}, \calC) \rightarrow \Fun( \Delta^1, \calC)$, which associates to each object $X \in \calC$ a zero morphism $0 \rightarrow X$. It follows that
$\coker$ preserves all colimits which exist in $\Fun(\Delta^1, \calC)$ (Proposition \toposref{adjointcol}). 
\end{remark}

\begin{definition}\label{stabl}\index{stable $\infty$-category}\index{$\infty$-category!stable}
An $\infty$-category $\calC$ is {\it stable} if it satisfies the
following conditions:

\begin{itemize}
\item[$(1)$] There exists a zero object $0 \in \calC$.
\item[$(2)$] Every morphism in $\calC$ admits a kernel and a cokernel.
\item[$(3)$] A triangle in $\calC$ is exact if and only if it is coexact.
\end{itemize}
\end{definition}

\begin{remark}
Condition $(3)$ of Definition \ref{stabl} is analogous to the
axiom for abelian categories which requires that the image of a
morphism be isomorphic to its coimage.
\end{remark}

\begin{example}
Recall that a {\it spectrum} consists of an infinite sequence of pointed topological spaces
$\{ X_i \}_{i \geq 0}$, together with homeomorphisms $X_i \simeq \Omega X_{i+1}$, where
$\Omega$ denotes the loop space functor. The collection of spectra can be organized into a stable $\infty$-category $\Spectra$. Moreover, $\Spectra$ is in some sense the universal example of a stable $\infty$-category. This motivates the terminology of Definition \ref{stabl}: an $\infty$-category $\calC$ is stable if it resembles the $\infty$-category $\Spectra$, whose homotopy category $\h{ \Spectra}$ can be identified with the classical {\em stable homotopy category}. We will return to the theory of spectra (using a slightly different definition) in \S \ref{stable8}. 
\end{example}

\begin{example}
Let $\calA$ be an abelian category. Under mild hypotheses, we can construct a stable $\infty$-category $\calD(\calA)$ whose homotopy category $\h{ \calD(\calA) }$ can be identified with the {\it derived category of $\calA$}, in the sense of classical homological algebra. We will outline the construction of $\calD(\calA)$ in \S \ref{stable10}.
\end{example}

\begin{remark}
If $\calC$ is a stable $\infty$-category, then the opposite
$\infty$-category $\calC^{op}$ is also stable.
\end{remark}

\begin{remark}\label{additivity}
One attractive feature of the theory of stable $\infty$-categories
is that stability is a property of $\infty$-categories, rather
than additional data. The situation for additive categories is similar. Although
additive categories are often presented as categories equipped
with additional structure (an abelian group structure on all
$\Hom$-sets), this additional structure is in fact determined by
the underlying category. If a category $\calC$ has a
zero object, finite sums, and finite products, then there always
exists a unique map $A \oplus B \rightarrow A \times B$ which can be described by the matrix $$ \left[ \begin{matrix} \id_A & 0  \cr 0 &
\id_B \cr \end{matrix} \right] \begin{matrix} \cr . \cr
\end{matrix}$$ If this morphism has an inverse $\phi_{A,B}$,
then we may define a sum of two morphisms $f,g: X \rightarrow Y$
by defining $f+g$ to be the composition $X \rightarrow X \times X
\stackrel{f,g}{\rightarrow} Y \times Y
\stackrel{\phi_{Y,Y}}{\rightarrow} Y \oplus Y \rightarrow Y$. This definition endows
each morphism set $\Hom_{\calC}(X,Y)$ with the structure of a commutative monoid. If
each $\Hom_{\calC}(X,Y)$ is actually a group (in other words, if every morphism
$f: X \rightarrow Y$ has an additive inverse), then $\calC$ is an additive category. This statement has an analogue in the setting of stable $\infty$-categories: any stable $\infty$-category $\calC$ is automatically enriched over the $\infty$-category of spectra. Since we do not wish to develop the language of enriched $\infty$-categories, we will not pursue this point further.
\end{remark}
 
 \section{The Homotopy Category of a Stable $\infty$-Category}\label{stable3} 
 
Our goal in this section is to show that if $\calC$ is a stable $\infty$-category, then the homotopy category $\h{\calC}$ is triangulated (Theorem \ref{surmite}). We begin by reviewing the definition of a triangulated category.

\begin{definition}[Verdier]\label{deftriangle}\index{category!triangulated}\index{triangulated category}\index{triangle!distinguished}\index{distinguished triangle}\index{Verdier's axioms}
A {\it triangulated category} consists of the following data:
\begin{itemize}
\item[$(1)$] An additive category $\calD$.
\item[$(2)$] A translation functor $$\calD \rightarrow \calD$$
$$ X \mapsto X[1], $$
which is an equivalence of categories.
\item[$(3)$] A collection of {\it distinguished triangles}
$$ X \stackrel{f}{\rightarrow} Y \stackrel{g}{\rightarrow} Z \stackrel{h}{\rightarrow} X[1].$$
\end{itemize}
These data are required to satisfy the following axioms:
\begin{itemize}
\item[$(TR1)$]
\begin{itemize}
\item[$(a)$] Every morphism $f: X \rightarrow Y$ in $\calD$ can be extended to distinguished triangle in $\calD$.
\item[$(b)$] The collection of distinguished triangles is stable under isomorphism.
\item[$(c)$] Given an object $X \in \calD$, the diagram
$$ X \stackrel{\id_{X}}{\rightarrow} X \rightarrow 0 \rightarrow X[1]$$
is a distinguished triangle.
\end{itemize}
\item[$(TR2)$] A diagram
$$ X \stackrel{f}{\rightarrow} Y \stackrel{g}{\rightarrow} Z \stackrel{h}{\rightarrow} X[1]$$
is a distinguished triangle if and only if the rotated diagram
$$ Y \stackrel{g}{\rightarrow} Z \stackrel{h}{\rightarrow} X[1] \stackrel{ - f[1]}{\rightarrow} Y[1]$$
is a distinguished triangle.
\item[$(TR3)$] Given a commutative diagram 
$$ \xymatrix{ X \ar[r] \ar[d]^{f} & Y \ar[r] \ar[d] & Z \ar@{-->}[d] \ar[r] & X[1] \ar[d]^{f[1]} \\
X' \ar[r] & Y' \ar[r] & Z' \ar[r] & X'[1] }$$
in which both horizontal rows are distinguished triangles, there exists a dotted arrow rendering the
entire diagram commutative.
\item[$(TR4)$] Suppose given three distinguished triangles\index{octahedral axiom}
$$ X \stackrel{f}{\rightarrow} Y \stackrel{u}{\rightarrow} Y/X \stackrel{d}{\rightarrow} X[1] $$
$$ Y \stackrel{g}{\rightarrow} Z \stackrel{v}{\rightarrow} Z/Y \stackrel{d'}{\rightarrow} Y[1] $$
$$ X \stackrel{g \circ f}{\rightarrow} Z \stackrel{w}{\rightarrow} Z/X \stackrel{d''}{\rightarrow} X[1]$$
in $\calD$. There exists a fourth distinguished triangle
$$ Y/X \stackrel{\phi}{\rightarrow} Z/X \stackrel{\psi}{\rightarrow} Z/Y \stackrel{\theta}{\rightarrow}
Y/X[1]$$
such that the diagram
$$ \xymatrix{ X \ar[rr]^{g \circ f} \ar[dr]^{f} & & Z \ar[dr]^{w} \ar[rr]^{v} & & Z/Y \ar[dr]^{d'}\ar[rr]^{\theta} & & Y/X[1] \\ 
& Y \ar[dr]^{u} \ar[ur]^{g} & & Z/X \ar[ur]^{\psi} \ar[dr]^{d''} & & Y[1] \ar[ur]^{u[1]} & \\
& & Y/X \ar[ur]^{\phi} \ar[rr]^{d} & & X[1] \ar[ur]^{f[1] } & & }$$
commutes.
\end{itemize}
\end{definition}

\begin{remark}
The theory of triangulated categories is an attempt to capture those features of stable $\infty$-categories which are visible at the level of homotopy categories. Triangulated categories which appear naturally in mathematics are usually equivalent to the homotopy categories of suitable stable $\infty$-categories.
\end{remark}

We now consider the problem of constructing a triangulated structure on the homotopy category of an $\infty$-category $\calC$. To begin the discussion, let us assume that $\calC$ is an arbitrary pointed $\infty$-category. We $\calM^{\Sigma}$ denote the full subcategory of 
$\Fun( \Delta^1 \times \Delta^1, \calC)$ spanned by those diagrams
$$ \xymatrix{ X \ar[r] \ar[d] & 0 \ar[d] \\
0' \ar[r] & Y }$$
which are pushout squares, and such that $0$ and $0'$ are zero objects of $\calC$.
If $\calC$ admits cokernels, then we can use Proposition \toposref{lklk} (twice) to conclude that evaluation at the initial vertex induces a trivial fibration $\calM^{\Sigma} \rightarrow \calC$. Let $s: \calC \rightarrow \calM^{\Sigma}$ be a section of this trivial fibration, and let $e: \calM^{\Sigma} \rightarrow \calC$ be the functor given by evaluation at the final vertex. The composition $e \circ s$ is a functor from $\calC$ to itself, which we will denote by $\Sigma: \calC \rightarrow \calC$ and refer to as the {\it suspension functor} on $\calC$. Dually, we define $\calM^{\Omega}$ to be the full
subcategory of $\Fun( \Delta^1 \times \Delta^1, \calC)$ spanned by diagrams as above which are pullback squares with $0$ and $0'$ zero objects of $\calC$. If $\calC$ admits kernels, then 
the same argument shows that evaluation at the final vertex induces a trivial fibration $\calM^{\Omega} \rightarrow \calC$. If we let $s'$ denote a section to this trivial fibration, then the composition of $s'$ with\index{suspension functor}\index{loop functor}
evaluation at the initial vertex induces a functor from $\calC$ to itself, which we will refer to as the {\em loop functor} and denote by $\Omega: \calC \rightarrow \calC$. If $\calC$ is stable, then
$\calM^{\Omega} = \calM^{\Sigma}$. It follows that $\Sigma$ and $\Omega$ are mutually inverse equivalences from $\calC$ to itself.

\begin{remark}\index{ZZZSigma@$\Sigma$}\index{ZZZOmega@$\Omega$}
If the $\infty$-category $\calC$ is not clear from context, then we will denote the suspension and loop functors $\Sigma, \Omega: \calC \rightarrow \calC$ by $\Sigma_{\calC}$ and $\Omega_{\calC}$, respectively.
\end{remark}

\begin{notation}\label{nutus}\index{shift functor}\index{ZZZ[]@$X[n]$}
If $\calC$ is a stable $\infty$-category and $n \geq 0$, we let
$$ X \mapsto X[n]$$
denote the $n$th power of the suspension functor $\Sigma: \calC \rightarrow \calC$
constructed above (this functor is well-defined up to canonical equivalence). If $n \leq 0$, we let $X \mapsto X[n]$ denote the $(-n)$th power of the loop functor $\Omega$. We will use the same notation to indicate the induced functors on the homotopy category $\h{\calC}$. 
\end{notation}

\begin{remark}
If the $\infty$-category $\calC$ is pointed but not necessarily stable, the suspension and loop space functors need not be homotopy inverses but are nevertheless {\em adjoint} to one another (provided that both functors are defined).
\end{remark}

If $\calC$ is a pointed $\infty$-category containing a pair of objects $X$ and $Y$, then the space
$\bHom_{\calC}(X,Y)$ has a natural base point, given by the zero map. Moreover, if $\calC$
admits cokernels, then the suspension functor $\Sigma_{\calC}: \calC \rightarrow \calC$ is essentially characterized by the existence of natural homotopy equivalences
$$ \bHom_{\calC}( \Sigma(X), Y) \rightarrow \Omega \bHom_{\calC}( X,Y).$$
In particular, we conclude that $\pi_0 \bHom_{\calC}( \Sigma(X), Y) \simeq
\pi_1 \bHom_{\calC}(X,Y)$, so that $\pi_0 \bHom_{\calC}( \Sigma(X),Y)$ has the structure of a group (here the fundamental group of $\bHom_{\calC}(X,Y)$ is taken with base point given by the zero map). Similarly, $\pi_0 \bHom_{\calC}( \Sigma^2(X),Y) \simeq \pi_2 \bHom_{\calC}(X,Y)$ has the structure of an {\em abelian group}. If the suspension functor $X \mapsto \Sigma(X)$ is an equivalence of $\infty$-categories, then
for every $Z \in \calC$ we can choose $X$ such that $\Sigma^2(X) \simeq Z$ to deduce the existence of an abelian group structure on $\bHom_{\calC}(Z,Y)$. It is easy to see that this group structure depends functorially
on $Z,Y \in \h{\calC}$. We are therefore most of the way to proving the following result:

\begin{lemma}\label{vival}
Let $\calC$ be a pointed $\infty$-category which admits cokernels, and suppose that the suspension functor $\Sigma: \calC \rightarrow \calC$ is an equivalence.
Then $\h{\calC}$ is an additive category.
\end{lemma}

\begin{proof}
The argument sketched above shows that $\h{\calC}$ is (canonically) enriched over the category of abelian groups. It will therefore suffice to prove that $\h{\calC}$ admits finite coproducts.
We will prove a slightly stronger statement: the $\infty$-category $\calC$ itself admits finite coproducts. Since $\calC$ has an initial object, it will suffice to treat the case of pairwise coproducts.
Let $X, Y \in \calC$, and let $\coker: \Fun(\Delta^1, \calC) \rightarrow \calC$ be a cokernel functor, so that we have equivalences $X \simeq \coker( X[-1] \stackrel{u}{\rightarrow} 0)$ and $Y \simeq \coker{ 0 \stackrel{v}{\rightarrow} Y }$. Proposition \toposref{limiteval} implies that $u$ and $v$ admit a coproduct in $\Fun(\Delta^1, \calC)$ (namely, the zero map $X[-1] \stackrel{0}{\rightarrow} Y$). Since the functor $\coker$ preserves coproducts (Remark \ref{coprodunt}), we conclude that
$X$ and $Y$ admit a coproduct (which can be constructed as the cokernel of the zero map
from $X[-1]$ to $Y$).
\end{proof}

Let $\calC$ be a pointed $\infty$-category which admits cokernels.
By construction, any diagram
$$ \xymatrix{ X \ar[r] \ar[d] & 0 \ar[d] \\
0' \ar[r] & Y }$$
which belongs to $\calM$ determines a canonical isomorphism $X[1] \rightarrow Y$ in the homotopy category $\h{\calC}$. We will need the following observation:

\begin{lemma}\label{purpleato}
Let $\calC$ be a pointed $\infty$-category which admits cokernels, and let
$$ \xymatrix{ X \ar[r]^{f} \ar[d]^{f'} & 0 \ar[d] \\
0' \ar[r] & Y }$$
be a diagram in $\calC$, classifying a morphism
$\theta \in \Hom_{ \h{\calC} }( X[1], Y)$. $($Here $0$ and $0'$ are zero objects of $\calC$.$)$
Then the transposed diagram
$$ \xymatrix{ X \ar[r]^{f'} \ar[d]^{f} & 0' \ar[d] \\
0 \ar[r] & Y }$$
classifies the morphism $-\theta \in \Hom_{\h{\calC}}(X[1], Y)$. Here 
$-\theta$ denotes the inverse of $\theta$ with respect to the group structure on
$\Hom_{\h{\calC}}(X[1], Y) \simeq \pi_1 \bHom_{\calC}(X,Y)$.
\end{lemma}

\begin{proof}
Without loss of generality, we may suppose that $0 = 0'$ and $f = f'$. Let
$\sigma: \Lambda^2_0 \rightarrow \calC$ be the diagram
$$ 0 \stackrel{f}{\leftarrow} X \stackrel{f}{\rightarrow} 0.$$
For every diagram $p: K \rightarrow \calC$, let $\calD(p)$ denote the Kan complex $\calC_{p/} \times_{\calC} \{Y\}$. Then 
$\Hom_{\h{\calC}}(X[1],Y) \simeq \pi_0 \calD(\sigma)$. We note that
$$ \calD(\sigma) \simeq \calD(f) \times_{ \calD(X) } \calD(f).$$
Since $0$ is an initial object of $\calC$, $\calD(f)$ is contractible. In particular, there
exists a point $q \in \calD(f)$. Let 
$$ \calD' = \calD(f) \times_{ \Fun( \{0\}, \calD(X) )} \Fun(\Delta^1, \calD(X))
\times_{ \Fun( \{1\}, \calD(X) )} \calD(f)$$
$$ \calD'' = \{q\} \times_{ \Fun( \{0\}, \calD(X)) } \Fun(\Delta^1, \calD(X))
\times_{ \Fun( \{1\}, \calD(X)) } \{q\}$$
so that we have canonical inclusions
$$ \calD'' \hookrightarrow \calD' \hookleftarrow \calD(\sigma).$$
The left map is a homotopy equivalence because $\calD(f)$ is contractible, and the right
map is a homotopy equivalence because the projection $\calD(f) \rightarrow \calD(X)$ is a Kan fibration. We observe that $\calD''$ can be identified with the simplicial loop space of 
$\Hom^{\lft}_{\calC}(X,Y)$ (taken with the base point determined by $q$, which we can identify with the zero map from $X$ to $Y$). Each of the Kan complexes $\calD(\sigma)$, $\calD'$, $\calD''$ is equipped with a canonical involution. On $\calD(\sigma)$, this involution corresponds to the transposition of diagrams as in the statement of the lemma. On $\calD''$, this involution corresponds to reversal of loops. The desired conclusion now follows from the observation that these involutions are compatible with the inclusions $\calD'', \calD(\sigma) \subseteq \calD'$.
\end{proof}

\begin{definition}\label{surtato}\index{triangle!distinguished}\index{distinguished triangle}
Let $\calC$ be a pointed $\infty$-category which admits cokernels.
Suppose given a diagram
$$ X \stackrel{f}{\rightarrow} Y \stackrel{g}{\rightarrow} Z \stackrel{h}{\rightarrow} X[1]$$
in the homotopy category $\h{\calC}$. We will say that this diagram is a {\it distinguished triangle}
if there exists a diagram $\Delta^1 \times \Delta^2 \rightarrow \calC$ as shown
$$ \xymatrix{ X \ar[r]^{\widetilde{f}} \ar[d] & Y \ar[d]^{\widetilde{g}} \ar[r] & 0 \ar[d] \\
0' \ar[r] & Z \ar[r]^{\widetilde{h}} & W, }$$
satisfying the following conditions:
\begin{itemize}
\item[$(i)$] The objects $0, 0' \in \calC$ are zero.
\item[$(ii)$] Both squares are pushout diagrams in $\calC$.
\item[$(iii)$] The morphisms $\widetilde{f}$ and $\widetilde{g}$ represent $f$ and $g$, respectively.
\item[$(iv)$] The map $h: Z \rightarrow X[1]$ is the composition of
(the homotopy class of) $\widetilde{h}$ with the isomorphism $W \simeq X[1]$ determined by the outer rectangle.
\end{itemize}
\end{definition}

\begin{remark}
We will generally only use Definition \ref{surtato} in the case where $\calC$ is a stable $\infty$-category. However, it will be convenient to have the terminology available in the case where $\calC$ is not yet known to be stable.
\end{remark}

The following result is an immediate consequence of Lemma \ref{purpleato}:

\begin{lemma}\label{tater}
Let $\calC$ be a stable $\infty$-category. Suppose given a diagram
$\Delta^2 \times \Delta^1 \rightarrow \calC$, depicted as
$$ \xymatrix{ X \ar[d]^{f} \ar[r] & 0 \ar[d] \\
Y \ar[d] \ar[r]^{g} & Z \ar[d]^{h} \\
0' \ar[r] & W, }$$
where both squares are pushouts and the objects $0,0' \in \calC$ are zero.
Then the diagram
$$ X \stackrel{f}{\rightarrow} Y \stackrel{g}{\rightarrow} Z \stackrel{-h'}{\rightarrow} 
X[1] $$
is a distinguished triangle in $\h{\calC}$, where $h'$ denotes the composition of $h$ with
the isomorphism $W \simeq X[1]$ determined by the outer square, and $-h'$ denotes the composition of $h'$ with the map $-\id \in \Hom_{\h{\calC}}( X[1],X[1] )
\simeq \pi_1 \bHom_{\calC}(X, X[1])$.
\end{lemma}

We can now state the main result of this section:

\begin{theorem}\label{surmite}
Let $\calC$ be a pointed $\infty$-category which admits cokernels, and suppose that the suspension functor $\Sigma$ is an equivalence. Then the translation functor of Notation \ref{nutus} and the class of distinguished triangles of Definition \ref{surtato} endow $\h{\calC}$ with the structure of a triangulated category.
\end{theorem}

\begin{remark}
The hypotheses of Theorem \ref{surmite} hold whenever $\calC$ is stable. In fact, the hypotheses of Theorem \ref{surmite} are {\em equivalent} to the stability of $\calC$: see Corollary \ref{charstut}.
\end{remark}

\begin{proof}
We must verify that Verdier's axioms $(TR1)$ through $(TR4)$ are satisfied.
\begin{itemize}
\item[$(TR1)$] Let $\calE \subseteq \Fun( \Delta^1 \times \Delta^2, \calC)$ be the full subcategory spanned by those diagrams 
$$ \xymatrix{ X \ar[r]^{f} \ar[d] & Y \ar[d] \ar[r] & 0 \ar[d] \\
0' \ar[r] & Z \ar[r] & W }$$
of the form considered in Definition \ref{surtato}, and let $e: \calE \rightarrow \Fun( \Delta^1, \calC)$
be the restriction to the upper left horizontal arrow. Repeated use of Proposition \toposref{lklk} 
implies $e$ is a trivial fibration. In particular, every morphism $f: X \rightarrow Y$
can be completed to a diagram belonging to $\calE$. This proves $(a)$. 
Part $(b)$ is obvious, and $(c)$ follows from the observation that
if $f = \id_{X}$, then the object $Z$ in the above diagram is a zero object of $\calC$.

\item[$(TR2)$] Suppose that
$$ X \stackrel{f}{\rightarrow} Y \stackrel{g}{\rightarrow} Z \stackrel{h}{\rightarrow} X[1]$$
is a distinguished triangle in $\h{\calC}$, corresponding to a diagram
$\sigma \in \calE$ as depicted above. Extend $\sigma$ to a diagram
$$ \xymatrix{ X \ar[r] \ar[d] & Y \ar[d] \ar[r] & 0 \ar[d] \\
0' \ar[r] \ar[d] & Z \ar[r] \ar[d] & W \ar[d]^{u} \\
0'' \ar[r] & 0''' \ar[r] & V }$$
where the lower right square is a pushout, and the objects $0'',0''' \in \calC$ are zero.
We have a map between the squares
$$ \xymatrix{ X \ar[r] \ar[d] & 0 \ar[d] & Y \ar[r] \ar[d] & 0 \ar[d] \\
0' \ar[r] & W & 0''' \ar[r] & V }$$
which induces a commutative diagram in the homotopy category $\h{\calC}$
$$ \xymatrix{ W \ar[d]^{u} \ar[r] & X[1] \ar[d]^{f[1]} \\
V \ar[r] & Y[1] }$$
where the horizontal arrows are isomorphisms.
Applying Lemma \ref{tater} to the rectangle on the right of the large diagram, 
we conclude that
$$ Y \stackrel{g}{\rightarrow} Z \stackrel{h}{\rightarrow} X[1] \stackrel{-f[1]}{\rightarrow} Y[1]$$
is a distinguished triangle in $\h{\calC}$.

Conversely, suppose that $$ Y \stackrel{g}{\rightarrow} Z \stackrel{h}{\rightarrow} X[1] \stackrel{-f[1]}{\rightarrow} Y[1]$$ is a distinguished triangle in $\h{\calC}$. Since the functor $\Sigma: \calC \rightarrow \calC$ is an equivalence, we conclude that the triangle
$$Y[-2] \stackrel{g[-2]}{\rightarrow} Z[-2] \stackrel{h[-2]}{\rightarrow} X[-1] \stackrel{-f[-1]}{\rightarrow} Y[-1]$$ is distinguished. Applying the preceding argument five times, we conclude that
the triangle $$ X \stackrel{f}{\rightarrow} Y \stackrel{g}{\rightarrow} Z \stackrel{h}{\rightarrow} X[1]$$ is distinguished, as desired.

\item[$(TR3)$] Suppose distinguished triangles
$$ X \stackrel{f}{\rightarrow} Y \rightarrow Z \rightarrow X[1] $$
$$ X' \stackrel{f'}{\rightarrow} Y' \rightarrow Z' \rightarrow X'[1]$$
in $\h{\calC}$. Without loss of generality, we may suppose that these triangles are induced
by diagrams $\sigma, \sigma' \in \calE$. Any commutative diagram
$$ \xymatrix{ X \ar[r]^{f} \ar[d] & Y \ar[d] \\
X' \ar[r]^{f'} & Y' }$$
in the homotopy category $\h{\calC}$ can be lifted (nonuniquely) to a square in
$\calC$, which we may identify with a morphism $\phi: e(\sigma) \rightarrow e(\sigma')$
in the $\infty$-category $\Fun( \Delta^1, \calC)$. Since $e$ is a trivial fibration of simplicial
sets, $\phi$ can be lifted to a morphism $\sigma \rightarrow \sigma'$ in $\calE$, which determines a natural transformation of distinguished triangles
$$ \xymatrix{ X \ar[d] \ar[r] & Y \ar[d] \ar[r] & Z \ar[r] \ar[d] & X[1] \ar[d] \\
X' \ar[r] & Y' \ar[r] & Z' \ar[r] & X'[1]. }$$

\item[$(TR4)$] \index{octahedral axiom}Let $f: X \rightarrow Y$ and $g: Y \rightarrow Z$ be morphisms in $\calC$.
In view of the fact that $e: \calE \rightarrow \Fun(\Delta^1, \calC)$ is a trivial fibration, any distinguished triangle in $\h{\calC}$ beginning with $f$, $g$, or $g \circ f$ is uniquely determined up to (nonunique) isomorphism. Consequently, it will suffice to prove that there exist
{\em some} triple of distinguished triangles which satisfies the conclusions of $(TR4)$. 
To prove this, we construct a diagram in $\calC$
$$ \xymatrix{ X \ar[r]^{f} \ar[d] & Y \ar[r]^{g} \ar[d] & Z \ar[r] \ar[d] & 0 \ar[d] & \\
0 \ar[r] & Y/X \ar[d] \ar[r] & Z/X \ar[d] \ar[r] & X' \ar[r] \ar[d] & 0 \ar[d] \\
& 0 \ar[r] & Z/Y \ar[r] & Y' \ar[r] & (Y/X)' }$$
where $0$ is a zero object of $\calC$, and each square in the diagram is a pushout (more precisely, we apply Proposition \toposref{lklk} repeatedly to construct a map from the nerve of the appropriate partially ordered set into $\calC$). Restricting to appropriate rectangles contained in the diagram, we obtain
isomorphisms $X' \simeq X[1]$, $Y' \simeq Y[1]$, $(Y/X)' \simeq Y/X[1]$, and
four distinguished triangles
$$ X \stackrel{f}{\rightarrow} Y \rightarrow Y/X \rightarrow X[1] $$
$$ Y \stackrel{g}{\rightarrow} Z \rightarrow Z/Y \rightarrow Y[1] $$
$$ X \stackrel{g \circ f}{\rightarrow} Z \rightarrow Z/X \rightarrow X[1]$$
$$ Y/X \rightarrow Z/X \rightarrow Z/Y \rightarrow Y/X[1].$$
The commutativity in the homotopy category $\h{\calC}$ required by $(TR4)$ follows from the (stronger) commutativity of the above diagram in $\calC$ itself.
\end{itemize}
\end{proof}

\begin{remark}
The definition of a stable $\infty$-category is quite
a bit simpler than that of a triangulated category. In particular,
the octahedral axiom $(TR4)$ is a consequence of $\infty$-categorical
principles which are basic and easily motivated.
\end{remark}

\begin{notation}\index{ZZZExt@$\Ext^n_{\calC}(X,Y)$}
Let $\calC$ be a stable $\infty$-category containing a pair of objects $X$ and $Y$. We
let $\Ext^n_{\calC}(X,Y)$ denote the abelian group $\Hom_{ \h{\calC}}( X[n], Y)$. 
If $n$ is negative, this can be identified with the homotopy group
$\pi_{-n} \bHom_{\calC}(X,Y)$. More generally, $\Ext^{n}_{\calC}(X,Y)$ can be identified with the $(-n)$th homotopy group of an appropriate {\em spectrum} of maps from $X$ to $Y$.
\end{notation}

\section{Properties of Stable $\infty$-Categories}\label{stable4}

According to Definition \ref{stabl}, a pointed $\infty$-category $\calC$ is stable if it admits certain pushout squares and certain pullback squares, which are required to coincide with one another. Our goal in this section is to prove that a stable $\infty$-category $\calC$ admits {\em all} finite limits and colimits, and that the pushout squares in $\calC$ coincide with the pullback squares in general (Proposition \ref{surose}). To prove this, we will need the following easy observation (which is quite useful in its own right):

\begin{proposition}\label{expfun}
Let $\calC$ be a stable $\infty$-category, and let $K$ be a simplicial set.
Then the $\infty$-category $\Fun(K, \calC)$ is stable.
\end{proposition}

\begin{proof}
This follows immediately from the fact that kernels and cokernels
in $\Fun(K, \calC)$ can be computed pointwise (Proposition \toposref{limiteval}). 
\end{proof}

\begin{definition}\index{stable subcategory}\index{subcategory!stable}
If $\calC$ is stable $\infty$-category, and $\calC_0$ is a
full subcategory containing a zero object and stable under the
formation of kernels and cokernels, then $\calC_0$ is itself stable. In this case, we will say that
$\calC_0$ is a {\it stable subcategory} of $\calC$.
\end{definition}

\begin{lemma}\label{surefoot}
Let $\calC$ be a stable $\infty$-category, and let $\calC' \subseteq \calC$ be 
a full subcategory which is stable under cokernels and under translation. Then
$\calC'$ is a stable subcategory of $\calC$.
\end{lemma}

\begin{proof}
It will suffice to show that $\calC'$ is stable under kernels. Let $f: X \rightarrow Y$ be a morphism in $\calC$. Theorem \ref{surmite} shows that there is a canonical equivalence $\ker(f) \simeq \coker(f)[-1]$. 
\end{proof}

\begin{proposition}\label{surose}
Let $\calC$ be a pointed $\infty$-category. Then $\calC$ is stable if and only if
the following conditions are satisfied:
\begin{itemize}
\item[$(1)$] The $\infty$-category $\calC$ admits finite limits and colimits.
\item[$(2)$] A square
$$ \xymatrix{ X \ar[r] \ar[d] & Y \ar[d] \\
X' \ar[r] & Y' }$$
in $\calC$ is a pushout if and only if it is a pullback.
\end{itemize}
\end{proposition}

\begin{proof}
Condition $(1)$ implies the existence of kernels and cokernels in $\calC$, and
condition $(2)$ implies that the exact triangles coincide with the coexact triangles.
This proves the ``if'' direction. 

Suppose now that $\calC$ is stable. We begin by proving $(1)$. It will suffice to show that $\calC$ admits finite colimits; the dual argument will show that $\calC$ admits finite limits as well.
According to Proposition \toposref{appendixdiagram}, it will suffice to show that $\calC$ admits coequalizers and finite coproducts. The existence of finite coproducts was established in Lemma \ref{vival}. We now conclude by observing that a coequalizer for a 
diagram
$$\xymatrix{ X \ar@<.4ex>[r]^{f} \ar@<-.4ex>[r]_{f'} & Y}$$
can be identified with $\coker( f-f')$. 

We now show that every pushout square in $\calC$ is a pullback; the converse will follow by a dual argument. Let $\calD \subseteq \Fun( \Delta^1 \times \Delta^1, \calC)$ be the full subcategory
spanned by the pullback squares. Then $\calD$ is stable under finite limits and under translations.
It follows from Lemma \ref{surefoot} that $\calD$ is a stable subcategory of
$\Fun( \Delta^1 \times \Delta^1, \calC)$.

Let $i: \Lambda^2_0 \hookrightarrow \Delta^1 \times \Delta^1$ be the inclusion, and let
$i_{!}: \Fun( \Lambda^2_0, \calC) \rightarrow \Fun( \Delta^1 \times \Delta^1, \calC)$ be a functor of left Kan extension. Then $i_{!}$ preserves finite colimits, and is therefore exact (Proposition \ref{funrose}). Let $\calD' = i_{!}^{-1} \calD$. Then
$\calD'$ is a stable subcategory of $\Fun( \Lambda^2_0, \calC)$; we wish to show that 
$\calD' = \Fun(\Lambda^2_0, \calC)$. To prove this, we observe that any diagram
$$ X' \leftarrow X \rightarrow X''$$
can be obtained as a (finite) colimit
$$ e'_{X'} \coprod_{ e'_{X} } e_{X} \coprod_{ e''_{X} } e''_{X''}$$
where $e_{X} \in \Fun( \Lambda^2_0, \calC)$ denotes the diagram
$ X \leftarrow X \rightarrow X,$
$e'_{Z} \in \Fun( \Lambda^2_0, \calC)$ denotes the diagram
$ Z \leftarrow 0 \rightarrow 0,$
and $e''_{Z} \in \Fun( \Lambda^2_0, \calC)$ denotes the diagram
$ 0 \leftarrow 0 \rightarrow Z.$
It will therefore suffice to prove that pushout of any of these five diagrams is also a pullback. This follows immediately from the following more general observation:
any pushout square
$$ \xymatrix{ A \ar[r] \ar[d]^{f} & A' \ar[d] \\
B \ar[r] & B' }$$
in an (arbitrary) $\infty$-category $\calC$ is also pullback square, provided that $f$ is an equivalence.
\end{proof}
 
\begin{proposition}\label{kappstable}
Let $\calC$ be a $($small$)$ stable $\infty$-category, and let $\kappa$ be a regular cardinal. Then the $\infty$-category $\Ind_{\kappa}(\calC)$ is stable.
\end{proposition}

\begin{proof}
The functor $j$ preserves finite limits and colimits (Propositions \toposref{yonedaprop} and \toposref{turnke}). It follows that $j(0)$ is a zero object of $\Ind_{\kappa}(\calC)$, so that $\Ind_{\kappa}(\calC)$ is pointed.

We next show that every morphism $f: X \rightarrow Y$ in $\Ind_{\kappa}(\calC)$ admits a kernel and a cokernel. According to Proposition \toposref{urgh1}, we may assume that $f$ is a $\kappa$-filtered colimit of morphisms $f_{\alpha}: X_{\alpha} \rightarrow Y_{\alpha}$ which belong to the essential image $\calC'$ of $j$. Since $j$ preserves kernels and cokernels, each of the maps $f_{\alpha}$ has a kernel and a cokernel in $\Ind_{\kappa}$. It follows immediately that $f$ has a cokernel (which can be written as a colimit of the cokernels of the maps $f_{\alpha}$). The existence of $\ker(f)$ is slightly more difficult. Choose a $\kappa$-filtered diagram $p: \calI \rightarrow \Fun( \Delta^1 \times \Delta^1, \calC')$, where each $p(\alpha)$ is a pullback square
$$ \xymatrix{ Z_{\alpha} \ar[r] \ar[d] & 0 \ar[d] \\
X_{\alpha} \ar[r]^{f_{\alpha}} & Y_{\alpha}. } $$
Let $\sigma$ be a colimit of the diagram $p$; we wish to show that $\sigma$ is a pullback diagram
in $\Ind_{\kappa}(\calC)$. Since $\Ind_{\kappa}(\calC)$ is stable under $\kappa$-small limits in $\calP(\calC)$, it will suffice to show that $\sigma$ is a pullback square in $\calP(\calC)$.
Since $\calP(\calC)$ is an $\infty$-topos, filtered colimits in $\calP(\calC)$ are left exact (Example \toposref{tucka}); it will therefore suffice to show that each $p(\alpha)$ is a pullback diagram
in $\calP(\calC)$. This is obvious, since the inclusion $\calC' \subseteq \calP(\calC)$ preserves all limits which exist in $\calC'$ (Proposition \toposref{yonedaprop}). 

To complete the proof, we must show that a triangle in $\Ind_{\kappa}(\calC)$ is exact if and only if it is coexact. Suppose given an exact triangle
$$ \xymatrix{ Z \ar[r] \ar[d] & 0 \ar[d] \\
X \ar[r] & Y }$$
in $\Ind_{\kappa}(\calC)$. The above argument shows that we can write this triangle as a filtered colimit of exact triangles
$$ \xymatrix{ Z_{\alpha} \ar[r] \ar[d] & 0 \ar[d] \\
X_{\alpha} \ar[r] & Y_{\alpha} }$$
in $\calC'$. Since $\calC'$ is stable, we conclude that these triangles are also coexact. The original triangle is therefore a filtered colimit of coexact triangles in $\calC'$, hence coexact. The converse follows by the same argument. 
\end{proof}

 \section{Exact Functors}\label{stable5}
 
Let $F: \calC \rightarrow \calC'$ be a functor between stable
$\infty$-categories. Suppose that $F$ carries zero objects into
zero objects. It follows immediately that $F$ carries triangles
into triangles. If, in addition, $F$ carries exact triangles into
exact triangles, then we will say that $F$ is {\it exact}. The exactness of a functor $F$ 
admits the following alternative characterizations:\index{exact!functor}\index{functor!exact}

\begin{proposition}\label{funrose}
Let $F: \calC \rightarrow \calC'$ be a functor between stable
$\infty$-categories. The following conditions are equivalent:
\begin{itemize}
\item[$(1)$] The functor $F$ is left exact. That is, $F$ commutes with
finite limits. 
\item[$(2)$] The functor $F$ is right exact. That is, $F$
commutes with finite colimits. 
\item[$(3)$] The functor $F$ is exact.
\end{itemize}
\end{proposition}

\begin{proof}
We will prove that $(2) \Leftrightarrow (3)$; the equivalence $(1) \Leftrightarrow (3)$ will follow by a dual argument. The implication $(2) \Rightarrow (3)$ is obvious. Conversely, suppose that $F$ is exact. The proof of Proposition \ref{surose} shows that $F$ preserves coequalizers, and the proof of
Lemma \ref{vival} shows that $F$ preserves finite coproducts. It follows that $F$ preserves
all finite colimits (see the proof of Proposition \toposref{appendixdiagram}). 
\end{proof}

The identity functor from any stable $\infty$-category to 
itself\index{ZZZCatiex@$\Cat_{\infty}^{\Ex}$} is exact, and a composition of exact functors is exact.
Consequently, there exists a subcategory $\Cat_{\infty}^{\Ex} \subseteq \Cat_{\infty}$ in which the objects are stable $\infty$-categories and the morphisms are the exact functors.
Our next few results concern the stability properties of this subcategory.

\begin{proposition}\label{hung1}
Suppose given a homotopy Cartesian diagram of $\infty$-categories
$$ \xymatrix{ \calC' \ar[r]^{G'} \ar[d]^{F'} & \calC \ar[d]^{F} \\
\calD' \ar[r]^{G} & \calD. }$$
Suppose further that $\calC$, $\calD'$, and $\calD$ are stable, and that the functors
$F$ and $G$ are exact. Then:
\begin{itemize}
\item[$(1)$] The $\infty$-category $\calC'$ is stable.
\item[$(2)$] The functors $F'$ and $G'$ are exact.
\item[$(3)$] If $\calE$ is a stable $\infty$-category, then a functor
$H: \calE \rightarrow \calC'$ is exact if and only if the functors
$F' \circ H$ and $G' \circ H$ are exact.
\end{itemize}
\end{proposition}

\begin{proof}
Combine Proposition \ref{surose} with Lemma \toposref{bird3}.
\end{proof}

\begin{proposition}\label{hung2}
Let $\{ \calC_{\alpha} \}_{ \alpha \in A }$ be a collection of stable $\infty$-categories. Then
the product $$\calC = \prod_{\alpha \in A} \calC_{\alpha} $$ is stable.
Moreover, for any stable $\infty$-category $\calD$, a functor
$F: \calD \rightarrow \calC$ is exact if and only if each of the compositions
$$ \calD \stackrel{F}{\rightarrow} \calC \rightarrow \calC_{\alpha}$$
is an exact functor.
\end{proposition}

\begin{proof}
This follows immediately from the fact that limits and colimits in $\calC$ are computed pointwise.
\end{proof}

\begin{theorem}\label{sunty}
The $\infty$-category $\Cat_{\infty}^{\Ex}$ admits small limits, and the inclusion
$$ \Cat_{\infty}^{\Ex} \subseteq \Cat_{\infty}$$
preserves small limits.
\end{theorem}

\begin{proof}
Using Propositions \ref{hung1} and \ref{hung2}, one can repeat the argument used to prove
Proposition \toposref{accprop}.
\end{proof}

We now prove an analogue of Theorem \ref{sunty}.

\begin{proposition}\label{puterstor}
Let $p: X \rightarrow S$ be an inner fibration of simplicial sets. Suppose
that:
\begin{itemize}
\item[$(i)$] For each vertex $s$ of $S$, the fiber $X_{s} = X \times_{S} \{s\}$ is a stable $\infty$-category.
\item[$(ii)$] For every edge $s \rightarrow s'$ in $S$, the restriction
$X \times_{S} \Delta^1 \rightarrow \Delta^1$ is a coCartesian fibration, associated
to an exact functor $X_{s} \rightarrow X_{s'}$.
\end{itemize}
Then:
\begin{itemize}
\item[$(1)$] The $\infty$-category $\bHom_{S}(S,X)$ of sections of $p$ is stable.
\item[$(2)$] If $\calC$ is an arbitrary stable $\infty$-category, and $f: \calC \rightarrow \bHom_{S}(S,X)$ induces an exact functor $\calC \stackrel{f}{\rightarrow} \bHom_{S}(S,X) \rightarrow X_{s}$ for every vertex $s$ of $S$, then $f$ is exact.
\item[$(3)$] For every set $\calE$ of edges of $S$, let $Y(\calE) \subseteq \bHom_{S}(S,X)$ be the full subcategory spanned by those sections $f: S \rightarrow X$ of $p$ with the following property:
\begin{itemize}
\item[$(\ast)$] For every $e \in \calE$, $f$ carries $e$ to a $p_{e}$-coCartesian edge of
the fiber product $X \times_{S} \Delta^1$, where $p_{e}: X \times_{S} \Delta^1 \rightarrow \Delta^1$ denotes the projection.
\end{itemize}
Then each $Y(\calE)$ is a stable subcategory of $\bHom_{S}(S,X)$.
\end{itemize}
\end{proposition}

\begin{proof}
Combine Proposition \toposref{prestorkus}, Theorem \ref{sunty}, and Proposition \ref{expfun}.
\end{proof}

\begin{proposition}
The $\infty$-category $\Cat_{\infty}^{\Ex}$ admits small filtered colimits, and the inclusion
$\Cat_{\infty}^{\Ex} \subseteq \Cat_{\infty}$ preserves filtered colimits.
\end{proposition}

\begin{proof}
Let $\calI$ be a filtered $\infty$-category, $p: \calI \rightarrow \Cat_{\infty}^{\Ex}$ a diagram, which we will indicate by $\{ \calC_{I} \}_{I \in \calI}$, and $\calC$ a colimit of the induced diagram $\calI \rightarrow \Cat_{\infty}$. We must prove:
\begin{itemize}
\item[$(i)$] The $\infty$-category $\calC$ is stable.
\item[$(ii)$] Each of the canonical functors $\theta_I: \calC_I \rightarrow \calC$ is exact.
\item[$(iii)$] Given an arbitrary stable $\infty$-category $\calD$, a functor
$f: \calC \rightarrow \calD$ is exact if and only if each of the composite functors
$\calC_I \stackrel{\theta_I}{\rightarrow} \calC \rightarrow \calD$ is exact.
\end{itemize}
In view of Proposition \ref{funrose}, $(ii)$ and $(iii)$ follow immediately from Proposition \toposref{unrose}. The same result implies that $\calC$ admits finite limits and colimits, and that each of the functors $\theta_I$ preserves finite limits and colimits.

To prove that $\calC$ has a zero object, we select an object $I \in \calI$. The functor
$\calI \rightarrow \calC$ preserves initial and final objects. Since $\calC_I$ has a zero object, so does
$\calC$.

We will complete the proof by showing that every exact triangle in $\calC$ is coexact (the converse follows by the same argument). Fix a morphism $f: X \rightarrow Y$ in $\calC$. Without loss of generality, we may suppose that there exists $I \in \calI$ and a morphism $\widetilde{f}: \widetilde{X} \rightarrow \widetilde{Y}$ in
$\calC_I$ such that $f = \theta_I( \widetilde{f} )$ (Proposition \toposref{grapeape}). Form a pullback diagram $\widetilde{\sigma}$
$$ \xymatrix{ \widetilde{W} \ar[r] \ar[d] & \widetilde{X} \ar[d] \\
0 \ar[r] & \widetilde{Y} }$$
in $\calC_I$. Since $\calC_I$ is stable, this diagram is also a pushout. It follows that 
$\theta_I(\sigma)$ is triangle $W \rightarrow X \stackrel{f}{\rightarrow} Y$ which is both exact and coexact in $\calC$. 
\end{proof}
 
\section{t-Structures and Localizations}\label{stable6}

Let $\calC$ be an $\infty$-category. Recall that we say a full subcategory $\calC' \subseteq \calC$
is a {\it localization} of $\calC$ if the inclusion functor $\calC' \subseteq \calC$ has a left adjoint (\S \toposref{locfunc}). In this section, we will introduce a special class of localizations, called {\it $t$-localizations}, in the case where $\calC$ is stable. We will further show that there is a bijective correspondence between t-localizations of $\calC$ and {\it t-structures} on the triangulated category $\h{\calC}$. We begin with a review of the classical theory of t-structures; for a more thorough introduction we refer the reader to \cite{BBD}.\index{localization!of a stable $\infty$-category}\index{t-structure}

\begin{definition}\label{deft}\index{ZZZCgz@$\calC_{\geq n}$}\index{ZZZClz@$\calC_{\leq n}$}
Let $\calD$ be a triangulated category. A {\it
t-structure} on $\calD$ is defined to be a pair of full
subcategories $\calD_{\geq 0}$, $\calD_{\leq 0}$ (always assumed
to be stable under isomorphism) having the following properties:

\begin{itemize}
\item[$(1)$] For $X \in \calD_{\geq 0}$ and $Y \in \calD_{\leq -1}$,
we have $\Hom_{\calD}(X,Y) = 0$.

\item[$(2)$] $\calD_{\geq 0}[1]
\subseteq \calD_{\geq 0}$, $\calD_{\leq 0}[-1] \subseteq
\calD_{\leq 0}$. 

\item[$(3)$] For any $X \in \calD$, there exists a
distinguished triangle $ X' \rightarrow X \rightarrow X''
\rightarrow X'[1]$ where $X' \in \calD_{\geq 0}$ and $X'' \in
\calD_{\leq 0}[-1]$.
\end{itemize}
\end{definition}

\begin{notation}
If $\calD$ is a triangulated category equipped with a $t$-structure, we will write
$\calD_{\geq n}$ for $\calD_{\geq 0}[n]$ and $\calD_{\leq n}$ for
$\calD_{\leq 0}[n]$. Observe that we use a {\em homological} indexing convention.
\end{notation}

\begin{remark}\label{quatthrust}
In Definition \ref{deft}, either of the full subcategories $\calD_{\geq 0}, \calD_{\leq 0} \subseteq \calC$ determines the other. For example, an object $X \in \calD$ belongs to $\calD_{\leq -1}$ if
and only if $\Hom_{\calD}(Y,X)$ vanishes for all $Y \in \calD_{\geq 0}$. 
\end{remark}

\begin{definition}
Let $\calC$ be a stable $\infty$-category. A {\it t-structure} on $\calC$ is a 
t-structure on the homotopy category $\h{\calC}$. If $\calC$ is equipped with a t-structure, we let
$\calC_{ \geq n}$ and $\calC_{ \leq n}$ denote the full subcategories of $\calC$ spanned by those objects which belong to $(\h{\calC})_{\geq n}$ and $(\h{\calC})_{\leq n}$, respectively.
\end{definition}

\begin{proposition}\label{surteen}
Let $\calC$ be a stable $\infty$-category equipped with a t-structure. For each
$n \in \Z$, the full subcategory $\calC_{\leq n}$ is a localization of $\calC$.
\end{proposition}

\begin{proof}
Without loss of generality, we may suppose $n = -1$. According to Proposition \toposref{testreflect}, it will suffice to prove that for each $X \in \calC$, there exists a map
$f: X \rightarrow X''$, where $X'' \in \calC_{\leq -1}$ and for each $Y \in \calC_{\leq -1}$, the map
$$ \bHom_{\calC}(X'', Y) \rightarrow \bHom_{\calC}(X,Y)$$
is a weak homotopy equivalence. Invoking part $(3)$ of Definition \ref{deft}, we can choose
$f$ to fit into a distinguished triangle
$$ X' \rightarrow X \stackrel{f}{\rightarrow} X'' \rightarrow X'[1]$$
where $X' \in \calC_{\geq 0}$. According to Whitehead's theorem, we need to show that for every $k \leq 0$, the map $$ \Ext^{k}_{\calC}(X'',Y) \rightarrow \Ext^{k}_{\calC}(X,Y)$$ is an isomorphism of abelian groups.
Using the long exact sequence associated to the exact triangle above, we are reduced to proving that the groups $\Ext^{k}_{\calC}(X',Y)$ vanish for $k \leq 0$.
We now use condition $(2)$ of Definition \ref{deft} to conclude that $X'[-k] \in \calC_{\geq 0}$. Condition $(1)$ of Definition \ref{deft} now implies that
$$\Ext^{k}_{\calC}(X',Y) \simeq \Hom_{ \h{\calC}}( X'[-k], Y) \simeq 0.$$
\end{proof}

\begin{corollary}
Let $\calC$ be a stable $\infty$-category equipped with a t-structure. The full subcategories
$ \calC_{ \leq n} \subseteq \calC$ are stable under all limits which exist in $\calC$. Dually, the full subcategories
$\calC_{ \geq 0} \subseteq \calC$ are stable under all colimits which exist in $\calC$.
\end{corollary}

\begin{notation}\label{sumogirl}\index{truncation functor}\index{ZZZtauk@$\tau_{\leq k}$}\index{ZZZtaugk@$\tau_{\geq k}$}
Let $\calC$ be a stable $\infty$-category equipped with a t-structure. We will let
$\tau_{\leq n}$ denote a left adjoint to the inclusion $\calC_{ \leq n} \subseteq \calC$, and
$\tau_{\geq n}$ a right adjoint to the inclusion $\calC_{ \geq n} \subseteq \calC$.
\end{notation}

\begin{remark}\label{sumochat}
Fix $n, m \in \Z$, and let $\calC$ be a stable $\infty$-category equipped with a t-structure. Then
the truncation functors $\tau_{\leq n}$, $\tau_{\geq n}$ map the full subcategory $\calC_{\leq m}$ to itself. To prove this, we first observe that $\tau_{ \leq n}$ is equivalent to the identity on
$\calC_{\leq m}$ if $m \leq n$, while if $m \geq n$ the essential image of $\tau_{\leq n}$ is contained in $\calC_{\leq n} \subseteq \calC_{\leq m}$. To prove the analogous result for $\tau_{\geq n}$, we observe that the proof of Proposition \ref{surteen} implies that for
each $X$, we have a distinguished triangle
$$ \tau_{\geq n} X \rightarrow X \stackrel{f}{\rightarrow} \tau_{\leq n-1} X \rightarrow (\tau_{\geq n} X)[1].$$
If $X \in \calC_{\leq m}$, then $\tau_{\leq n-1} X$ also belongs to $\calC_{\leq m}$, so that
$\tau_{ \geq n} X \simeq \ker(f)$ belongs to $\calC_{\leq m}$ since $\calC_{\leq m}$ is stable under limits.
\end{remark}

\begin{warning}\label{umper}
In \S \toposref{truncintro}, we introduced for every $\infty$-category $\calC$ a full subcategory
$\tau_{\leq n} \calC$ of {\it $n$-truncated objects} of $\calC$. In that context, we used the symbol $\tau_{\leq n}$ to denote a left adjoint to the inclusion $\tau_{\leq n} \calC \subseteq \calC$.
This is {\em not} compatible with Notation \ref{sumogirl}. In fact, if $\calC$ is a stable $\infty$-category, then it has no nonzero truncated objects at all: if $X \in \calC$ is nonzero, then
the identity map from $X$ to itself determines a nontrivial homotopy class in 
$\pi_{n} \bHom_{\calC}( X[-n], X)$, for all $n \geq 0$. Nevertheless, the two notations are
consistent when restricted to $\calC_{\geq 0}$, in view of the following fact:
\begin{itemize}
\item Let $\calC$ be a stable $\infty$-category equipped with a t-structure. An object $X \in \calC_{\geq 0}$ is $k$-truncated (as an object of $\calC_{\geq 0}$) if and only if $X \in \calC_{\leq k}$. 
\end{itemize}
In fact, we have the following more general statement: for any $X \in \calC$ and $k \geq -1$, $X$ belongs to $\calC_{\leq k}$ if and only if $\bHom_{\calC}(Y,X)$ is $k$-truncated for every $Y \in \calC_{\geq 0}$. Because the latter condition is equivalent to the vanishing of $\Ext^{n}_{\calC}(Y,X)$ for $n < -k$, we can use the shift functor to reduce to the case where $n=0$ and $k=-1$, which is covered by Remark \ref{quatthrust}.
\end{warning}

Let $\calC$ be a stable $\infty$-category equipped with a t-structure, and let $n, m \in \Z$. Remark \ref{sumochat} implies that we have a commutative diagram of simplicial sets
$$ \xymatrix{ \calC_{ \geq n} \ar@{^{(}->}[r] \ar[d]^{\tau_{\leq m}} & \calC \ar[d]^{\tau_{\leq m} } \\
\calC_{\geq n} \cap \calC_{\leq m} \ar@{^{(}->}[r] & \calC_{\leq m}. }$$
As explained in \S \toposref{propertopoi}, we get an induced transformation
of functors $$ \theta: \tau_{\leq m} \circ \tau_{ \geq n} \rightarrow \tau_{\geq n} \circ \tau_{\leq m}.$$

\begin{proposition}
Let $\calC$ be a stable $\infty$-category equipped with t-structure. Then the natural transformation
$$ \theta: \tau_{\leq m} \circ \tau_{ \geq n} \rightarrow \tau_{\geq n} \circ \tau_{\leq m}$$
is an equivalence of functors $\calC \rightarrow \calC_{\leq m} \cap \calC_{\geq n}$. 
\end{proposition}

\begin{proof}
This is a classical fact concerning triangulated categories; we include a proof for completeness.
Fix $X \in \calC$; we wish to show that
$$\theta(X): \tau_{\leq m} \tau_{\geq n} X \rightarrow \tau_{\geq n} \tau_{\leq m} X$$
is an isomorphism in the homotopy category of $\calC_{\leq m} \cap \calC_{\geq n}$. If $m < n$, then both sides are
zero and there is nothing to prove; let us therefore assume that $m \geq n$.
Fix $Y \in \calC_{\leq m} \cap \calC_{\geq n}$; it will suffice to show that composition with $\theta(X)$ induces an isomorphism
$$ \Ext^0( \tau_{\geq n} \tau_{\leq m} X, Y) \rightarrow \Ext^0( \tau_{\leq m} \tau_{\geq n} X, Y)
\simeq \Ext^0( \tau_{\geq n} X,Y).$$
We have a map of long exact sequences
$$ \xymatrix{ \Ext^0( \tau_{\leq n-1} \tau_{\leq m} X, Y) \ar[d] \ar[r]^{f_0} & \Ext^0( \tau_{\leq n-1} X, Y) \ar[d] \\
\Ext^0( \tau_{\leq m} X,Y) \ar[r] \ar[d] \ar[r]^{f_1} \ar[d] & \Ext^0( X, Y) \ar[d] \\
\Ext^0( \tau_{\geq n} \tau_{\leq m} X,Y) \ar[r]^{f_2} \ar[d] & \Ext^0( \tau_{\geq n} X, Y) \ar[d] \\
\Ext^1( \tau_{\leq n-1} \tau_{\leq m} X,Y) \ar[d] \ar[r]^{f_3} & \Ext^1( \tau_{\leq n-1} X, Y) \ar[d] \\
\Ext^1( \tau_{\leq m} X, Y) \ar[r]^{f_4} & \Ext^1(X,Y). }$$
Since $m \geq n$, the natural transformation $\tau_{ \leq n-1} \rightarrow \tau_{ \leq n-1} \tau_{ \leq m}$ is an equivalence of functors; this proves that $f_0$ and $f_3$ are bijective.
Since $Y \in \calC_{\leq m}$, $f_1$ is bijective and $f_4$ is injective. It follows from the
``five lemma'' that $f_2$ is bijective, as desired.
\end{proof}

\begin{definition}\index{heart}\index{ZZZCheart@$\heart{\calC}$}
Let $\calC$ be a stable $\infty$-category equipped with a $t$-structure. The {\it heart}
$\heart{\calC}$ of $\calC$ is the full subcategory $\calC_{\geq 0} \cap \calC_{\leq 0} \subseteq \calC$. For each $n \in \Z$, we let $\pi_{0}: \calC \rightarrow \heart{\calC}$ denote the
functor $\tau_{\leq 0} \circ \tau_{\geq 0} \simeq \tau_{\geq 0} \circ \tau_{\leq 0}$, and
we let $\pi_{n}: \calC \rightarrow \heart{\calC}$ denote the composition of $\pi_0$ with the shift functor $X \mapsto X[-n]$. 
\end{definition}

\begin{remark}
Let $\calC$ be a stable $\infty$-category equipped with a t-structure, and let
$X,Y \in \heart{\calC}$. The homotopy group $\pi_n \bHom_{\calC}(X,Y) \simeq \Ext^{-n}_{\calC}(X,Y)$ vanishes for $n > 0$. It follows that $\heart{\calC}$ is equivalent to (the nerve of) its homotopy category $\h{\heart{\calC}}$. Moreover, the category $\h{\heart{\calC}}$ is abelian (\cite{BBD}). \end{remark}

Let $\calC$ be a stable $\infty$-category. In view of Remark \ref{quatthrust}, t-structures
on $\calC$ are determined by the corresponding localizations $\calC_{\leq 0} \subseteq \calC$. 
However, not every localization of $\calC$ arises in this way.
Recall (see \S \toposref{invloc}) that every localization of $\calC$ has the form $S^{-1} \calC$, where
$S$ is an appropriate collection of morphisms of $\calC$. Here
$S^{-1} \calC$ denotes the full subcategory of $\calC$ spanned by $S$-local objects, where
an object $X \in \calC$ is said to be $S$-local if and only if, for each $f: Y' \rightarrow Y$
in $S$, composition with $f$ induces a homotopy equivalence
$$ \bHom_{\calC}(Y,X) \rightarrow \bHom_{\calC}(Y',X).$$
If $\calC$ is stable, then we extend the morphism $f$ to a distinguished triangle
$$ Y' \rightarrow Y \rightarrow Y'' \rightarrow Y'[1],$$
and we have an associated long exact sequence
$$ \ldots \rightarrow \Ext^i_{\calC}(Y'', X) \rightarrow \Ext^i_{\calC}(Y, X)
\stackrel{ \theta_{i}}{\rightarrow} \Ext^i_{\calC}(Y',X)
\rightarrow \Ext^{i+1}_{\calC}( Y'',X) \rightarrow \ldots$$
The requirement that $X$ be $\{ f \}$-local amounts to the condition that $\theta_i$ be
an isomorphism for $i \leq 0$. Using the long exact sequence, we see that if
$X$ is $\{ f \}$-local, then $\Ext^i_{\calC}(Y'',X) = 0$ for $i \leq 0$. Conversely, if
$\Ext^i_{\calC}(Y'',X) = 0$ for $i \leq 1$, then $X$ is $\{ f \}$-local.
Experience suggests that it is usually more natural to require the vanishing of the groups $\Ext^i_{\calC}(Y'',X)$ than it is to require that the maps $\theta_i$ to be isomorphisms.
Of course, if $Y'$ is a zero object of $\calC$, then the distinction between these conditions disappears.

\begin{definition}\label{surtut}
Let $\calC$ be an $\infty$-category which admits pushouts. We will say that a collection $S$ of 
morphisms of $\calC$ is {\it quasisaturated} if it satisfies the following conditions:
\begin{itemize}
\item[$(1)$] Every equivalence in $\calC$ belongs to $S$.
\item[$(2)$] Given a $2$-simplex $\Delta^2 \rightarrow \calC$
$$ \xymatrix{ X \ar[rr]^{h} \ar[dr]^{f} & & Z \\
& Y, \ar[ur]^{g} & }$$
if any two of $f$, $g$, and $h$ belongs to $S$, then so does the third.
\item[$(3)$] Given a pushout diagram
$$ \xymatrix{ X \ar[d]^{f} \ar[r] & X' \ar[d]^{f'} \\
Y \ar[r] & Y', }$$
if $f \in S$, then $f' \in S$.
\end{itemize}
Any intersection of quasisaturated collections of morphisms is weakly saturated. Consequently, for any collection of morphisms $S$ there is a smallest quasisaturated collection $\overline{S}$ containing $S$. We will say that $\overline{S}$ is the {\it quasisaturated collection of morphisms generated by $S$}.
\end{definition}

\begin{definition}\index{closed!under extensions}
Let $\calC$ be a stable $\infty$-category. A full subcategory $\calC' \subseteq \calC$
is {\it closed under extensions} if, for every distinguished triangle
$$ X \rightarrow Y \rightarrow Z \rightarrow X[1]$$
such that $X$ and $Z$ belong to $\calC'$, the object $Y$ also belongs to $\calC'$.
\end{definition}

We observe that if $\calC$ is as in Definition \ref{surtut} and $L: \calC \rightarrow \calC$ is a localization functor, then the collection of all morphisms $f$ of $\calC$ such that $L(f)$ is an equivalence is quasisaturated.

\begin{proposition}\label{condit}
Let $\calC$ be a stable $\infty$-category, let $L: \calC \rightarrow \calC$ be a localization functor, and let $S$ be the collection of morphisms $f$ in $\calC$ such that $L(f)$ is an equivalence. The following conditions are equivalent:

\begin{itemize}
\item[$(1)$] There exists a collection of morphisms $\{ f: 0 \rightarrow X \}$ which
generates $S$ $($as a quasisaturated collection of morphisms$)$.

\item[$(2)$] The collection of morphisms $\{ 0 \rightarrow X : L(X) \simeq 0 \}$ generates
$S$ $($as a quasisaturated collection of morphisms$)$.

\item[$(3)$] The essential image of $L$ is closed under extensions.

\item[$(4)$] For any $A \in \calC$, $B \in L\calC$, the natural map
$\Ext^1(LA,B) \rightarrow \Ext^1(A,B)$ is injective.

\item[$(5)$] The full subcategories $\calC_{\geq 0} = \{ A: LA \simeq 0
\}$ and $\calC_{\leq -1} = \{A: LA \simeq A \}$ determine a
t-structure on $\calC$.

\end{itemize}
\end{proposition}

\begin{proof}
The implication $(1) \Rightarrow (2)$ is obvious. We next prove that $(2) \Rightarrow (3)$.
Suppose given an exact triangle
$$ X \rightarrow Y \rightarrow Z$$
where $X$ and $Z$ are both $S$-local. We wish to prove that $Y$ is $S$-local. In view of assumption $(2)$, it will suffice to show that $\bHom_{\calC}( A, Y)$ is contractible, provided
that $L(A) \simeq 0$. In other words, we must show that $\Ext^i_{\calC}(A,Y) \simeq 0$ for
$i \leq 0$. We now observe that there is an exact sequence
$$ \Ext^i_{\calC}(A,X) \rightarrow \Ext^i_{\calC}(A,Y) \rightarrow \Ext^i_{\calC}(A,Z)$$
where the outer groups vanish, since $X$ and $Z$ are $S$-local and the map
$0 \rightarrow A$ belongs to $S$.

We next show that $(3) \Rightarrow (4)$. 
Let $B \in L \calC$, and let $\eta \in \Ext^1_{\calC}(LA,B)$ classify a
distinguished triangle
$$B \rightarrow C \stackrel{g}{\rightarrow} LA \stackrel{\eta}{\rightarrow} B[1].$$
Condition $(3)$ implies that $C \in L \calC$. If the image of $\eta$ in
$\Ext^1_{\calC}(A,B)$ is trivial, then the localization map $A \rightarrow LA$ factors as a composition
$$ A \stackrel{f}{\rightarrow} C \stackrel{g}{\rightarrow} LA.$$
Applying $L$ to this diagram (and using the fact that $C$ is local) we conclude that the
map $g$ admits a section, so that $\eta = 0$.

We now claim that $(4) \Rightarrow (5)$. Assume $(4)$, and define
$\calC_{\geq 0}$, $\calC_{\leq -1}$ as in $(5)$. We will show that the axioms of
Definition \ref{deft} are satisfied:
\begin{itemize}
\item If $X \in \calC_{ \geq 0}$ and $Y \in \calC_{\leq -1}$, then
$\Ext^0_{\calC}(X,Y) \simeq \Ext^{0}_{\calC}(LX, Y) \simeq \Ext^0_{\calC}(0,Y) \simeq 0$.

\item Since $\calC_{\leq -1}$ is a localization of $\calC$, it is stable under limits, so that
$\calC_{\leq -1}[-1] \subseteq \calC_{\leq -1}$. Similarly, since the functor $L: \calC \rightarrow \calC_{\leq -1}$ preserves all colimits which exist in $\calC$, the subcategory $\calC_{ \geq 0}$ is stable under finite colimits, so that $\calC_{\geq 0} [1] \subseteq \calC_{\geq 0}$.

\item Let $X \in \calC$, and form a distinguished triangle
$$ X' \rightarrow X \rightarrow LX \rightarrow X'[1].$$
We claim that $X' \in \calC_{\geq 0}$; in other words, that $LX' = 0$. 
For this, it suffices to show that for all $Y \in L \calC$, the morphism space
$$ \Ext^0_{\calC}(LX', Y) = 0.$$ Since $Y$ is local, we have isomorphisms
$$ \Ext^0_{\calC}(LX',Y) \simeq \Ext^0_{\calC}(X',Y) \simeq \Ext^1_{\calC}(X'[1],Y).$$ 
We now observe that there is a long exact sequence
$$ \Ext^0(LX,Y) \stackrel{f}{\rightarrow} \Ext^0(X,Y) \rightarrow \Ext^1_{\calC}(X'[1], Y) \rightarrow \Ext^1_{\calC}(LX,Y) \stackrel{f'}{\rightarrow} \Ext^1_{\calC}(X,Y).$$
Here $f$ is bijective (since $Y$ is local) and $f'$ is injective (in virtue of assumption $(4)$).
\end{itemize}

We conclude by showing that $(5) \Rightarrow (1)$. Let $S'$ be the smallest quasisaturated collection of morphisms which contains the zero map $0 \rightarrow A$, for every $A \in \calC_{\geq 0}$. We wish to prove that $S = S'$. For this, we choose an arbitrary morphism $u: X \rightarrow Y$ belonging to $S$. Then $Lu: LX \rightarrow LY$ is an equivalence, so we have a pushout diagram
$$ \xymatrix{ X' \ar[r]^{u'} \ar[d] & Y' \ar[d] \\
X \ar[r]^{u} & Y, }$$
where $X'$ and $Y'$ are kernels of the respective localization maps $X \rightarrow LX$, $Y \rightarrow LY$. Consequently, it will suffice to prove that $u' \in S'$. Since $X', Y' \in \calC_{\geq 0}$, this follows from the two-out-of-three property, applied to the diagram
$$ \xymatrix{ & X' \ar[dr]^{u'} & \\
0 \ar[rr] \ar[ur] & & Y'. }$$
\end{proof}

\begin{proposition}\label{charprojjj}\index{module!projective}
Let $\calC$ be a stable $\infty$-category equipped with a left-complete t-structure. Let $P \in \calC_{\geq 0}$. The following conditions are equivalent:

\begin{itemize}
\item[$(1)$] The object $P$ is projective in $\calC_{\geq 0}$.
\item[$(2)$] For every $Q \in \calC_{ \leq 0}$, the abelian group
$\Ext^{1}_{\calC}(P,Q)$ vanishes.
\item[$(3)$] Given a distinguished triangle
$$ N' \rightarrow N \rightarrow N'' \rightarrow N'[1],$$
where $N', N, N'' \in \calC_{\geq 0}$, the induced map $\Ext^0_{\calC}( P, N) \rightarrow \Ext^0_{\calC}(P, N'')$
is surjective.
\end{itemize}
\end{proposition}

\begin{proof}
It follows from Lemma \ref{surput} that $\calC_{\geq 0}$ admits geometric realizations for simplicial objects, so that condition $(1)$ makes sense. We first show that $(1) \Rightarrow (2)$. Let $f: \calC \rightarrow \SSet$ be the functor corepresented by $P$. Let $M_{\bigdot}$ be a \Cech nerve for the morphism
$0 \rightarrow Q[1]$, so that $M_{n} \simeq Q^{n} \in \calC_{\geq 0}$. 
Then $Q[1]$ can be identified with the geometric realization
$|M_{\bigdot}|$. Since $P$ is projective, $f(Q[1])$ is equivalent to the geometric realization
$| f(M_{\bigdot}) |$. We have a surjective map $\ast \simeq \pi_0 f(M_0) \rightarrow \pi_0 | f(M_{\bigdot}) |$, so that $\pi_0 f(Q[1]) = \Ext^1_{\calC}(P,Q) = 0$.

We now show that $(2) \Rightarrow (1)$.  
Proposition \ref{urtusk21} implies that $f$ is homotopic to a composition
$$ \calC \stackrel{F}{\rightarrow} \Spectra \stackrel{\Omega^{\infty}}{\rightarrow} \SSet,$$
where $F$ is an exact functor. Applying $(2)$, we deduce that $F$
is right t-exact (Definition \ref{hurtman}). Lemma \ref{surput} implies
that the induced map $\calC_{\geq 0} \rightarrow \connSpectra$ preserves
geometric realizations of simplicial objects. Applying Proposition \ref{denkmal}, we conclude that $f | \calC_{\geq 0}$ preserves geometric realizations as well.

The implication $(2) \Rightarrow (3)$ follows immediately from the exactness
of the sequence
$$ \Ext^0_{\calC}(P,N) \rightarrow \Ext^0_{\calC}(P,N'') \rightarrow \Ext^1_{\calC}(P, N').$$
Conversely, suppose that $(3)$ is satisfied, and let $\eta \in \Ext^1_{\calC}(P,Q)$. 
Then $\eta$ classifies a distinguished triangle
$$ Q \rightarrow Q' \stackrel{g}{\rightarrow} P \rightarrow Q[1].$$
Since $Q, P \in \calC_{\geq 0}$, we have $Q' \in \calC_{\geq 0}$ as well. Invoking
$(3)$, we deduce that $g$ admits a section, so that $\eta = 0$.
\end{proof}

\section{Boundedness and Completeness}\label{stable6.5}

Let $\calC$ be a stable $\infty$-category equipped with a $t$-structure.
We let $\calC^{\dminus} = \bigcup \calC_{ \leq n} \subseteq \calC$, $\calC^{\dplus} = \bigcup \calC_{ \geq -n}$, and $\calC^{b} = \calC^{\dminus} \cap \calC^{\dplus}$. It is easy to see that $\calC^{\dplus}$, $\calC^{\dminus}$, and $\calC^{b}$ are stable subcategories of $\calC$. We will say that $\calC$ is {\it left bounded}
if $\calC = \calC^{\dminus}$, {\it right bounded} if $\calC = \calC^{\dplus}$, and {\it bounded} if $\calC = \calC^{b}$.\index{ZZZCplus@$\calC^{\dminus}$}\index{ZZZCminus@$\calC^{\dplus}$}\index{ZZZCbound@$\calC^{b}$}\index{left bounded}\index{right bounded}\index{bounded t-structure}\index{t-structure!left bounded}\index{t-structure!right bounded}\index{t-structure!bounded}

At the other extreme, given a stable $\infty$-category $\calC$\index{ZZZCcomp@$\widehat{\calC}$}\index{completion!left}\index{completion!right}\index{t-structure!left complete}\index{t-structure!right complete} equipped with a t-structure, we define the {\it left completion}
$\widehat{\calC}$ of $\calC$ to be homotopy limit of the tower
$$\ldots \rightarrow \calC_{\leq 2}
\stackrel{\tau_{\leq 1}}{\rightarrow} \calC_{\leq 1}
\stackrel{\tau_{\leq 0}}{\rightarrow} \calC_{\leq 0} \stackrel{\tau_{\leq -1}}{\rightarrow} \ldots$$ 
Using the results of \S \toposref{catlim}, we can obtain a very concrete description
of this inverse limit: it is the full subcategory of $\Fun( \Nerve(\Z), \calC)$ spanned by
those functors $F: \Nerve(\Z) \rightarrow \calC$ with the following properties:
\begin{itemize}
\item[$(1)$] For each $n \in \Z$, $F(n) \in \calC_{ \leq -n}$.
\item[$(2)$] For each $m \leq n \in \Z$, the associated map
$F(m) \rightarrow F(n)$ induces an equivalence
$\tau_{ \leq -n} F(m) \rightarrow F(n)$.
\end{itemize}
We will denote this inverse limit by $\widehat{\calC}$, and refer to it as
the {\it left completion} of $\calC$.

\begin{proposition}\label{sparrow}
Let $\calC$ be a stable $\infty$-category equipped with a t-structure. Then:
\begin{itemize}
\item[$(1)$] The left completion $\widehat{\calC}$ is also stable.
\item[$(2)$] Let $\widehat{\calC}_{ \leq 0 }$ and $\widehat{\calC}_{\geq 0}$
be the full subcategories of $\widehat{\calC}$ spanned by those functors
$F: \Nerve(\Z) \rightarrow \calC$ which factor through $\calC_{\leq 0}$ and
$\calC_{\geq 0}$, respectively. Then these subcategories determine a t-structure on $\widehat{\calC}$.
\item[$(3)$] There is a canonical functor $\calC \rightarrow \widehat{\calC}$. This
functor is exact, and induces an equivalence $\calC_{\leq 0} \rightarrow \widehat{\calC}_{\leq 0}$.
\end{itemize}
\end{proposition}

\begin{proof}
We observe that $\widehat{\calC}$ can be identified with the homotopy inverse limit of the tower
$$ \ldots \rightarrow \calC_{\leq 0} \stackrel{ \tau_{\leq 0} \Sigma}{\rightarrow} \calC_{\leq 0} 
\stackrel{\tau_{\leq 0} \Sigma}{\rightarrow} \calC_{\leq 0}.$$
In other words, $\widehat{\calC}^{op} \simeq \Spectra(\calC^{op})$ (see \S \ref{stable10}). 
Assertion $(1)$ now follows from Proposition \ref{sinkle}.

We next prove $(2)$. We begin by observing that, 
if we identify $\widehat{\calC}$ with a full subcategory of
$\Fun(\Nerve(\Z), \calC)$, then the shift functors on $\widehat{\calC}$ can be defined by the formula $$ (F[n])(m) = F(m-n)[n].$$ This proves immediately that $\widehat{\calC}_{\geq 0}[1] \subseteq
\widehat{\calC}_{\geq 0}$ and $\widehat{\calC}_{\leq 0}[-1] \subseteq \widehat{\calC}_{\leq 0}$. 
Moreover, if $X \in \widehat{\calC}_{\geq 0}$ and $Y \in \widehat{\calC}_{\leq -1} = \widehat{\calC}_{\leq 0}[-1]$, then $\bHom_{\widehat{\calC}}( X,Y)$ can be identified with the homotopy limit of a tower
of spaces 
$$\ldots \rightarrow \bHom_{\calC}( X(n), Y(n) ) \rightarrow \bHom_{\calC}( X(n-1), Y(n-1)) \rightarrow \ldots$$ Since each of these spaces is contractible, we conclude that $\bHom_{\widehat{\calC}}(X,Y) \simeq \ast$; in particular,
$\Ext^0_{\widehat{\calC}}(X,Y) = 0$. Finally, we consider an arbitrary $X \in \widehat{\calC}$. Let
$X'' = \tau_{\leq -1} \circ X: \Nerve(\Z) \rightarrow \calC$, and let $u: X \rightarrow X''$ be the induced map. It is easy to check that $X'' \in \widehat{\calC}_{\leq -1}$ and that $\ker(u) \in \widehat{\calC}_{\geq 0}$. This completes the proof of $(2)$.

To prove $(3)$, we let $\calD$ denote the full subcategory of $\Nerve(\Z) \times \calC$ spanned by pairs $(n,C)$ such that $C \in \calC_{\leq -n}$. Using Proposition \toposref{testreflect}, we deduce that the inclusion $\calD \subseteq \Nerve(\Z) \times \calC$ admits a left adjoint $L$.
The composition
$$ \Nerve(\Z) \times \calC \stackrel{L}{\rightarrow} \calD \subseteq
\Nerve(\Z) \times \calC) \rightarrow \calC$$
can be identified with a functor $\theta: \calC \rightarrow \Fun( \Nerve(\Z), \calC)$ which factors
through $\widehat{\calC}$. To prove that $\theta$ is exact, it suffices to show that $\theta$ is right
exact (Proposition \ref{funrose}). Since the truncation functors
$\tau_{\leq n}: \calC_{\leq n+1} \rightarrow \calC_{\leq n}$ are right exact, finite colimits in
$\widehat{\calC}$ are computed pointwise. Consequently, it suffices to prove that each of
compositions
$$ \calC \stackrel{\theta}{\rightarrow} \widehat{\calC} \rightarrow \tau_{\leq n} \calC$$
is right exact. But this composition can be identified with the functor $\tau_{\leq n}$.

Finally, we observe that $\widehat{\calC}_{\leq 0}$ can be identified with a homotopy limit of the essentially constant tower
$$ \ldots \calC_{\leq 0} \stackrel{\id}{\rightarrow} \calC_{\leq 0} 
\stackrel{\id}{\rightarrow} \calC_{\leq 0} \stackrel{ \tau_{\leq -1} }{\rightarrow} \calC_{\leq -1}
\rightarrow \ldots,$$
and that $\theta$ induces an identification of this homotopy limit with $\calC_{\leq 0}$. 
\end{proof}

If $\calC$ is a stable $\infty$-category equipped with a t-structure, then we will say that
$\calC$ is {\it left complete} if the functor $\calC \rightarrow \widehat{\calC}$ described in Proposition \ref{sparrow} is an equivalence.

\begin{remark}
Let $\calC$ be as in Proposition \ref{sparrow}. Then the inclusion
$ \calC^{\dminus} \subseteq \calC$ induces an equivalence
$\widehat{ \calC}^{\dminus} \rightarrow \widehat{\calC}$, and
the functor $\calC \rightarrow \widehat{\calC}$ induces an equivalence
$ \calC^{\dminus} \rightarrow \widehat{\calC}^{\dminus}$. Consequently, the constructions
$$ \calC \mapsto \widehat{\calC}$$
$$ \calC \mapsto \calC^{\dminus}$$
furnish an equivalence between the theory of left bounded stable $\infty$-categories and the theory of left complete stable $\infty$-categories.
\end{remark}

We conclude this section with a useful criterion for establishing left completeness.

\begin{proposition}\label{cosparrow}
Let $\calC$ be a stable $\infty$-category equipped with a t-structure.
Suppose that $\calC$ admits countable products, and that
$\calC_{\geq 0}$ is stable under countable products. The
following conditions are equivalent:
\begin{itemize}
\item[$(1)$] The $\infty$-category $\calC$ is left complete.
\item[$(2)$] The full subcategory $\calC_{ \geq \infty } = \bigcap \calC_{ \geq n}
\subseteq \calC$ consists only of zero objects of $\calC$.
\end{itemize}
\end{proposition}

\begin{proof}
We first observe every tower of objects
$$ \ldots \rightarrow X_{n} \rightarrow X_{n-1} \rightarrow \ldots $$
in $\calC$ admits a limit $\varprojlim \{ X_{n} \}$: we can compute this limit as
the kernel of an appropriate map
$$ \prod X_{n} \rightarrow \prod X_{n}.$$
Moreover, if each $X_{n}$ belongs to $\calC_{\geq 0}$, then $\varprojlim \{ X_{n} \}$
belongs to $\calC_{ \geq -1}$. 

The functor $F: \calC \rightarrow \widehat{\calC}$ of Proposition \ref{sparrow} admits a
right adjoint $G$, given by 
$$ f \in \widehat{\calC} \subseteq \Fun( \Nerve(\Z), \calC) \mapsto \varprojlim(f).$$
Assertion $(1)$ is equivalent to the statement that the unit and counit maps
$$ u: F \circ G \rightarrow \id_{ \widehat{\calC}}$$
$$ v: \id_{\calC} \rightarrow G \circ F $$
are equivalences. If $v$ is an equivalence, then any object $X \in \calC$
can be recovered as the limit of the tower $\{ \tau_{ \leq n} X \}$. In particular,
this implies that $X = 0$ if $X \in \calC_{\geq \infty}$, so that $(1) \Rightarrow (2)$.

Now assume $(2)$; we will prove that $u$ and $v$ are both equivalences. To prove that
$u$ is an equivalence, we must show that for every
$f \in \widehat{\calC}$, the natural map
$$ \theta: \varprojlim(f) \rightarrow f(n)$$
induces an equivalence $\tau_{\leq -n} \varprojlim(f) \rightarrow f(n)$. In other words, we must
show that the kernel of $\theta$ belongs to $\calC_{\geq -n+1}$. To prove this, we first
observe that $\theta$ factors as a composition
$$ \varprojlim(f) \stackrel{\theta'}{\rightarrow} f(n-1) \stackrel{\theta''}{\rightarrow} f(n).$$
The octahedral axiom ($(TR4)$ of Definition \ref{deftriangle}) implies the existence
of an exact triangle
$$ \ker(\theta') \rightarrow \ker(\theta) \rightarrow \ker(\theta'').$$
Since $\ker(\theta'')$ clearly belongs to $\calC_{\geq -n+1}$, it will suffice to show that
$\ker(\theta')$ belongs to $\calC_{\geq -n+1}$. We observe that $\ker(\theta')$ can be identified
with the limit of a tower $\{ \ker( f(m) \rightarrow f(n-1) ) \}_{m < n}$. It therefore suffices to show that
each $\ker(f(m) \rightarrow f(n-1))$ belongs to $\calC_{\geq -n+2}$, which is clear. 

We now show prove that $v$ is an equivalence. Let $X$ be an object of $\calC$, and
$v_X: X \rightarrow (G \circ F)(X)$ the associated map. Since $u$ is an equivalence of functors, we conclude that $\tau_{\leq n}( v_X )$ is an equivalence for all $n \in \Z$. It follows that
$\coker(v_X) \in \calC_{ \geq n+1}$ for all $n \in \Z$. Invoking assumption $(2)$, we conclude that
$\coker(v_X) \simeq 0$, so that $v_X$ is an equivalence as desired.
\end{proof}

\begin{remark}
The ideas introduced above can be dualized in an obvious way, so that we can speak of {\it right completions} and {\it right completeness} for a stable $\infty$-category equipped with a t-structure.
\end{remark}

\section{Stabilization}\label{stable9.1}

In this section, we will describe a method for constructing stable $\infty$-categories:
namely, for any $\infty$-category $\calC$ which admits finite limits, one can consider an $\infty$-category $\Spectra(\calC)$ of {\it spectrum objects} of $\calC$. In the case where
$\calC$ is the $\infty$-category of spaces, we recover classical stable homotopy theory,
which we will discuss in \S \ref{stable8}.

\begin{definition}\label{spinn}
Let $\calC$ be an $\infty$-category. A {\it prespectrum object} of $\calC$
is a functor $X: \Nerve( \Z \times \Z) \rightarrow \calC$ with the following property:
for every pair of integers $i \neq j$, the value $X(i,j)$ is a zero object of $\calC$.
We let $\PreSpectra(\calC)$ denote the full subcategory of $\Fun( \Nerve(\Z \times \Z), \calC)$ spanned by the prespectrum objects of $\calC$.

For every integer $n$, evaluation at $(n,n) \in \Z \times \Z$ induces a functor
$\PreSpectra(\calC) \rightarrow \calC$. We will refer to this functor as the {\it $n$th space
functor} and denote it by $\Omega^{\infty-n}_{\calC}$.
\end{definition}

\begin{remark}
The partially ordered set $\Z \times \Z$ is isomorphic to its opposite, via the map
$(i,j) \mapsto (-i, -j)$. Composing with this map, we obtain an equivalence
$$ \PreSpectra(\calC)^{op} \simeq \PreSpectra(\calC^{op}).$$
\end{remark}

\begin{remark}
Let $X$ be a prespectrum object of an $\infty$-category $\calC$. Since the objects
$X(i,j) \in \calC$ are zero for $i \neq j$, it is customary to ignore them and instead emphasize
the objects $X(n,n) \in \calC$ lying along the diagonal. We will often denote
$X(n,n) = \Omega^{\infty-n}_{\calC} X$ by $X[n]$. For each $n \geq 0$, the diagram
$$ \xymatrix{ X(n,n) \ar[r] \ar[d] & X(n, n+1) \ar[d] \\
X(n+1, n) \ar[r] & X(n+1, n+1) }$$
determines an (adjoint) pair of morphisms
$$ \alpha: \Sigma_{\calC} X[n] \rightarrow X[n+1] \quad \quad \beta: X[n] \rightarrow \Omega_{\calC} X[n+1].$$
\end{remark}

\begin{definition}
Let $X$ be a prespectrum object of a pointed $\infty$-category $\calC$, and $n$ an integer.
We will say that $X$ is a {\it spectrum below $n$} if the canonical map
$\beta: X[m-1] \rightarrow \Omega_{\calC} X[m]$ is an equivalence for each $m \leq n$.
We say that $X$ is a {\it suspension prespectrum above $n$} if the canonical map
$\alpha: \Sigma_{\calC} X[m] \rightarrow X[m+1]$ is an equivalence for all $m \geq n$.
We say that $X$ is an {\it $n$-suspension prespectrum} if it is a
suspension prespectrum above $n$ and a spectrum below $n$.
We say that $X$ is a {\it spectrum object} if it is a spectrum object below $n$ for
all integers $n$. We let $\Spectra(\calC)$ denote the full subcategory of $\PreSpectra(\calC)$ spanned by the spectrum objects of $\calC$.

If $\calC$ is an arbitrary $\infty$-category, we let $\Stab(\calC) = \Spectra( \calC_{\ast} )$. Here
$\calC_{\ast}$ denotes the $\infty$-category of pointed objects of $\calC$. We will
refer to $\Stab(\calC)$ as the {\it stabilization} of $\calC$.
\end{definition}

\begin{remark}
Suppose that $\calC$ is a pointed $\infty$-category. Then the forgetful functor
$\calC_{\ast} \rightarrow \calC$ is a trivial Kan fibration, which induces a trivial
Kan fibration $\Stab(\calC) \rightarrow \Spectra(\calC)$.
\end{remark}

\begin{example}
Let $\Q$ be the ring of rational numbers, let $\bfA$ be the category of simplicial commutative $\Q$-algebras, viewed as simplicial model category (see Proposition \toposref{sutcoat}), and let
$\calC = \Nerve( \bfA^{\degree})$ be the underlying $\infty$-category. 
Suppose that $R$ is a commutative $\Q$-algebra, regarded as an object of $\calC$. Then
$\Stab( \calC_{/R} )$ is a stable $\infty$-category, whose homotopy category is equivalent to the (unbounded) derived category of $R$-modules. The loop functor
$\Omega^{\infty}: \Stab(\calC_{/R}) \rightarrow \calC_{/R}$ admits a left adjoint $\Sigma^{\infty}: \calC_{/R} \rightarrow \Stab(\calC_{/R})$ (Proposition \ref{urtusk22}). This left adjoint assigns to
each morphism of commutative rings $S \stackrel{\phi}{\rightarrow} R$ an object
$\Sigma^{\infty}(\phi) \in \Stab(\calC_{/R})$, which can be identified with $L_{S} \otimes_{S} R$, where $L_{S}$ denotes the (absolute) {\em cotangent complex} of $S$. We will discuss this example in greater detail in \cite{deformation}; see also \cite{schwede} for discussion.\index{cotangent complex}
\end{example}

\begin{remark}
Let $\calC$ be a pointed presentable $\infty$-category. Using Lemmas \toposref{stur2}, \toposref{stur3}, and \toposref{stur1}, we deduce that $\PreSpectra(\calC)$ and $\Spectra(\calC)$ are accessible localizations of $\Fun( \Nerve( \Z \times \Z), \calC)$. It follows that $\PreSpectra(\calC)$
and $\Spectra(\calC)$ are themselves presentable $\infty$-categories. Moreover, the inclusion
$\Spectra(\calC) \subseteq \PreSpectra(\calC)$ admits an accessible left adjoint $L_{\calC}$, which we will refer to as the {\it spectrification functor}. We will give a more direct construction of $L_{\calC}$ below in the case where $\calC$ satisfies some mild hypotheses.
\end{remark}

\begin{remark}\label{ulker}
Suppose that $\calC$ is a pointed $\infty$-category which admits finite limits and countable colimits, and that the loop functor $\Omega_{\calC}: \calC \rightarrow \calC$ preserves sequential colimits.
Then the collection of prespectrum objects of $\calC$ which are spectra below
$n$ is closed under sequential colimits.
\end{remark}

\begin{remark}\label{postulker}
The hypotheses of Remark \ref{ulker} are always satisfied in any of the following cases:
\begin{itemize}
\item[$(1)$] The $\infty$-category $\calC$ is pointed and compactly generated. 
\item[$(2)$] The $\infty$-category $\calC$ is an $\infty$-topos (Example \toposref{tucka}).
\item[$(3)$] The $\infty$-category $\calC$ is stable and admits countable coproducts.
In this case, Proposition \toposref{appendixdiagram} guarantees that $\calC$ admits all countable colimits, and the functor $\Omega_{\calC}$ is an equivalence and therefore preserves all colimits which exist in $\calC$.
\end{itemize}
\end{remark}

In order to work effectively with prespectrum objects, it is convenient to introduce a bit
of additional terminology.

\begin{notation}\label{jirl}
For $- \infty \leq a \leq b \leq \infty$, we let $Q(a,b) =
\{ (i,j) \in \Z \times \Z: (i \neq j) \vee (a \leq i = j \leq b) \}$. If
$\calC$ is an $\infty$-category, we let $\PreSpectra^{b}_{a}(\calC)$ denote the
full subcategory of $\Fun( \Nerve(Q(a,b)), \calC)$ spanned by those functors
$X$ such that $X(i,j)$ is a zero object of $\calC$ for $i \neq j$.
\end{notation}

\begin{lemma}\label{atwas1}
Let $\calC$ be a pointed $\infty$-category which admits finite limits,
and suppose given $\infty < a \leq b \leq \infty$. Let
$X_0 \in \PreSpectra^{b}_{a}(\calC)$. Then:
\begin{itemize}
\item[$(1)$] There exists an object $X \in \PreSpectra^{b}_{a-1}(\calC)$ which is a
right Kan extension of $X_0$.
\item[$(2)$] An arbitrary object $X \in \PreSpectra^{b}_{a-1}(\calC)$ which extends
$X_0$ is a right Kan extension of $X_0$ if and only if the induced map
$X[a-1] \rightarrow X[a]$. 
\end{itemize}
\end{lemma}

\begin{proof}
Note that $Q(a-1,b)$ is obtained from
$Q(a,b)$ by adjoining a single additional object $(a-1, a-1)$. It now suffices to observe that the the inclusion of $\infty$-categories
$$ \Nerve( \{ (a-1,a), (a,a), (a, a-1) \})^{op} \subseteq \Nerve( \{ (i,j) \in Q(a,b) | (a-1 \leq i,j) \})^{op}$$ is cofinal, which follows immediately from the criterion of Theorem \toposref{hollowtt}.
\end{proof}

\begin{lemma}\label{atwas2}
Let $\calC$ be a pointed $\infty$-category which admits finite limits
and suppose given $\infty < a \leq b \leq \infty$. Let
$X_0 \in \PreSpectra^{b}_{a}(\calC)$. Then:
\begin{itemize}
\item[$(1)$] There exists an object $X \in \PreSpectra^{b}_{- \infty}(\calC)$ which
is a right Kan extension of $X_0$.
\item[$(2)$] An arbitrary object $X \in \PreSpectra^{b}_{- \infty}(\calC)$ extending
$X_0$ is a right Kan extension of $X_0$ if and only if $X$ is a spectrum object
below $a$.
\end{itemize}
\end{lemma}

\begin{proof}
Combine Lemma \ref{atwas1} with Proposition \toposref{acekan}.
\end{proof}

\begin{lemma}\label{atwas3}
Let $\calC$ be a pointed $\infty$-category and $n$ an integer. Then
evaluation at $(n,n)$ induces a trivial Kan fibration
$\PreSpectra^{n}_{n}(\calC) \rightarrow \calC$.
\end{lemma}

\begin{proof}
Let $Q' = \{ (i,j) \in \Z \times \Z: (i < j \leq n) \vee (j < i \leq n) \vee (i=j=n) \}$,
and let $\calD \subseteq \Fun( \Nerve(Q'), \calC)$ denote the full subcategory
spanned by those functors $X$ such that $X(i,j)$ is a zero object of $\calC$ for
$i \neq j$. We observe the following:
\begin{itemize}
\item[$(a)$] A functor $X: \Nerve(Q') \rightarrow \calC$ belongs to $\calD$ if
and only if $X$ is a left Kan extension of $X | \{ (n,n) \}$.
\item[$(b)$] A functor $X: \Nerve( Q(n,n)) \rightarrow \calC$ belongs to
$\PreSpectra^{n}_{n}(\calC)$ if and only if $X| \Nerve(Q') \in \calD$ and
$X$ is a right Kan extension of $X | \Nerve(Q')$.
\end{itemize}
It now follows from Proposition \toposref{lklk} that the restriction functors
$$ \PreSpectra^{n}_{n}(\calC) \rightarrow \calD \rightarrow \calC$$
are trivial Kan fibrations, and the composition is given by evaluation at
$(n,n)$.
\end{proof}

We can now describe the stabilization $\Spectra(\calC)$ of a pointed $\infty$-category
$\calC$ in more conceptual terms:

\begin{proposition}\label{camer}
Let $\calC$ be a pointed $\infty$-category which admits finite limits. Then
the $\infty$-category $\Spectra(\calC)$ can be identified with the homotopy inverse limit of the tower
$$ \ldots \rightarrow \calC \stackrel{ \Omega_{\calC}}{\rightarrow} \calC \stackrel{ \Omega_{\calC}}{\rightarrow} \calC.$$
\end{proposition}

\begin{proof}
For every nonnegative integer $n$, let $\calD_n$ denote the full subcategory of
$\PreSpectra^{n}_{-\infty}(\calC)$ spanned by those functors $X$ such that the diagram
$$ \xymatrix{ X(m,m) \ar[r] \ar[d] & X(m,m+1) \ar[d] \\
X(m+1, m) \ar[r] & X(m+1, m+1) }$$
is a pullback square for each $m < n$. We note that Lemma \ref{atwas2} and Proposition
\toposref{lklk} imply that the composition
$$ \calD_n \subseteq \PreSpectra^{n}_{- \infty}(\calC) \rightarrow
\PreSpectra^{n}_{n}(\calC)$$
is a trivial Kan fibration. Combining this with Lemma \ref{atwas3}, we deduce that
evaluation at $(n,n)$ induces a trivial Kan fibration $\psi_n: \calD_n \rightarrow \calC$.
Let $s_n$ denote a section to $\psi_n$. We observe that the composite functor
$$ \calC \stackrel{s_n}{\rightarrow} \calD_n \rightarrow \calD_{n-1}
\stackrel{\psi_{n}}{\rightarrow} \calC$$
can be identified with the loop functor $\Omega_{\calC}$. It follows that the tower
$$ \ldots \rightarrow \calC \stackrel{ \Omega_{\calC}}{\rightarrow} \calC \stackrel{ \Omega_{\calC}}{\rightarrow} \calC.$$
is equivalent to the tower of restriction maps
$$ \ldots \rightarrow \calD_2 \rightarrow \calD_1 \rightarrow \calD_0.$$
This tower consists of categorical fibrations between $\infty$-categories, so its
homotopy inverse limit coincides with the actual inverse limit
$\varprojlim \{ \calD_n \}_{ n \geq 0} \simeq \Spectra(\calC)$.
\end{proof}

We now study the ``spectrification functor'' $\PreSpectra(\calC) \rightarrow \Spectra(\calC)$
in the case where $\calC$ is well-behaved.

\begin{lemma}\label{katen}
Let $P = \{ (i,j,k) \in \Z \times \Z \times \Z |
(i \neq j) \vee (i=j \geq k) \}$. Let $X: \Nerve(P) \rightarrow \calC$ be a
functor, where $\calC$ is a pointed $\infty$-category which admits finite limits. Suppose that $X(i,j,k)$ is a zero object of $\calC$ for all $i \neq j$. Then:

\begin{itemize}
\item[$(1)$] Let $\overline{X}: \Nerve( \Z \times \Z \times \Z) \rightarrow \calC$ be an extension of $X$. The following conditions are equivalent:
\begin{itemize}
\item[$(i)$] The functor $\overline{X}$ is a right Kan extension of $X$.
\item[$(ii)$] For each $k \geq 0$, the induced functor
$\overline{X} | \Nerve( \Z \times \Z \times \{k\})$ is a right
Kan extension of $X | \Nerve( Q(k, \infty) \times \{k\} )$.
\item[$(iii)$] For each $k \geq 0$, the functor $\overline{X} | \Nerve( \Z \times \Z \times \{k\})$ is a spectrum object below $k$.
\end{itemize}
\item[$(2)$] There exists a functor
$\overline{X}: \Nerve( \Z \times \Z \times \Z) \rightarrow \calC$
satisfying the equivalent conditions of $(1)$.
\end{itemize}
\end{lemma}

\begin{proof}
Let $\overline{X}$ be as in the statement of $(1)$. To prove the equivalence of conditions $(i)$ and $(ii)$, it will suffice to prove the following more precise claim: for every triple of nonnegative integers $(i,j,k)$, the functor $\overline{X}$ is a right Kan extension of $X$ at $(i,j,k)$ if and only if
$\overline{X} | \Nerve( \Z \times \Z \times \{k\})$ is a right Kan extension of
$X | \Nerve( Q(k, \infty) \times \{k\})$ at $(i,j,k)$. This follows from the observation that the inclusion
$$ \Nerve(\{ (i',j',k) \in P | (i' \geq i) \wedge (j' \geq j) \})^{op} \subseteq
\Nerve(\{ (i',j',k') \in P | (i' \geq i) \wedge (j' \geq j) \wedge (k' \geq k) \})^{op}$$
is cofinal (since it admits left adjoint, given by $(i',j',k') \mapsto (i',j',k)$).
The equivalence of $(ii)$ and $(iii)$ follows from Lemma \ref{atwas2}.

To prove $(2)$, we define a sequence of subsets
$$ P = P(0) \subseteq P(1) \subseteq P(2) \subseteq \ldots \subseteq \Z \times \Z \times \Z$$
by the formula $P(m) = \{ (i,j,k) \in \Z \times \Z \times \Z| (i \neq j) \vee
(i = j \geq k-m) \}.$
Using Proposition \toposref{acekan}, we deduce that if $\overline{X}: \Nerve( \Z \times \Z \times \Z) \rightarrow \calC$ is an extension of $X$, then $\overline{X}$ is a right
Kan extension of $X$ if and only if each restriction $\overline{X} | \Nerve( P(m+1) )$
is a right Kan extension of $\overline{X}| \Nerve(P(m) )$. Consequently, $(2)$
is a consequence of the following more precise assertion:
\begin{itemize}
\item[$(2')$] Every functor $Y_0: \Nerve( P(m) ) \rightarrow \calC$ extending $X$ admits a
right Kan extension $Y: \Nerve( P(m+1) ) \rightarrow \calC$ satisfying $(1')$.
\end{itemize}
To prove these assertions, we note that every element of $P(m+1)$ which does not belong
to $P(m)$ has the form $(n-1,n-1,n+m)$ for some integer $n$. It now suffices to
observe that the inclusion $\Nerve( P'_0)^{op} \subseteq \Nerve(P')^{op}$ is cofinal,
where $$P' = \{(i,j,k) \in P(m): (i,j \geq n-1) \wedge (k \geq n+m) \}$$
$$P'_0 = \{ (n-1, n, n+m), (n, n-1, n+m), (n, n, n+m) \} \subseteq P'.$$
This follows immediately from the criterion of Theorem \toposref{hollowtt}.
\end{proof}

\begin{corollary}\label{camber}
Let $\calC$ be a pointed $\infty$-category which admits finite limits. Then
there exists a sequence of functors 
$$ \id \rightarrow  L_0 \rightarrow L_1 \rightarrow L_2 \rightarrow \ldots$$
from $\PreSpectra(\calC)$ to itself such that the following conditions are satisfied
for every prespectrum object $X$ of $\calC$ and every $n \geq 0$:
\begin{itemize}
\item[$(1)$] The prespectrum $L_n(X)$ is a spectrum below $n$.

\item[$(2)$] The map
$\alpha: X \rightarrow L_{n}(X)$ induces an equivalence
$X[m] \rightarrow L_n(X)[m]$ for $m \geq n$.

\item[$(3)$] Suppose that $X$ is a spectrum below $n$. Then the map $\alpha: X \rightarrow L_n(X)$ is an equivalence.

\item[$(4)$] Let $Y$ be any prespectrum object of $\calC$ which is a spectrum below $n$.
Then composition with $\alpha$ induces a homotopy equivalence
$$ \bHom_{\PreSpectra(\calC)}( L_n(X), Y) \rightarrow \bHom_{ \PreSpectra(\calC)}(X, Y).$$
\end{itemize}
\end{corollary}

\begin{proof}
Let $P$ be defined as in Lemma \ref{katen}, let $\calD$ denote the full subcategory of
$\Fun( \Nerve(P), \calC)$ spanned by those functors $F$ such that $F(i,j,k)$ is a zero object of
$\calC$ for $i \neq j$, and let $\calD' \subseteq \Fun( \Nerve( \Z \times \Z \times \Z), \calC)$ denote the full subcategory spanned by those objects $F$ such that
$F$ is a right Kan extension of $F | \Nerve(P) \in \calD$. Using Lemma \ref{katen} and Proposition \toposref{lklk}, we deduce that the restriction map $\calD' \rightarrow \calD$ is a trivial Kan fibration.
Let $p$ denote a section of this restriction map.

Let $q: \PreSpectra(\calC) \rightarrow \calD$ denote the functor induced by the map
of partially ordered sets $P \rightarrow \Z \times \Z$
$(i,j,k) \mapsto (i,j)$. The composition $p \circ q: \PreSpectra(\calC) \rightarrow \calD'$
determines a map $\PreSpectra(\calC) \times \Nerve( \Z_{\geq 0}) \rightarrow \PreSpectra(\calC)$,
which we can identify with a sequence of functors $\{ L_n \}_{n \in \Z}$ from
$\PreSpectra(\calC)$ to itself. By construction, we also have a canonical map
$\id \rightarrow L_0$.
Assertions $(1)$ and $(2)$ are immediate consequences of the construction, and assertion $(3)$ follows from $(1)$ and $(2)$. To prove $(4)$, it will suffice (by virtue of $(3)$) to show that composition
with $\alpha$ induces a homotopy equivalence
$$ \bHom_{\PreSpectra(\calC)}( L_n(X), L_n(Y) ) \rightarrow \bHom_{ \PreSpectra(\calC)}(X, L_n(Y) ).$$
Since $\alpha$ induces an equivalence $X | \Nerve(Q(n,\infty)) \simeq L_{n}(X) | \Nerve(Q(n,\infty))$, it will suffice to show that $L_n(Y)$ is a right Kan extension of $L_n(Y)| \Nerve(Q(n,\infty))$, which follows from the equivalence of $(i)$ and $(ii)$ in Lemma \ref{atwas2}.
\end{proof}

\begin{corollary}\label{postzing}
Let $\calC$ be a pointed $\infty$-category satisfying the hypotheses of Remark \ref{ulker}. Let
$\{ L_n: \PreSpectra(\calC) \rightarrow \PreSpectra(\calC) \}_{n \geq 0}$ be the localization
functors of Corollary \ref{camber}. Then $L = \varinjlim_{n} L_n$ is a localization functor
from $\PreSpectra(\calC)$ to itself, whose essential image is the collection of spectrum objects
of $\calC$.
\end{corollary}

\begin{proof}
It will suffice to prove the following assertions for every prespectrum object $X$ of
$\calC$:
\begin{itemize}
\item[$(1)$] The prespectrum $LX$ is a spectrum object of $\calC$.
\item[$(2)$] For every spectrum object $Y$ of $\calC$, composition with the map
$\alpha: X = L_0 X \rightarrow LX$ induces a homotopy equivalence
$$ \bHom_{ \PreSpectra(\calC)}(LX, Y) \rightarrow \bHom_{ \PreSpectra(\calC)}(X,Y).$$
\end{itemize}
To prove $(1)$, it suffices to show that $LX$ is a prespectrum below $n$, for each $n \geq 0$.
Since $LX$ is a colimit of the sequence of prespectra $\{ L_m X\}_{m \geq n}$, each of which is a spectrum below $n$, this follows from Remark \ref{ulker}.

To prove $(2)$, we note that $\bHom_{ \PreSpectra(\calC)}( LX,Y)$ is given by the homotopy
inverse limit of a tower of spaces $\{ \bHom_{ \PreSpectra(\calC)}( L_n X, Y) \}_{n \geq 0}$.
Consequently, it will suffice to prove that each of the canonical maps 
$ \bHom_{ \PreSpectra(\calC)}(L_n X, Y) \rightarrow \bHom_{ \PreSpectra(\calC)}(X,Y)$
is a homotopy equivalence. This follows from Corollary \ref{camber}, since $Y$ is a spectrum below $n$.
\end{proof}

\begin{remark}\label{jilly}
Let $\calC$ be a pointed $\infty$-category satisfying the hypotheses of Remark \ref{ulker}, and let 
$f: X \rightarrow Y$ be a morphism between prespectrum objects of $\calC$.
Suppose that there exists an integer $n \geq 0$ such that $f$ induces an equivalence
$X[m] \rightarrow Y[m]$ for $m \geq n$. It follows that $L_{m}(f): L_m X \rightarrow L_m Y$
is an equivalence for $m \geq n$, so that $L(f) = \varinjlim_{m} L_m(f)$ is an equivalence
in $\Spectra(\calC)$.
\end{remark}

\begin{remark}\label{post}
Let $\calC$ be a presentable pointed $\infty$-category satisfying the hypotheses of Remark \ref{ulker}. Corollary \ref{postzing} asserts
that if $X$ is a prespectrum object of $\calC$, then the associated spectrum $X'$ is computed in the usual way: the $n$th space $X'[n]$ is given as a colimit $\colim_{m \geq 0} \Omega_{\calC}^{m} X[n+m]$.
\end{remark}

Suppose that $\calC$ is a presentable pointed $\infty$-category. We note that
$\Spectra(\calC)$ is closed under small limits in $\Fun( \Nerve(\Z \times \Z), \calC)$, so
that limits in $\Spectra(\calC)$ are computed pointwise. It follows that the evaluation functors
$\Omega^{\infty-n}_{\calC}: \Spectra(\calC) \rightarrow \calC$ preserve small limits. Since
these functors are also accessible, Corollary \toposref{adjointfunctor} guarnatees that
$\Omega^{\infty-n}_{\calC}$ admits a left adjoint, which we will denote by
$\Sigma^{\infty-n}_{\calC}: \calC \rightarrow \Spectra(\calC)$. Our next goal is to
describe this functor in more explicit terms.

\begin{lemma}\label{atwas4}
Let $\calC$ be a pointed $\infty$-category which admits finite limits and colimits.
Then:
\begin{itemize}
\item[$(1)$] A prespectrum object $X$ of $\calC$ is a suspension prespectrum above $n$ if and only if $X$ is a left Kan extension of $X | Q(-\infty,n)$.
\item[$(2)$] For every $X_0 \in \PreSpectra^{n}_{-\infty}(\calC)$, there exists
an extension $X \in \PreSpectra(\calC)$ of $X_0$ which satisfies the equivalent conditions of
$(1)$.
\item[$(3)$] Let $\calD$ denote the full subcategory of $\PreSpectra(\calC)$ spanned by those prespectrum objects of $\calC$ which are suspension prespectra above $n$. 
Then $\calD$ is a colocalization of $\PreSpectra(\calC)$.
Moreover, a morphism of prespectra $X \rightarrow Y$ exhibits $X$ as a $\calD$-colocalization of $Y$ if and only if $X \in \calD$ and the the induced map $X[k] \rightarrow Y[k]$ is an equivalence for $k \leq n$.
\item[$(4)$] Let $\calD_0$ denote the full subcategory of $\PreSpectra(\calC)$ spanned by
the $n$-suspension prespectrum objects, and let $\calE$ denote the full subcategory
of $\PreSpectra^{n}_{-\infty}(\calC)$ spanned by those functors $X$ such
that the induced map $X[k] \rightarrow \Omega_{\calC} X[k+1]$ is an equivalence for $k < n$.
Then the restriction maps $\calD \rightarrow \PreSpectra(\calC)$ and
$\calD_0 \rightarrow \calE$ are trivial Kan fibrations.
\item[$(5)$] The $\infty$-category $\calE$ is a localization of
$\PreSpectra^{n}_{-\infty}(\calC)$. Moreover, a morphism $X \rightarrow Y$ in
$\PreSpectra^{n}_{- \infty}(\calC)$ exhibits $Y$ as an $\calE$-localization of $X$
if and only if $Y \in \calE$ and the map $X[n] \rightarrow Y[n]$ is an equivalence.
\item[$(6)$] The $\infty$-category $\calD_0$ is a localization of
$\calD$. Moreover, a morphism $X \rightarrow Y$ in
$\calD$ exhibits $Y$ as an $\calD_0$-localization of $X$
if and only if $Y \in \calD_0$ and the map $X[n] \rightarrow Y[n]$ is an equivalence.
\end{itemize}
\end{lemma}

\begin{proof}
Assertions $(1)$ and $(2)$ follow by applying Lemma \ref{atwas2} to the opposite
$\infty$-category $\calC^{op}$. Assertion $(4)$ follows from $(1)$ and $(2)$
together with Proposition \toposref{lklk}, and assertion $(6)$ follows immediately from
$(5)$ and $(4)$. We will give the proof of $(3)$; the proof of $(5)$ is similar.

Consider an arbitrary object $Y \in \PreSpectra(\calC)$. Let
$Y_0 = Y | \Nerve(Q(-\infty,n))$, and let $X \in \PreSpectra(\calC)$ be a left Kan extension of $Y_0$
(whose existence is guaranteed by $(2)$). Then
$X \in \calD$ and we have a canonical map
$\alpha: X \rightarrow Y$. We claim that $\alpha$ exhibits $X$ as a a
$\calD$-colocalization of $Y$. To prove this, let us consider
an arbitrary object $W \in \calD$. We have a commutative diagram
$$ \xymatrix{ \bHom_{ \PreSpectra(\calC)}(W,X) \ar[r] \ar[d] & 
\bHom_{ \PreSpectra^{n}_{-\infty}(\calC)}( W | \Nerve(Q(-\infty,n)), X | \Nerve(Q(-\infty,n))) \ar[d] \\
\bHom_{ \PreSpectra(\calC)}(W,Y) \ar[r] & 
\bHom_{ \PreSpectra^{n}_{-\infty}(\calC)}( W | \Nerve(Q(-\infty,n)), Y | \Nerve(Q(-\infty,n))). }$$
Since $X$ and $Y$ have the same restriction to $\Nerve(Q(-\infty,n))$, the right vertical map
is a homotopy equivalence. The horizontal maps are homotopy equivalences since
$W$ is a left Kan extension of $W | \Nerve(Q(-\infty,n))$, by virtue of $(1)$. This completes the
proof that $\alpha$ exhibits $X$ as a $\calD$-colocalization of $Y$,
and the proof that $\calD$ is a colocalization of $\PreSpectra(\calC)$.

To complete the proof of $(3)$, let us consider an arbitrary map $\beta: X' \rightarrow Y$, where
$X' \in \calD$. We wish to show that $\beta$ exhibits
$X'$ as a $\calD$-colocalization of $Y$ if and only if the induced map
$X'[m] \rightarrow Y[m]$ is an equivalence for $m \leq n$. The above argument shows that
$\beta$ fits into a commutative triangle
$$ \xymatrix{ & X \ar[dr]^{\alpha} & \\
X' \ar[rr]^{\beta} \ar[ur]^{\gamma} & & Y. }$$
Since $\alpha$ exhibits $X$ as a $\calD$-colocalization of $Y$,
and $\alpha$ induces equivalences $X[m] \rightarrow Y[m]$ for $m \leq n$, we can restate
the desired assertion as follows: the map $\gamma$ is an equivalence if and only if
$\gamma$ induces equivalences $X[m] \rightarrow X'[m]$ for $m \leq n$. The ``only if'' part of the assertion is obvious, and the converse follows from the fact that both $X$ and $X'$ are left Kan
extensions of their restrictions to $\Nerve(Q(-\infty,n))$ (by virtue of $(1)$).
\end{proof}

\begin{proposition}\label{swagger}
Let $\calC$ be a pointed $\infty$-category which admits finite limits and colimits, and let
$\calD_0 \subseteq \PreSpectra(\calC)$ denote the full subcategory spanned by the
$n$-suspension prespectra. Then evaluation at $(n,n)$ induces a trivial Kan fibration
$\calD_0 \rightarrow \calC$.
\end{proposition}

\begin{proof}
Let $\calE \subseteq \PreSpectra^{n}_{-\infty}(\calC)$ be defined as in
Lemma \ref{atwas4}. The evaluation functor factors as a composition
$$ \calD_0 \stackrel{\phi_0}{\rightarrow} \calE
\stackrel{\phi_1}{\rightarrow} \PreSpectra^{n}_{n}(\calC)
\stackrel{\phi_2}{\rightarrow} \calC.$$
Here $\phi_0$ is a trivial fibration by Lemma \ref{atwas4}, the map
$\phi_1$ is a trivial fibration by virtue of Lemma \ref{atwas2} and Proposition
\toposref{lklk}, and the map $\phi_2$ is a trivial fibration by Lemma \ref{atwas3}.
\end{proof}

\begin{notation}
Let $\calC$ be a pointed $\infty$-category which admits finite limits and colimits.
We let $\widetilde{\Sigma}^{\infty - n}_{\calC}: \calC \rightarrow \PreSpectra(\calC)$ denote a
section of the trivial Kan fibration $\calD_0 \rightarrow \calC$ of
Proposition \ref{swagger}.
\end{notation}

\begin{remark}\label{inwise}
Let $\calC$ be a pointed $\infty$-category which admits finite limits and colimits,
let $C \in \calC$ be an object and let $X \in \PreSpectra(\calC)$ be a spectrum below $n$.
Then the canonical map
$$ e: \bHom_{ \PreSpectra(\calC)}( \widetilde{\Sigma}^{\infty - n}_{\calC} C, X)
\rightarrow \bHom_{ \calC}( C, \Omega^{\infty-n}_{\calC} X)$$
is a homotopy equivalence. To prove this, we observe that
$e$ factors as a composition (using the conventions of Notation \ref{jirl})
\begin{eqnarray*}
\bHom_{ \PreSpectra(\calC)}( \widetilde{ \Sigma}^{\infty -n}_{\calC} C, X)
& \stackrel{\phi_0}{\rightarrow} & \bHom_{ \PreSpectra^{n}_{-\infty}(\calC)}( 
(\widetilde{ \Sigma}^{\infty -n}_{\calC} C) | \Nerve(Q(-\infty,n)), X | \Nerve(Q(-\infty,n)) ) \\
& \stackrel{\phi_1}{\rightarrow} &  \bHom_{ \PreSpectra^{n}_{n}(\calC)}( 
(\widetilde{ \Sigma}^{\infty -n}_{\calC} C) | \Nerve(Q(n,n)), X | \Nerve(Q(n,n)) ) \\
& \stackrel{\phi_2}{\rightarrow} & \bHom_{\calC}( C, X[n] ). \end{eqnarray*}
It will therefore suffice to prove that the maps $\phi_0$, $\phi_1$, and
$\phi_2$ are homotopy equivalences. For $\phi_0$, the desired result
follows from our assumption that
$\widetilde{\Sigma}^{\infty-n}_{\calC} C$ is a left Kan extension of its restriction
to $\Nerve(Q(-\infty,n))$. For $\phi_1$, we invoke the fact that$X | \Nerve(Q(-\infty,n))$ is a right Kan extension of its restriction to $\Nerve(Q(n,n))$. For $\phi_2$, we apply Lemma \ref{atwas3}.
\end{remark}

\begin{proposition}\label{kuup}
Let $\calC$ be a presentable pointed $\infty$-category, and let
$L: \PreSpectra(\calC) \rightarrow \Spectra(\calC)$ denote a left adjoint to the inclusion.
Then the evaluation functor $\Omega^{\infty-n}_{\calC}: \Spectra(\calC) \rightarrow \calC$
admits a left adjoint, given by the composition
$$ \calC \stackrel{ \widetilde{\Sigma}^{\infty-n}_{\calC}}{\rightarrow} \PreSpectra(\calC)
\stackrel{L}{\rightarrow} \Spectra(\calC).$$
\end{proposition}

\begin{proof}
The canonical natural transformation $\id_{ \PreSpectra(\calC)} \rightarrow L$
induces a transformation 
$$\alpha: \id_{\calC} = \Omega^{\infty-n}_{\calC} \circ \widetilde{\Sigma}^{\infty-n}_{\calC}
\rightarrow \Omega^{\infty-n}_{\calC} \circ (L \circ \widetilde{\Sigma}^{\infty-n}_{\calC}).$$
We claim that $\alpha$ is the unit of an adjunction. To prove this, it suffices to show that
for every object $C \in \calC$ and every spectrum object $X \in \Spectra(\calC)$, the
composite map
$$ \bHom_{ \Spectra(\calC) }( L \widetilde{\Sigma}^{\infty-n}_{\calC} C, X)
\stackrel{\phi}{\rightarrow} \bHom_{ \PreSpectra(\calC)}( \widetilde{\Sigma}^{\infty-n}_{\calC} C, X)
\stackrel{\psi}{\rightarrow} \bHom_{ \calC}( C, X[n])$$ 
is a homotopy equivalence. It now suffices to observe that $\phi$ is a homotopy equivalence because
$X$ is a spectrum object, and $\psi$ is a homotopy equivalence by virtue of Remark \ref{inwise}.
\end{proof}

We close this section by discussing the {\em shift} functor on prespectrum objects of
an $\infty$-category $\calC$. We observe that precomposition with the map $(i,j) \mapsto (i+1, j+1)$
determines a functor $S: \PreSpectra(\calC) \rightarrow \PreSpectra(\calC)$, which
restricts to a functor $\Spectra(\calC) \rightarrow \Spectra(\calC)$ which we will also denote by $S$.
By construction, this functor is an equivalence (in fact, an isomorphism of simplicial sets).
We observe that if $X$ is a spectrum object of $\calC$, then we have canonical
equivalences $\Omega_{\calC} S(X)[n] = \Omega_{\calC} X[n+1] \simeq X$;
this strongly suggests that $S$ is a homotopy inverse to the loop functor
$\Omega_{ \Spectra(\calC) }$ given by pointwise composition with $\Omega_{\calC}$.
To prove this (in a slightly stronger form), we need to introduce a bit of notation.

\begin{notation}\label{corpae}
Consider the order-preserving maps $s_{+}, s_{-}: \Z \times \Z \rightarrow \Z \times \Z$ defined
by the formulae
$$ s_{+}(i,j) = \begin{cases} (i,j) & \text{if } i \neq j \\
(i+1, j) & \text{if } i = j.\end{cases} \quad \quad s_{-}(i,j) = \begin{cases} (i,j) & \text{if } i \neq j \\
(i, j+1) & \text{if } i = j.\end{cases}$$
For every $\infty$-category $\calC$, composition with $s_{+}$ and $s_{-}$ induces functors
$S_{+}, S_{-}: \PreSpectra(\calC) \rightarrow \PreSpectra(\calC)$, fitting into a commutative
diagram
$$ \xymatrix{ \id \ar[r] \ar[d] & S_{+} \ar[d] \\
S_{-} \ar[r] & S. }$$
Note that the images of $s_{+}$ and $s_{-}$ are disjoint from the diagonal
$\{ (n,n) \}_{n \geq 0} \subseteq \Z_{\geq 0} \times \Z_{\geq 0}$, so that
$S_{+}(X)$ and $S_{-}(X)$ are zero objects of $\PreSpectra(\calC)$ for every
$X \in \PreSpectra(\calC)$. If $\calC$ admits finite limits, then
the above diagram determines a morphism
$\alpha: X \rightarrow \Omega_{\PreSpectra(\calC)} S(X)$. If
$\calC$ also admits finite colimits, then $\alpha$ admits an adjoint
$\beta: \Sigma_{ \PreSpectra(\calC)} X \rightarrow S(X)$.
\end{notation}

\begin{lemma}\label{kija}
Let $\calC$ be a small pointed $\infty$-category, and let
$\calP_{\ast}(\calC)$ denote the full subcategory of
$\calP(\calC) = \Fun( \calC^{op}, \SSet)$ spanned by those functors which
carry zero objects of $\calC$ to final objects of $\SSet$.
Then:
\begin{itemize}
\item[$(1)$] Let $S$ denote the set consisting of a single morphism from
an initial object of $\calP(\calC)$ to a final object of $\calP(\calC)$. Then
$\calP_{\ast}(\calC) = S^{-1} \calP(\calC)$.
\item[$(2)$] The $\infty$-category $\calP_{\ast}(\calC)$ is an accessible localization
of $\calP(\calC)$. In particular, $\calP_{\ast}(\calC)$ is presentable.
\item[$(3)$] The Yoneda embedding $\calC \rightarrow \calP(\calC)$ factors through
$\calP_{\ast}(\calC)$, and the induced embedding $j: \calC \rightarrow \calP_{\ast}(\calC)$
preserves zero objects.
\item[$(4)$] Let $\calD$ be an $\infty$-category which admits small colimits, and let
$\Fun^{L}( \calP_{\ast}(\calC), \calD)$ denote the full subcategory of $\Fun( \calP_{\ast}(\calC), \calD)$
spanned by those functors which preserve small colimits. Then composition with
$j$ induces an equivalence of $\infty$-categories
$\Fun^{L}( \calP_{\ast}(\calC), \calD) \rightarrow \Fun_0( \calC, \calD)$, where
$\Fun_0(\calC, \calD)$ denotes the full subcategory of $\Fun( \calC, \calD)$ spanned by
those functors which carry zero objects of $\calC$ to initial objects of $\calD$.

\item[$(5)$] The $\infty$-category $\calP_{\ast}(\calC)$ is pointed.

\item[$(6)$] The full subcategory $\calP_{\ast}(\calC) \subseteq \calP(\calC)$ is closed
under small limits and under small colimits parametrized by weakly contractible simplicial sets.
In particular, $\calP_{\ast}(\calC)$ is stable under small filtered colimits in $\calP(\calC)$.

\item[$(7)$] The functor $j: \calC \rightarrow \calP_{\ast}(\calC)$ preserves all small limits
which exist in $\calC$.

\item[$(8)$] The $\infty$-category $\calP_{\ast}(\calC)$ is compactly generated.
\end{itemize}
\end{lemma}

\begin{proof}
For every object $X \in \SSet$, let $F_X \in \calP(\calC)$ denote the constant functor taking the
value $X$. Then $F_X$ is a left Kan extension of $F_X | \{0\}$, where $0$ denotes a zero object
of $\calC$. It follows that for any object $G \in \calP(\calC)$, evaluation at $0$ induces a homotopy equivalence
$$ \bHom_{ \calP(\calC)}( F_X, G) \rightarrow \bHom_{\SSet}( F_X(0), G(0) ) = \bHom_{\SSet}(X, G(0) ).$$
We observe that the inclusion $\emptyset \subseteq \Delta^0$ induces a map
$F_{\emptyset} \rightarrow F_{\Delta^0}$ from an initial object of $\calP(\calC)$ to a final object
of $\calP(\calC)$. It follows that an object $G$ of $\calP(\calC)$ is $S$-local if and only if
the induced map
$$ G(0) \simeq \bHom_{\SSet}( \Delta^0, G(0) ) \rightarrow \bHom_{\SSet}( \emptyset, G(0) ) \simeq \Delta^0$$
is a homotopy equivalence: that is, if and only if $G \in \calP_{\ast}(\calC)$. This proves $(1)$.

Assertion $(2)$ follows immediately from $(1)$, and assertion $(3)$ is obvious. Assertion
$(4)$ follows from $(1)$, Theorem \toposref{charpresheaf}, and Proposition \toposref{unichar}.
To prove $(5)$, we observe that $F_{\Delta^0}$ is a final object of $\calP(\calC)$, and therefore
a final object of $\calP_{\ast}(\calC)$. It therefore suffices to show that $F_{\Delta^0}$ is an initial
object of $\calP_{\ast}(\calC)$. This follows from the observation that for every $G \in \calP(\calC)$, we have homotopy equivalences $\bHom_{ \calP(\calC)}(F_{\Delta^0}, G) \simeq \bHom_{\SSet}( \Delta^0, G(0)) \simeq G(0)$ so that the mapping space $\bHom_{ \calP(\calC)}(F_{\Delta^0}, G)$ is contractible
if $G \in \calP_{\ast}(\calC)$.

Assertion $(6)$ is obvious, and $(7)$ follows from $(6)$ together with Proposition \toposref{yonedaprop}. We deduce $(8)$ from $(6)$ together with Corollary \toposref{starmin}.
\end{proof}

\begin{proposition}\label{sinkle}
Let $\calC$ be a pointed $\infty$-category which admits finite limits. Then:
\begin{itemize}
\item[$(1)$] For every object $X \in \Spectra(\calC)$, the canonical map
$X \rightarrow \Omega_{ \PreSpectra(\calC)} S(X)$ is an equivalence.
\item[$(2)$] The shift functor $S: \Spectra(\calC) \rightarrow \Spectra(\calC)$
is a homotopy inverse to the loop functor $\Omega_{ \Spectra(\calC) }$.
\item[$(3)$] The $\infty$-category $\Spectra(\calC)$ is stable.
\end{itemize}
\end{proposition}

\begin{proof}
Assertion $(1)$ is an immediate consequence of the definitions.
We note that $(1)$ implies that $S$ is a right homotopy inverse to
$\Omega_{ \Spectra(\calC)}$. Since $S$ is invertible, it follows that
$S$ is also a left homotopy inverse to $\Omega_{ \Spectra(\calC)}$.
In particular, $\Omega_{ \Spectra(\calC)}$ is invertible.

To prove $(3)$, we may assume without loss of generality
that $\calC$ is small. Lemma \ref{kija} implies that there exists
a fully faithful left exact functor $j: \calC \rightarrow \calD$, where
$\calD$ is a compactly generated pointed $\infty$-category
(this that the functor $\Omega_{\calD}$ preserves sequential colimits; see Remark \ref{postulker}).
Then $\Spectra(\calC)$ is equivalent to a full subcategory
of $\Spectra(\calD)$, which is closed under finite limits and shifts.
Consequently, it will suffice to show that $\Spectra(\calD)$ is stable,
which is a consequence of Corollary \ref{taskwise} (proven in \S \ref{stable9.2}).
\end{proof}

\begin{corollary}\label{charstut}
Let $\calC$ be a pointed $\infty$-category. The
following conditions are equivalent:

\begin{itemize}
\item[$(1)$] The $\infty$-category $\calC$ is stable.

\item[$(2)$] The $\infty$-category $\calC$ admits finite colimits and the
suspension functor $\Sigma: \calC \rightarrow \calC$ is an equivalence.

\item[$(3)$] The $\infty$-category $\calC$ admits finite limits and the loop
functor $\Omega: \calC \rightarrow \calC$ is an equivalence.
\end{itemize}
\end{corollary}

\begin{proof}
We will show that $(1) \Leftrightarrow (3)$; the dual argument will prove that $(1) \Leftrightarrow (2)$. The implication $(1) \Rightarrow (3)$ is clear. Conversely, suppose that $\calC$ admits finite limits and that $\Omega$ is an equivalence. Lemma \toposref{pointerprime} asserts that the forgetful functor $\calC_{\ast} \rightarrow \calC$ is a trivial fibration. Consequently, $\Spectra(\calC)$ can be identified with the homotopy inverse limit of the tower
$$ \ldots \stackrel{\Omega}{\rightarrow} \calC \stackrel{\Omega}{\rightarrow} \calC.$$
By assumption, the loop functor $\Omega$ is an equivalence, so this tower is essentially constant. It follows that $\Omega^{\infty}: \Spectra(\calC) \rightarrow \calC$ is an equivalence of $\infty$-categories. Since $\Spectra(\calC)$ is stable (Proposition \ref{sinkle}), so is $\calC$. 
\end{proof}

For later use, we record also the following result:

\begin{proposition}\label{makesig}
Let $\calC$ be a pointed $\infty$-category satisfying the hypotheses of Remark \ref{ulker}, and let $L: \PreSpectra(\calC) \rightarrow \Spectra(\calC)$ denote a left adjoint to the inclusion. Let $X$ be a prespectrum object of
$\calC$, and let $\beta: \Sigma_{ \PreSpectra(\calC)} X \rightarrow S(X)$ denote the map described in
Notation \ref{corpae}. Then $L(\beta)$ is an equivalence in $\Spectra(\calC)$.
\end{proposition}

\begin{proof}
Since $L$ is a left adjoint, it commutes with suspensions. It will therefore suffice to show that $\beta$
induces an equivalence $\Sigma_{ \Spectra(\calC)} LX \rightarrow L S(X)$: in other words, we must
show that the diagram $\sigma:$
$$ \xymatrix{ LX \ar[r] \ar[d] & L S_{+}(X) \ar[d] \\
L S_{-}(X) \ar[r] & L S(X) }$$
is a pushout square in the $\infty$-category $\Spectra(\calC)$. Since $\Spectra(\calC)$ is stable,
it suffices to show that $\sigma$ is a pullback square. In other words, we must show that
for each $n \geq 0$, the diagram
$$ \xymatrix{ LX[n] \ar[r] \ar[d] & LS_{+}(X)[n] \ar[d] \\
L S_{-}(X)[n] \ar[r] & L S(X)[n]}$$
is a pullback square in $\calC$. 

Let $P$ be defined as in Lemma \ref{katen}, let $\overline{X}_0: \Nerve(P) \rightarrow \calC$
be given by the composition
$$ \Nerve(P) \rightarrow \Nerve( \Z \times \Z) \stackrel{X}{\rightarrow} \calC,$$
and let $\overline{X}: \Nerve( \Z \times \Z \times \Z) \rightarrow \calC$
be a right Kan extension of $\overline{X}_0$. 
Let $\overline{S}( \overline{X})$, $\overline{S}_{+}( \overline{X})$, and $\overline{S}_{-}( \overline{X})$
be obtained from $\overline{X}$ by composing with the maps $\Z \times \Z \times \Z \rightarrow \Z \times \Z \times \Z$ given by
$(i,j,k) \mapsto (i+1, j+1, k)$, $s_{+} \times \id$, and $s_{-} \times \id$. We have a commutative diagram
$$ \xymatrix{ \overline{X} \ar[r] \ar[d] & \overline{S}_{+}(\overline{X} ) \ar[d] \\
\overline{S}_{-}( \overline{X} ) \ar[r] & \overline{S}( \overline{X} ) }$$
which we can think of as encoding a sequence of commutative squares 
$\{ \sigma_k: \Delta^1 \times \Delta^1 \rightarrow \PreSpectra(\calC) \}_{k \geq 0}$.
We can identify $\sigma$ with the colimit of this sequence. Consequently, it will suffice
to prove that for every integer $n$, the diagram
$$ \xymatrix{ \overline{X}(n,n,k) \ar[r] \ar[d] & \overline{S}_{+}( \overline{X})(n,n,k) \ar[d] \\
\overline{S}_{-}( \overline{X} )(n,n,k) \ar[r] & \overline{S}( \overline{X} )(n,n,k) }$$
is a pullback square in $\calC$ for all sufficiently large $k$. In fact, this is true for
all $k > n$, by virtue of Lemma \ref{katen}.
\end{proof}
 
\section{The $\infty$-Category of Spectra}\label{stable8}

In this section, we will discuss what is perhaps the most important example of a stable $\infty$-category: the $\infty$-category of spectra. In classical homotopy theory, one defines a spectrum to be a sequence of pointed spaces $\{ X_n \}_{n \geq 0}$, equipped with homotopy equivalences (or homeomorphisms, depending on the author) $X_{n} \rightarrow \Omega(X_{n+1})$ for all $n \geq 0$. 
This admits an immediate $\infty$-categorical translation:

\begin{definition}
A {\it spectrum} is a spectrum object of the $\infty$-category $\SSet_{\ast}$ of pointed spaces.
We let $\Spectra = \Spectra( \SSet_{\ast}) = \Stab(\SSet)$ denote the $\infty$-category of spectra.
\end{definition}

\begin{proposition}\label{specster1}
\item[$(1)$] The $\infty$-category
$\Spectra$ is stable.

\item[$(2)$] Let $(\Spectra)_{ \leq -1 }$ denote the full subcategory of
$\Spectra$ spanned by those objects $X$ such that $\Omega^{\infty}(X) \in \SSet$
is contractible. Then $(\Spectra)_{ \leq -1}$ determines an accessible t-structure
on $\Spectra$ (see \S \ref{stable16}). 

\item[$(3)$] The t-structure on $\Spectra$ is both left complete and right complete, and the heart
$\SSet_{\infty}^{\heartsuit}$ is canonically equivalent to the (nerve of the) category of abelian groups.

\end{proposition}

\begin{proof}
Assertion $(1)$ follows immediately from Proposition \ref{sinkle}.
Assertion $(2)$ is a special case of Proposition \ref{tmonster}, which will be established in \S \ref{stable16}. We will prove $(3)$. Note that a spectrum $X$ can be identified with a sequence of pointed spaces $\{ X(n) \}$, equipped with equivalences $X(n) \simeq \Omega X(n+1)$ for all $n \geq 0$. We observe that $X \in (\Spectra)_{\leq m}$ if and only if each $X(n)$ is $(n+m)$-truncated. In general, the sequence $\{ \tau_{\leq n+m} X(n) \}$ itself determines a spectrum, which we can identify with the truncation
$\tau_{\leq m} X$. It follows that $X \in (\Spectra)_{\geq m+1}$ if and only if each
$X(n)$ is $(n+m+1)$-connective. In particular, $X$ lies in the heart of $\Spectra$ if and only if
each $X(n)$ is an Eilenberg-MacLane object of $\SSet$ of degree $n$ (see Definition \toposref{gropab}). It
follows that the heart of $\Spectra$ can be identified with the homotopy inverse limit of the tower of $\infty$-categories
$$ \ldots \stackrel{\Omega}{\rightarrow} \EM_1(\SSet) \stackrel{\Omega}{\rightarrow} \EM_0(\SSet),$$
where $\EM_n(\SSet)$ denotes the full subcategory of $\SSet_{\ast}$ spanned by the Eilenberg-MacLane objects of degree $n$. Proposition \toposref{EM} asserts that after the second term, this tower is equivalent to the constant diagram taking the value $\Nerve( \Ab)$, where $\Ab$ is category of abelian groups.

It remains to prove that $\Spectra$ is both right and left complete. We begin by observing that
if $X \in \Spectra$ is such that $\pi_n X \simeq 0$ for all $n \in \Z$, then $X$ is a zero object of
$\Spectra$ (since each $X(n) \in \SSet$ has vanishing homotopy groups, and is therefore contractible by Whitehead's theorem). Consequently, both $\bigcap (\Spectra)_{\leq -n}$ and $\bigcap (\Spectra)_{\geq n}$ coincide with the collection of zero objects of $\Spectra$. It follows that
$$ (\Spectra)_{\geq 0} = \{ X \in \Spectra : (\forall n < 0) [ \pi_n X \simeq 0] \}$$
$$ (\Spectra)_{\leq 0} = \{ X \in \Spectra : (\forall n > 0) [\pi_n X \simeq 0 ] \}.$$
According to Proposition \ref{cosparrow}, to prove that $\Spectra$ is left and right complete it will suffice to show that the subcategories $( \Spectra)_{\geq 0}$ and $( \Spectra )_{\leq 0}$ are stable under products and coproducts. In view of the above formulas, it will suffice to show that the homotopy group functors $\pi_n: \Spectra \rightarrow \Nerve(\Ab)$ preserve products and coproducts. Since $\pi_n$ obviously commutes with finite coproducts, it will suffice to show that $\pi_n$ commutes with products and filtered colimits. Shifting if necessary, we may reduce to the case $n=0$. Since products and filtered colimits in the category of abelian groups can be computed at the level of the underlying sets, we are reduced to proving that the composition
$$ \Spectra \stackrel{\Omega^{\infty}}{\rightarrow} \SSet \stackrel{\pi_0}{\rightarrow} \Nerve(\Set)$$
preserves products and filtered colimits. This is clear, since each of the factors individually preserves products and filtered colimits.
\end{proof}

Our next goal is to prove that the $\infty$-category $\Spectra$ is compactly generated.
To prove this, we need to review a bit of the theory of finite spaces.

\begin{notation}\index{ZZZFinnSpace@$\FinnSpace$}\index{space!finite}\index{ZZFinSpace@$\FinSpace$}\label{finner}
Let $\SSet_{\ast}$ denote the $\infty$-category of pointed objects of $\SSet$. That is,
$\SSet_{\ast}$ denotes the full subcategory of $\Fun(\Delta^1, \SSet)$ spanned by those morphisms
$f: X \rightarrow Y$ for which $X$ is a final object of $\SSet$ (Definition \toposref{gropab}).  
Let $\FinnSpace$ denote the smallest full subcategory of $\SSet$ which contains
the final object $\ast$ and is stable under finite colimits. We will refer to $\FinnSpace$ as the {\it $\infty$-category of finite spaces}. We let $\FinSpace \subseteq \SSet_{\ast}$ denote the $\infty$-category of pointed objects of $\FinnSpace$. We observe that the suspension functor $\Sigma: \SSet_{\ast} \rightarrow \SSet_{\ast}$ carries $\FinSpace$ to itself. For each $n \geq 0$, we let $S^n \in \SSet_{\ast}$ denote a representative for the (pointed) $n$-sphere.
\end{notation}

\begin{remark}\label{yogause}\index{ZZZFunRexx@$\Fun^{\Rexx}(\calC, \calD)$}
It follows from Remark \toposref{poweryoga} and Proposition \toposref{lklk} that
$\FinnSpace$ is characterized by the following universal property: for every $\infty$-category $\calD$ which admits finite colimits, evaluation at $\ast$ induces an equivalence of $\infty$-categories $\Fun^{\Rexx}(\FinnSpace, \calD) \rightarrow \calD$. Here $\Fun^{\Rexx}( \FinnSpace, \calD)$ denotes the full subcategory of $\Fun( \FinnSpace, \calD)$ spanned by the right exact functors.
\end{remark}

\begin{lemma}\label{comtil}
\begin{itemize}
\item[$(1)$] Each object of $\FinSpace$ is compact in $\SSet_{\ast}$.
\item[$(2)$] The inclusion $\FinSpace \subseteq \SSet_{\ast}$ induces an equivalence
$\Ind( \FinSpace ) \rightarrow \SSet_{\ast}$. In particular, $\SSet_{\ast}$ is compactly generated.
\item[$(3)$] The subcategory $\FinSpace \subseteq \SSet_{\ast}$ is the smallest full subcategory
which contains $S^0$ and is stable under finite colimits.
\end{itemize}
\end{lemma}

\begin{proof}
Since $\FinnSpace$ consists of compact objects of $\SSet$, Proposition \toposref{accessforwardslice} implies that $\FinSpace$ consists of compact objects of
$\SSet_{\ast}$. This proves $(1)$. 

We next observe that $\FinSpace$ is stable under finite colimits in $\SSet_{\ast}$. Using the proof of Corollary \toposref{allfin}, we may reduce to showing that $\FinSpace$ is stable under pushouts and contains an initial object of $\SSet_{\ast}$. The second assertion is obvious, and the first follows from the fact that the forgetful functor $\SSet_{\ast} \rightarrow \SSet$ commutes with pushouts (Proposition \toposref{goeselse}). 

We now prove $(3)$. Let $\SSet'_{\ast}$ be the smallest full subcategory which contains $S^0$ and is stable under finite colimits. The above argument shows that $\SSet'_{\ast} \subseteq \FinSpace$.
To prove the converse, we let
$f: \SSet \rightarrow \SSet_{\ast}$ be a left adjoint to the forgetful functor, so that $f(X) \simeq X \coprod \ast$. Then $f$ preserves small colimits. Since $f(\ast) \simeq S^0 \in \SSet'_{\ast}$, we conclude that $f$ carries $\FinnSpace$ into $\SSet'_{\ast}$. 
If $x: \ast \rightarrow X$ is a pointed object of $\SSet$, then $x$ can be written as a coproduct $f(X) \coprod_{S^0} \ast$. In particular, if $x \in \FinSpace$, then $X \in \FinnSpace$, so that
$f(X), S^0, \ast \in \SSet'_{\ast}$. Since $\SSet'_{\ast}$ is stable under pushouts, we conclude that
$x \in \SSet'_{\ast}$; this completes the proof of $(3)$.

We now prove $(2)$. Part $(1)$ and Proposition \toposref{uterr} imply that we have a fully faithful functor $\theta: \Ind(\FinSpace) \subseteq \SSet_{\ast}$. Let $\SSet''_{\ast}$ be the essential image of $\theta$. Proposition \toposref{sumatch} implies that $\SSet''_{\ast}$ is stable under small colimits. Since $S^0 \in \SSet''_{\ast}$ and $f$ preserves small colimits, we conclude that
$f(X) \in \SSet''_{\ast}$ for all $X \in \SSet$. Since $\SSet''_{\ast}$ is stable under pushouts, we conclude that $\SSet''_{\ast} = \SSet_{\ast}$, as desired.
\end{proof}

\begin{warning}\label{burner}\index{Wall finiteness obstruction}
The $\infty$-category $\FinnSpace$ does not coincide with the $\infty$-category of compact objects $\SSet^{\omega} \subseteq \SSet$. Instead, there is an inclusion $\FinnSpace \subseteq \SSet^{\omega}$, which realizes $\SSet^{\omega}$ as an idempotent completion of $\FinnSpace$. 
An object of $X \in \SSet^{\omega}$ belongs to $\FinnSpace$ if and only if its {\em Wall finiteness obstruction} vanishes.
\end{warning}

\begin{proposition}\label{specster2}
The $\infty$-category of spectra is compactly generated. Moreover, an object
$X \in \Spectra$ is compact if and only if it is a retract of $\Sigma^{\infty-n} Y$, for
some $Y \in \FinSpace$ and some integer $n$.
\end{proposition}

\begin{proof}
Let $\LLPres{\omega}$, $\RRPres{\omega}$, and $\Cat^{\Rexx}_{\infty}$ be defined as \S \toposref{compactgen}. According to 
Proposition \toposref{sumer}, we can view the construction of $\Ind$-categories as determining a localization functor $\Ind: \Cat_{\infty}^{\Rexx} \rightarrow \LLPres{\omega}$. 
Let $\FinSpectra$ denote the colimit of the sequence
$$ \FinSpace \stackrel{\Sigma}{\rightarrow} \FinSpace \stackrel{\Sigma}{\rightarrow} \ldots $$
in $\Cat^{\Rexx}_{\infty}$. Since $\SSet_{\ast} \simeq \Ind( \FinSpace )$ (Lemma \ref{comtil}) and the functor $\Ind$ preserves colimits, we conclude that
$\Ind(\FinSpectra)$ can be identified with the colimit of the sequence
$$ \SSet_{\ast} \stackrel{\Sigma}{\rightarrow} \SSet_{\ast} \stackrel{\Sigma}{\rightarrow} \ldots $$
in $\LLPres{\omega}$. Invoking the equivalence between $\LLPres{\omega}$ and $(\RRPres{\omega})^{op}$ (see Notation \toposref{funnote}), we can identify $\Ind(\FinSpectra)$ with the {\em limit} of the tower
$$ \ldots \stackrel{\Omega}{\rightarrow} \SSet_{\ast} \stackrel{\Omega}{\rightarrow} \SSet_{\ast}$$
in $\RRPres{\omega}$. Since the inclusion functor $\RRPres{\omega} \subseteq \widehat{\Cat}_{\infty}$ preserves limits (Proposition \toposref{cnote}), we conclude that there is an equivalence $F: \Ind(\FinSpectra) \simeq \Spectra$ (Proposition \ref{camer}). This proves that $\Spectra$ is compactly generated, and that the compact objects of $\Spectra$ are precisely those which appear as retracts of $F(Y)$, for some $Y \in \Ind( \FinSpectra)$. To complete the proof, we observe that
$Y$ itself lies in the image of one of the maps $\FinSpace \rightarrow \FinSpectra$, and that
the composite maps
$$ \FinSpace \rightarrow \FinSpectra \rightarrow \Ind( \FinSpectra) \stackrel{F}{\rightarrow} \Spectra$$
are given by restricting the suspension spectrum functors $\Sigma^{\infty-n}: \FinSpace \rightarrow \Spectra$.
\end{proof}

\begin{remark}
The proof of Proposition \ref{specster2} implies that we can identify $\FinSpectra$ with a full subcategory of the compact objects of $\Spectra$. In fact, every compact object of $\Spectra$ belongs to this full subcategory. The proof of this is not completely formal (especially in view of Warning \ref{burner}); it relies on the fact that the ring of integers $\Z$ is a principal ideal domain, so that every finitely generated projective $\Z$-module is free.
\end{remark}

\begin{remark}
It is possible to use the proof of Proposition \ref{specster2} to prove directly
that the $\infty$-category $\Spectra$ is stable, without appealing to the general
results on stabilization proved in \S \ref{stable9.1}. Indeed, by virtue of 
Proposition \ref{kappstable}, it suffices to show that the $\infty$-category
$\FinSpectra$ of finite spectra is stable. The essence of the matter is now to show that
every pushout square in $\FinSpectra$ is also a pullback square. Every pushout square
is obtained from a pushout diagram
$$ \xymatrix{ W \ar[r] \ar[d] & X \ar[d] \\
Y \ar[r] & Z }$$
in the $\infty$-category $\FinSpace$ of finite pointed spaces. This pushout square will typically not be homotopy Cartesian $\FinSpace$, but will be {\em approximately} homotopy Cartesian
if the spaces involved are highly connected: this follows from the Blakers-Massey homotopy excision theorem (see for example \cite{hatcher}, p. 360). Using the fact that the approximation gets better and better as we iterate the suspension functor $\Sigma$ (which increases the connectivity of spaces), one can deduce that the image of the above square is a pullback in $\FinSpace$. 
\end{remark}

\begin{remark}\label{sunrise}
Let $\Ab$ denote the category of abelian groups. For each $n \in \Z$, we let
$\pi_n: \Spectra \rightarrow \Nerve(\Ab)$ be the composition of the shift
functor $X \mapsto X[-n]$ with the equivalence $\SSet_{\infty}^{\heartsuit} \simeq \Nerve(\Ab)$.
Note that if
$n \geq 2$, then $\pi_n$ can be identified with the composition
$$ \Spectra \stackrel{\Omega^{\infty}_{\ast}}{\rightarrow} \SSet_{\ast} \stackrel{\pi_n}{\rightarrow} \Nerve(\Ab)$$
where the second map is the usual homotopy group functor. Since $\Spectra$ is both left and right complete, we conclude that a map $f: X \rightarrow Y$ of spectra is an equivalence if and only if
it induces isomorphisms $\pi_n X \rightarrow \pi_n Y$ for all $n \in \Z$. 
\end{remark}

\begin{proposition}\label{denkmal}
The functor $\Omega^{\infty}: (\Spectra)_{\geq 0} \rightarrow \SSet$ preserves geometric realizations of simplicial objects. 
\end{proposition}

\begin{proof}
Since the simplicial set $\Nerve(\cDelta)$ is weakly contractible, the forgetful functor
$\SSet_{\ast} \rightarrow \SSet$ preserves geometric realizations of simplicial objects (Proposition \toposref{goeselse}). It will therefore suffice to prove that the functor $\Omega^{\infty}_{\ast} | (\Spectra)_{\geq 0} \rightarrow \SSet_{\ast}$ preserves geometric realizations of simplicial objects.

For each $n \geq 0$, let $\SSet^{\geq n}$ denote the full subcategory of $\SSet$ spanned by the $n$-connective objects, and let $\SSet^{\geq n}_{\ast}$ be the $\infty$-category of pointed objects of $\SSet^{\geq n}$. We observe that $(\Spectra)_{\geq 0}$ can be identified with the homotopy inverse limit of the tower
$$ \ldots \stackrel{\Omega}{\rightarrow} \SSet_{\ast}^{\geq 1} \stackrel{\Omega}{\rightarrow} \SSet_{\ast}^{\geq 0}.$$
It will therefore suffice to prove that for every $n \geq 0$, the loop functor
$\Omega: \SSet_{\ast}^{\geq n+1} \rightarrow \SSet_{\ast}^{\geq n}$ preserves geometric realizations of simplicial objects. 

The $\infty$-category $\SSet^{\geq n}$ is the preimage (under $\tau_{\leq n-1}$) of
the full subcategory of $\tau_{\leq n-1} \SSet$ spanned by the final objects. Since this full subcategory is stable under geometric realizations of simplicial objects and since
$\tau_{\leq n-1}$ commutes with all colimits, we conclude that
$\SSet^{\geq n} \subseteq \SSet$ is stable under geometric realizations of simplicial objects. 

According to Lemmas \toposref{preEM} and \toposref{postEM}, there is an equivalence of
$\SSet_{\ast}^{\geq 1}$ with the $\infty$-category of group objects $\Group(\SSet_{\ast})$. This restricts to an equivalence of $\SSet_{\ast}^{\geq n+1}$ with $\Group( \SSet_{\ast}^{\geq n} )$ for all $n \geq 0$. Moreover, under this equivalence, the loop functor $\Omega$ can be identified with the composition
$$ \Group( \SSet_{\ast}^{\geq n} ) \subseteq \Fun( \Nerve(\cDelta)^{op}, \SSet_{\ast}^{\geq n})
\rightarrow \SSet_{\ast}^{\geq n},$$
where the second map is given by evaluation at the object $[1] \in \cDelta$. This evaluation map commutes with geometric realizations of simplicial objects ( Proposition \toposref{limiteval}). 
Consequently, it will suffice to show that $\Group( \SSet_{\ast}^{\geq n}) \subseteq \Fun( \Nerve(\cDelta)^{op}, \SSet_{\ast}^{\geq n})$ is stable under geometric realizations of simplicial objects. 

Without loss of generality, we may suppose $n = 0$; now we are reduced to showing that
$\Group( \SSet_{\ast}) \subseteq \Fun( \Nerve(\cDelta)^{op}, \SSet_{\ast} )$ is stable under geometric realizations of simplicial objects. In view of Lemma \toposref{postEM}, it will suffice to show that $\Group(\SSet) \subseteq \SSet_{\Delta}$ is stable under geometric realizations of simplicial objects. Invoking Proposition \toposref{tinner}, we are reduced to proving that
the formation of geometric realizations in $\SSet$ commutes with finite products, which follows
from Lemma \toposref{bale2}.
\end{proof}

\section{Excisive Functors}\label{stable9.2}
 



In order to study the relationship between an $\infty$-category $\calC$ and its stabilization $\Stab(\calC)$, we need to introduce a bit of terminology.

\begin{definition}\index{weakly excisive}\index{excisive}\index{functor!excisive}\index{functor!weakly excisive}\label{cripp}
Let $F: \calC \rightarrow \calD$ be a functor between $\infty$-categories.
\begin{itemize}
\item[$(i)$] If $\calC$ has an initial object $\emptyset$, then we will say that $F$ is {\it weakly excisive} if $F(\emptyset)$ is a final object of $\calD$. We let $\Fun_{\ast}(\calC, \calD)$ denote the full subcategory of $\Fun(\calC, \calD)$ spanned by the weakly excisive functors.
\item[$(ii)$] If $\calC$ admits finite colimits, then we will say that $F$ is {\it excisive} if it is weakly excisive, and $F$ carries pushout squares in $\calC$ to pullback squares in $\calD$. We let $\Exc(\calC, \calD)$ denotes the full subcategory of $\Fun(\calC, \calD)$ spanned by the excisive functors.\index{ZZZFunast@$\Fun_{\ast}(\calC, \calD)$}\index{ZZZFunExc@$\Exc(\calC, \calD)$}
\end{itemize}
\end{definition}

\begin{warning}
Definition \ref{cripp} is somewhat nonstandard: most authors do not require the property
the preservation of zero objects in the definition of excisive functors.
\end{warning}

\begin{remark}\label{youg}
Let $F: \calC \rightarrow \calD$ be a functor between $\infty$-categories, and suppose that $\calC$ admits finite colimits. If $\calC$ is stable, then $F$ is excisive if and only if it is left exact (Proposition \ref{surose}). If instead $\calD$ is stable, then $F$ is excisive if and only if it is right exact. In particular, if both $\calC$ and $\calD$ are stable, then $F$ is excisive if and only if it is exact (Proposition \ref{funrose}).
\end{remark}

\begin{lemma}\label{surtit}
Let $\calC$ and $\calD$ be $\infty$-categories, and assume that $\calC$ has an initial object. Then: 
\begin{itemize}
\item[$(1)$] The forgetful functor $\theta: \Fun_{\ast}(\calC, \calD_{\ast}) \rightarrow \Fun_{\ast}(\calC, \calD)$ is a trivial fibration of simplicial sets. 
\item[$(2)$] If $\calC$ admits finite colimits, then the forgetful functor
$\theta': \Exc(\calC, \calD_{\ast}) \rightarrow \Exc(\calC, \calD)$ is a trivial fibration of simplicial sets.
\end{itemize}
\end{lemma}

\begin{remark}
If the $\infty$-category $\calD$ does not have a final object, then the conclusion of Lemma \ref{surtit} is valid, but degenerate: both of the relevant $\infty$-categories of functors are empty.
\end{remark}

\begin{proof}
To prove $(1)$, we first
observe that objects of $\Fun_{\ast}(\calC, \calD_{\ast})$ can be identified with
maps $F: \calC \times \Delta^1 \rightarrow \calD$ with the following properties:
\begin{itemize}
\item[$(a)$] For every initial object $C \in \calC$, $F(C,1)$ is a final object of $\calD$.
\item[$(b)$] For every object $C \in \calC$, $F(C,0)$ is a final object of $\calD$. 
\end{itemize}
Assume for the moment that $(a)$ is satisfied, and let $\calC' \subseteq \calC \times \Delta^1$ be the full subcategory spanned by those objects $(C,i)$ for which either $i = 1$, or $C$ is an initial object of $\calC$. We observe that $(b)$ is equivalent to the following pair of conditions:
\begin{itemize}
\item[$(b')$] The functor $F | \calC'$ is a right Kan extension of $F| \calC \times \{1\}$.
\item[$(b'')$] The functor $F$ is a left Kan extension of $F| \calC'$. 
\end{itemize}
Let $\calE$ be the full subcategory of $\Fun( \calC \times \Delta^1, \calD)$ spanned by those functors which satisfy conditions $(b')$ and $(b'')$. Using Proposition \toposref{lklk}, we deduce that
the projection $\overline{\theta}: \calE \rightarrow \Fun( \calC \times \{1\}, \calD)$ is a trivial Kan fibration. Since $\theta$ is a pullback of $\overline{\theta}$, we conclude that $\theta$ is a trivial Kan fibration. This completes the proof of $(1)$.

To prove $(2)$, we observe that $\theta'$ is a pullback of $\theta$ (since Proposition \toposref{needed17} asserts that a square in
$\calD_{\ast}$ is a pullback if and only if the underlying square in $\calD$ is a pullback).
\end{proof}

\begin{remark}\label{neeko}
Let $\calC$ be a pointed $\infty$-category which admits finite colimits, and $\calD$ a pointed $\infty$-category which admits finite limits. Let $F: \Fun( \calC, \calD) \rightarrow \Fun(\calC, \calD)$ be given by composition with the suspension functor $\calC \rightarrow \calC$, and let
$G: \Fun( \calC, \calD) \rightarrow \Fun( \calC, \calD)$ be given by composition with the loop functor $\Omega: \calD \rightarrow \calD$. Then $F$ and $G$ restrict to give homotopy inverse equivalences $$ \Adjoint{F}{\Exc(\calC,\calD)}{\Exc(\calC,\calD)}{G}.$$
\end{remark}

\begin{notation}\label{unwind}
Let $F: \calC \rightarrow \calD$ be a functor between $\infty$-categories, 
and assume that $\calD$ admits finite limits. For every commutative square $\tau$:
$$ \xymatrix{ W \ar[r] \ar[d] & X \ar[d] \\
Y \ar[r] & Z }$$
in $\calC$, we obtain a commutative square $F(\tau)$:
$$ \xymatrix{ F(W) \ar[r] \ar[d] & F(X) \ar[d] \\
F(Y) \ar[r] & F(Z) }$$ in $\calD$. This diagram determines a map
$\eta_{\tau}: F(W) \rightarrow F(X) \times_{ F(Z) } F(Y)$ in the $\infty$-category $\calD$, which is well-defined up to homotopy. If we suppose further that $Y$ and $Z$ are zero objects of $\calC$,
that $F(Y)$ and $F(Z)$ are zero objects of $\calD$, and that $\tau$ is a pushout diagram, then
we obtain a map $F(W) \rightarrow \Omega F( \Sigma W)$, which we will denote simply by $\eta_{W}$.
\end{notation}

\begin{proposition}\label{coople}
Let $\calC$ be a pointed $\infty$-category which admits finite colimits, $\calD$ a pointed
$\infty$-category which admits finite limits, and let $F: \calC \rightarrow \calD$ be a functor which carries zero objects of $\calC$ to zero objects of $\calD$.
The following conditions are equivalent:
\begin{itemize}
\item[$(1)$] The functor $F$ is excisive: that is, $F$ carries pushout squares in $\calC$ to pullback squares in $\calD$.

\item[$(2)$] For every object $X \in \calC$, the canonical map $\eta_{X}: F(X) \rightarrow \Omega F(\Sigma X)$ is an equivalence in $\calD$ $($see Notation \ref{unwind}$)$.
\end{itemize}
\end{proposition}

\begin{corollary}\label{cuuple}
Let $F: \calC \rightarrow \calD$ be a functor between stable $\infty$-categories. Then
$F$ is exact if and only if the following conditions are satisfied:
\begin{itemize}
\item[$(1)$] The functor $F$ carries zero objects of $\calC$ to zero objects of $\calD$.
\item[$(2)$] For every object $X \in \calC$, the canonical map
$\Sigma F(X) \rightarrow F( \Sigma X)$ is an equivalence in $\calD$.
\end{itemize}
\end{corollary}

\begin{corollary}\label{taskwise}
Let $\calC$ be a pointed $\infty$-category which admits finite limits and colimits.
Then:
\begin{itemize}
\item[$(1)$] If the suspension functor $\Sigma_{\calC}$ is fully faithful, then
every pushout square in $\calC$ is a pullback square.
\item[$(2)$] If the loop functor $\Omega_{\calC}$ is fully faithful, then every
pullback square in $\calC$ is a pushout square.
\item[$(3)$] If the loop functor $\Omega_{\calC}$ is an equivalence of
$\infty$-categories, then $\calC$ is stable.
\end{itemize}
\end{corollary}

\begin{proof}
Assertion $(1)$ follows by applying Proposition \ref{coople} to the identity functor
$\id_{\calC}$, and assertion $(2)$ follows from $(1)$ by passing to the opposite $\infty$-category.
Assertion $(3)$ is an immediate consequence of $(1)$ and $(2)$.
\end{proof}

The proof of Proposition \ref{coople} makes use of the following lemma:

\begin{lemma}\label{goodwillie}
Let $\calC$ be a pointed $\infty$-category which admits finite colimits, $\calD$ a pointed
$\infty$-category which admits finite limits, and $F: \calC \rightarrow \calD$ a functor which
carries zero objects of $\calC$ to zero objects of $\calD$. Suppose given a pushout diagram $\tau$:
$$ \xymatrix{ W \ar[r] \ar[d] & X \ar[d] \\
Y \ar[r] & Z }$$
in $\calC$. Then there exists a map
$\theta_{\tau}: F(X) \times_{ F(Z) } F(Y) \rightarrow \Omega F( \Sigma W)$
with the following properties:
\begin{itemize}
\item[$(1)$] The composition $\theta_{\tau} \circ \eta_{\tau}$ is homotopic to 
$\eta_{W}$. Here $\eta_{\tau}$ and $\eta_{W}$ are defined as in Notation \ref{unwind}.

\item[$(2)$] Let $\Sigma(\tau)$ denote the induced diagram
$$ \xymatrix{ \Sigma W \ar[r] \ar[d] & \Sigma X \ar[d] \\
\Sigma Y \ar[r] & \Sigma Z.}$$
Then there is a pullback square
$$ \xymatrix{ \eta_{ \Sigma(\tau)} \circ \theta_{\tau} \ar[r] \ar[d] & \eta_{ X } \ar[d] \\
\eta_{ Y } \ar[r] & \eta_{ Z } }$$
in the $\infty$-category $\Fun( \Delta^1, \calD)$ of morphisms in $\calD$.
\end{itemize}
\end{lemma}

\begin{proof}
In the $\infty$-category $\calC$, we have the following commutative diagram (in which every square is a pushout):
$$ \xymatrix{ W \ar[r] \ar[d] & X \ar[r] \ar[d] & 0 \ar[d] & \\
Y \ar[r] \ar[d] & X \coprod_{W} Y \ar[d] \ar[r] & 0 \coprod_{W} Y \ar[r] \ar[d] & 0 \ar[d] \\
0 \ar[r] & X \coprod_{W} 0 \ar[r] \ar[d] & \Sigma W \ar[r] \ar[d] & \Sigma Y \ar[d] \\
& 0 \ar[r] & \Sigma X \ar[r] & \Sigma( X \coprod_{W} Y). }$$
Applying the functor $F$, and replacing the upper left square by a pullback, we obtain a new diagram
$$ \xymatrix{ F(X) \times_{ F(Z)} F(Y) \ar[r] \ar[d] & F(X) \ar[r] \ar[d] & 0 \ar[d] & \\
F(Y) \ar[d]  \ar[r] & F(Z) \ar[d] \ar[r] & F(0 \coprod_{W} Y) \ar[r] \ar[d] & 0 \ar[d] \\
0 \ar[r] & F(X \coprod_{W} 0) \ar[r] \ar[d] & F( \Sigma(W) ) \ar[r] \ar[d] & F( \Sigma Y) \ar[d] \\
& 0 \ar[r] & F( \Sigma X ) \ar[r] & F( \Sigma Z). }$$
Restricting attention to the large square in the upper left, we obtain the desired map
$\theta_{\tau}: F(X) \times_{ F(Z)} F(Y) \rightarrow \Omega F( \Sigma W)$.
It is easy to verify that $\theta_{\tau}$ has the desired properties.
\end{proof}

\begin{proof}[Proof of Proposition \ref{coople}]
The implication $(1) \Rightarrow (2)$ is obvious. Conversely, suppose that $(2)$ is satisfied.
We must show that for every pushout square $\tau:$
$$ \xymatrix{ X \ar[r] \ar[d] & Y \ar[d] \\
Z \ar[r] & Y \coprod_{X} Z  }$$
in the $\infty$-category $\calC$, the induced map $\eta_{\tau}$ is an equivalence in $\calD$.
Let $\theta_{\tau}$ be as in the statement of Lemma \ref{goodwillie}. Then
$\theta_{\tau} \circ \eta_{\tau}$ is homotopic to $\eta_{X}$, and is therefore an equivalence
(in virtue of assumption $(2)$). It will therefore suffice to show that $\eta_{\tau}$ is an equivalence.
The preceding argument shows that $\theta_{\tau}$ has a right homotopy inverse. To show that $\theta_{\tau}$ admits a left homotopy inverse, it will suffice to show that $\eta_{ \Sigma \tau} \circ \theta_{\tau}$ is an equivalence. This follows from the second assertion of Lemma \ref{goodwillie}, since the maps
$\eta_{Y}$, $\eta_{Z}$, and $\eta_{Y \coprod_{X} Z}$ are equivalences (by assumption $(2)$, again).
\end{proof}

Let $\calC$ be a (small) pointed $\infty$-category. Let
$\calP_{\ast}(\calC)$ be defined as in Lemma \ref{kija}. Lemma \ref{surtit} implies that the canonical map $\Fun_{\ast}( \calC^{op}, \SSet_{\ast}) \rightarrow \calP_{\ast}(\calC)$ is a trivial fibration. Consequently, the Yoneda embedding lifts to a fully faithful functor $j': \calC \rightarrow \Fun_{\ast}(\calC^{op}, \SSet_{\ast})$, which we will refer to as the {\it pointed Yoneda embedding}. Our terminology is slightly abusive: the functor $j'$ is only well-defined up to a contractible space of choices; we will ignore this ambiguity.

\begin{proposition}\label{urtusk21}
Let $\calC$ be a pointed $\infty$-category which admits finite colimits and $\calD$ an
$\infty$-category which admits finite limits. Then composition with
the canonical map $\Stab(\calD) \rightarrow \calD$ induces an equivalence of $\infty$-categories
$$\theta: \Exc( \calC, \Stab(\calD) ) \rightarrow \Exc(\calC, \calD).$$
\end{proposition}

\begin{proof}
Since the loop functor $\Omega_{\calD}: \calD \rightarrow \calD$ is left exact, the domain of $\theta$ can be identified with a homotopy limit of the tower
$$ \ldots \stackrel{ \circ \Omega_{\calD}}{\rightarrow} \Exc( \calC, \calD_{\ast} ) \stackrel{ \circ \Omega_{\calD}}{\rightarrow} \Exc(\calC, \calD_{\ast}).$$ 
Remark \ref{neeko} implies that this tower is essentially constant. Consequently, it will suffice
to show that the canonical map $\Exc(\calC, \calD_{\ast}) \rightarrow \Exc(\calC, \calD)$ is a trivial fibration of simplicial sets, which follows from Lemma \ref{surtit}.
\end{proof}

\begin{example}\label{expire}
Let $\calC$ be an $\infty$-category which admits finite limits, and $K$ an arbitrary simplicial set. Then $\Fun(K, \calC)$ admits finite limits (Proposition \toposref{limiteval}). We have a canonical isomorphism $\Fun(K, \calC)_{\ast} \simeq \Fun(K, \calC_{\ast})$, and the loop functor
on $\Fun(K,\calC)_{\ast}$ can be identified with the functor given by composition with
$\Omega: \calC_{\ast} \rightarrow \calC_{\ast}$. It follows that there is a canonical equivalence of $\infty$-categories
$$ \Stab( \Fun(K, \calC) ) \simeq \Fun( K, \Stab(\calC) ).$$
In particular, $\Stab( \calP(K) )$ can be identified with $\Fun(K, \Spectra)$. 
\end{example}

We can apply Proposition \ref{urtusk21} to give another description of the $\infty$-category $\Stab(\calC)$. 

\begin{lemma}\label{prancer}
Let $\calC$ be an $\infty$-category which admits finite colimits, let $f: \calC \rightarrow \calC_{\ast}$ be a left adjoint to the forgetful functor, and let $\calD$ be a stable $\infty$-category. Then
composition with $f$ induces an equivalence of $\infty$-categories
$\phi: \Exc( \calC_{\ast}, \calD) \rightarrow \Exc(\calC, \calD)$.
\end{lemma}

\begin{proof}
Consider the composition
$$ \theta: \Fun(\calC,\calD) \times \calC_{\ast} \subseteq \Fun(\calC,\calD) \times \Fun(\Delta^1, \calC) \rightarrow \Fun(\Delta^1, \calD) \stackrel{\coker}{\rightarrow} \calD.$$
We can identify $\theta$ with a map $\Fun(\calC, \calD) \rightarrow \Fun( \calC_{\ast}, \calD)$.
Since the collection of pullback squares in $\calD$ is a stable subcategory of $\Fun(\Delta^1 \times \Delta^1, \calD)$, we conclude $\theta$ restricts to a map
$\psi: \Exc(\calC, \calD) \rightarrow \Exc(\calC_{\ast}, \calD)$. It is not difficult to verify that
$\psi$ is a homotopy inverse to $\phi$.
\end{proof}

\begin{proposition}\label{surina}
Let $\calC$ be an $\infty$-category which admits finite colimits, and let $\calD$ be an $\infty$-category which admits finite limits. Then there is a canonical isomorphism
$\Exc( \calC_{\ast}, \calD) \simeq \Exc( \calC, \Stab(\calD))$ in the homotopy category of $\infty$-categories.
\end{proposition}

\begin{proof}
Combining Lemma \ref{prancer} and Proposition \ref{urtusk21}, we obtain
a diagram of equivalences
$$ \Exc( \calC_{\ast}, \calD) \leftarrow \Exc(\calC_{\ast}, \Stab(\calD)) \rightarrow \Exc(\calC, \Stab(\calD)).$$
\end{proof}

\begin{corollary}\label{surritt}
Let $\calD$ be an $\infty$-category which admits finite limits. Then there is a canonical
equivalence $\Stab(\calD) \simeq \Exc( \FinSpace, \calD)$ in the homotopy category of $\infty$-categories.
\end{corollary}

\begin{proof}
Combine Proposition \ref{surina}, Remark \ref{yogause}, and Remark \ref{youg}. 
\end{proof}

\begin{corollary}\label{tintrusj}
The $\infty$-category of spectra is equivalent to the $\infty$-category
$\Exc( \SSet^{\fin}_{\ast}, \SSet)$.
\end{corollary}

\begin{remark}\label{surrit}\index{homology theory}
Corollary \ref{tintrusj} provides a very explicit model for spectra. Namely, we can identify a spectrum with an excisive functor $F: \FinSpace \rightarrow \SSet$. We should think of $F$ as a {\it homology theory} $A$. More precisely, given a pair of finite spaces $X_0 \subseteq X$, we can define the relative homology group $A_n( X,X_0)$ to be $\pi_n F( X/X_0)$, where $X/X_0$ denotes the pointed space obtained from $X$ by collapsing $X_0$ to a point (here the homotopy group is taken with base point provided by the map $\ast \simeq F(\ast) \rightarrow F(X/X_0)$ ). 
The assumption that $F$ is excisive is precisely what is needed to guarantee the existence of the usual excision exact sequences for the homology theory $A$. 
\end{remark}
 
\section{Filtered Objects and Spectral Sequences}\label{filttt}

Suppose given a sequence of objects
$$ \ldots \rightarrow X(-1) \stackrel{f^0}{\rightarrow} X(0) \stackrel{f^1}{\rightarrow} X(1) \rightarrow \ldots.$$ 
in a stable $\infty$-category $\calC$. Suppose further that $\calC$ is equipped with a t-structure, and that the heart of $\calC$ is equivalent to the nerve of an abelian category $\calA$. In this section, we will construct a spectral sequence taking values in the abelian category $\calA$, with the $E_1$-page described by the formula
$$ E_{1}^{p,q} = \pi_{p+q} \coker( f^p ) \in \calA.$$
Under appropriate hypotheses, we will see that this spectral sequence converges to the homotopy groups of the colimit $\varinjlim( X(i) )$. 

Our first step is to construct some auxiliary objects in $\calC$. 

\begin{definition}\label{sumptuous}\index{$\calI$-complex}\index{complex}
Let $\calC$ be a pointed $\infty$-category, and let $\calI$ be a linearly ordered set. We let
$\calI^{[1]}$ denote the partially ordered set of pairs of elements $i \leq j$ of $\calI$, where
$(i,j) \leq (i',j')$ if $i \leq j$ and $i' \leq j'$. An {\it $\calI$-complex} in $\calC$
is a functor $F: \Nerve( \calI^{[1]}) \rightarrow \calC$ with the following properties:

\begin{itemize}
\item[$(1)$] For each $i \in \calI$, $F(i,i)$ is a zero object of $\calC$. 

\item[$(2)$] For every $i \leq j \leq k$, the associated diagram
$$ \xymatrix{ F(i,j) \ar[r] \ar[d] & F(i,k) \ar[d] \\
F(j,j) \ar[r] & F(j,k)}$$ is a pushout square in $\calC$.
\end{itemize}

We let $\Gap( \calI, \calC)$ denote the full subcategory of $\Fun( \Nerve( \calI^{[1]} ), \calC)$
spanned by the $\calI$-complexes in $\calC$.\index{ZZZGap@$\Gap(\calI, \calC)$}
\end{definition}

\begin{remark}\label{makeplex}
Let $F \in \Gap( \Z, \calC)$ be a $\Z$-complex in a stable $\infty$-category $\calC$. For each $n \in \Z$, the functor $F$ determines pushout square
$$ \xymatrix{ F(n-1, n) \ar[r] \ar[d] & F(n-1, n+1) \ar[d] \\
0 \ar[r] & F(n, n+1), }$$
hence a boundary $\delta: F(n, n+1) \rightarrow F(n-1, n)[1]$. 
If we set $C_{n} = F(n-1,n)[-n]$, then we obtain a sequence of maps
$$ \ldots \rightarrow C_{1} \stackrel{d_1}{\rightarrow} C_{0} \stackrel{d_0}{\rightarrow} C_{-1} \rightarrow \ldots$$
in the homotopy category $\h{\calC}$. The commutative diagram
$$ \xymatrix{ F(n,n+1) \ar[r]^{\delta} & F(n-2,n)[1] \ar[d] \ar[r] & F(n-1,n)[1] \ar[d]^{\delta} \\
0 \ar[r]^-{\sim} & F(n-1,n-1)[1] \ar[r] & F(n-2, n-1)[2]. }$$ 
shows that $d_{n-1} \circ d_{n} \simeq 0$, so that $( C_{\bigdot}, d_{\bigdot} )$ can be viewed
as a chain complex in the triangulated category $\h{\calC}$. This motivates the terminology of Definition \ref{sumptuous}. 
\end{remark}

\begin{lemma}\label{loij}
Let $\calC$ be a pointed $\infty$-category which admits pushouts.
Let $\calI = \calI_0 \cup \{ - \infty \}$ be a linearly ordered set containing a least element
$- \infty$. We regard $\calI_0$ as a linearly ordered subset of $\calI^{[1]}$ via the
embedding
$$ i \mapsto (-\infty, i).$$
Then the restriction map $\Gap( \calI, \calC) \rightarrow \Fun( \Nerve(\calI_0), \calC)$ is an equivalence of $\infty$-categories.
\end{lemma}

\begin{proof}
Let $\calJ = \{ (i,j) \in \calI^{[1]} : ( i = - \infty) \vee (i=j) \}$. We now make the following observations:

\begin{itemize}
\item[$(1)$] A functor $F: \Nerve(\calI^{[1]}) \rightarrow \calC$ is a complex if and only if
$F$ is a left Kan extension of $F | \Nerve(\calJ)$, and $F(i,i)$ is a zero object of $\calC$ for
all $i \in \calI$.
\item[$(2)$] Any functor $F_0: \Nerve(\calJ) \rightarrow \calC$ admits a left Kan extension
to $\Nerve( \calI^{[1]})$ (use Lemma \toposref{kan2} and the fact that $\calC$ admits pushouts).
\item[$(3)$] A functor $F_0: \Nerve(\calJ) \rightarrow \calC$ has the property
that $F_0(i,i)$ is a zero object, for every $i \in \calI$, if and only if 
$F_0$ is a right Kan extension of $F_0 | \Nerve( \calI_0 )$. 
\item[$(4)$] Any functor $F_0: \Nerve(\calI_0) \rightarrow \calC$ admits a right
Kan extension to $\Nerve(\calJ)$ (use Lemma \toposref{kan2} and the fact that $\calC$ has a final object). 
\end{itemize}
The desired conclusion now follows immediately from Proposition \toposref{lklk}. 
\end{proof}

\begin{remark}\index{K-theory!of an $\infty$-category}\index{Waldhausen K-theory}
Let $\calC$ be a pointed $\infty$-category which admits pushouts (for example, a stable $\infty$-category). For each $n \geq 0$, let $\Gap^{0}( [n], \calC)$ be the largest Kan complex contained in 
$\Gap( [n], \calC)$. Then the assignment
$$[n] \mapsto \Gap( [n], \calC)$$
determines a simplicial object in the category $\Kan$ of Kan complexes. We can then define
the {\it Waldhausen $K$-theory of $\calC$} to be a geometric realization of this bisimplicial set (for example, the associated diagonal simplicial set). In the special case where $A$ is an $A_{\infty}$-ring and $\calC$ is the smallest stable subcategory of $\Mod_{A}$ which contains $A$, this definition recovers the usual $K$-theory of $A$. We refer the reader to \cite{toenK} for a related construction.
\end{remark}

\begin{construction}\label{sseq}\index{spectral sequence}\index{ZZZEpqr@$E^{p,q}_{r}$}
Let $\calC$ be a stable $\infty$-category equipped with a t-structure, such that the
heart of $\calC$ is equivalent to the nerve of an abelian category $\calA$. 
Let $X \in \Gap( \Z, \calC)$. We observe that for every triple of integers
$i \leq j \leq k$, there is a long exact sequence
$$ \ldots \rightarrow \pi_{n} X(i,j) \rightarrow \pi_{n} X(i,k)
\rightarrow \pi_{n} X(j,k) \stackrel{ \delta }{\rightarrow} \pi_{n-1} X(i,j) \rightarrow \ldots $$
in the abelian category $\calA$. For every $p,q \in \Z$ and every $r \geq 1$, we define
the object $E^{p,q}_{r} \in \calA$ by the formula
$$ E^{p,q}_{r} = \im( \pi_{p+q} X(p-r, p) \rightarrow \pi_{p+q} X(p-1, p+r-1)).$$ 
There is a differential $d_{r}: E^{p,q}_{r} \rightarrow E^{p - r, q + r -1}_{r}$, uniquely determined by the requirement that the diagram
$$ \xymatrix{ \pi_{p+q} X(p-r,p) \ar[r] \ar[d]^{\delta} & E^{p,q}_{r} \ar[r] \ar[d]^{d_r} &
\pi_{p+q} X(p-1, p+r-1) \ar[d]^{\delta} \\
\pi_{p+q-1} X(p-2r, p-r) \ar[r] & E^{p-r, q+r-1}_{r} \ar[r] & \pi_{p+q-1} X(p-r-1, p-1) }$$ 
be commutative.
\end{construction}

\begin{proposition}\label{sest}
Let $X \in \Gap( \Z, \calC)$ be as in Construction \ref{sseq}. Then:
\begin{itemize}
\item[$(1)$] For each $r \geq 1$, the composition $d_r \circ d_r$ is zero.
\item[$(2)$] There are canonical isomorphisms
$$ E^{p,q}_{r+1} \simeq \ker( d_r: E^{p, q}_{r} \rightarrow E^{p-r, q+r-1}_{r} )
/ \im( d_r: E^{p+r, q-r+1}_{r} \rightarrow E^{p,q}_{r}).$$
\end{itemize}
Consequently, $\{ E^{p,q}_{r}, d_r \}$ is a spectral sequence $($with values in the abelian category
$\calA$ $)$.
\end{proposition}

\begin{remark}\label{saltine2}
For fixed $q \in \Z$, the complex $(E^{\bigdot, q}_{1}, d_1)$ in $\calA$ can be obtained from
the $\h{\calC}$-valued chain complex $C_{\bigdot}$ described in Remark \ref{makeplex} by applying the cohomological functor $\pi_{q}$. 
\end{remark}

\begin{proof}
We have a commutative diagram 
$$ \xymatrix{ & \pi_{p+q} X(p-r-1,p) \ar[d] & \\
 \pi_{p+q+1} X(p,p+r) \ar[r]^{\delta} \ar[d] &  \pi_{p+q} X(p-r,p) \ar[d] \ar[r]^{\delta} 
&   \pi_{p+q-1} X(p-2r, p-r) \ar[d] \\
E_{r}^{p+r, q-r+1} \ar[r]^{d_r} \ar[d] & E^{p,q}_{r} \ar[d] \ar[r]^{d_r} & E^{p-r, q+r-1}_{r} \ar[d] \\
\pi_{p+q+1} X(p+r-1, p+2r-1) \ar[r]^{\delta} & \pi_{p+q} X(p-1, p+r-1) \ar[r]^{\delta} \ar[d] &
\pi_{p+q-1} X(p-r-1, p-1) \\
&  \pi_{p+q} X(p-1, p+r). & }$$
Since the upper left vertical map is an epimorphism, $(1)$ will follow provided that we can show that the composition 
$$ \pi_{p+q+1} X(p,p+r) \stackrel{\delta}{\rightarrow} \pi_{p+q} X(p-r,p)
\stackrel{\delta}{\rightarrow} \pi_{p+q-1} X(p-2r, p-r)$$
is zero. This follows immediately from the commutativity of the diagram 
$$ \xymatrix{ X(p,p+r) \ar[r]^{\delta} & X(p-2r,p)[1] \ar[d] \ar[r] & X(p-r,p)[1] \ar[d]^{\delta} \\
0 \ar[r]^-{\sim} & X(p-r,p-r)[1] \ar[r] & X(p-2r, p-r)[2]. }$$ 

We next claim that the composite map
$$ \phi: \pi_{p+q} X(p-r-1, p) \rightarrow E^{p,q}_{r} \stackrel{d_r}{\rightarrow} E^{p-r, q+r-1}_{r}$$
is zero. Because $E^{p-r, q+r-1}_{r} \rightarrow \pi_{p+q-1} X(p-r-1,p-1)$ is a monomorphism, this follows from the commutativity of the diagram
$$ \xymatrix{ \pi_{p+q} X(p-r-1,p) \ar[d] \ar[r] & \pi_{p+q-1} X(p-2r, p-r-1) \ar[d] \\
 \pi_{p+q} X(p-r,p) \ar[d] \ar[r] &  \pi_{p+q-1} X(p-2r, p-r) \ar[d] \\
 E^{p,q}_{r} \ar[r]^{d_r} & E^{p-r, q+r-1}_{r} \ar[d] \\
 & \pi_{p+q-1} X(p-r-1, p-1), }$$
since the composition of the left vertical line factors through $\pi_{p+q-1} X(p-r-1,p-r-1) \simeq 0$. 
A dual argument shows that the composition
$$ E^{p+r, q+r-1}_{r} \stackrel{d_r}{\rightarrow} E^{p,q}_{r} \rightarrow
\pi_{p+q} X(p-1, p+r)$$
is zero as well.

Let $Z = \ker( d_r: E^{p,q}_{r} \rightarrow E^{p-r,q+r-1}_{r})$ and
$B = \im( d_r: E^{p+r, q-r+1}_{r} \rightarrow E^{p,q}_{r} )$. The above arguments yield a sequence of morphisms
$$ \pi_{p+q} X(p-r-1,p) \stackrel{\phi}{\rightarrow} 
Z \stackrel{\phi'}{\rightarrow} Z/B \stackrel{\psi'}{\rightarrow}
E^{p,q}_{r} / B \stackrel{\psi}{\rightarrow} \pi_{p+q} X(p-1, p+r).$$
To complete the proof of $(2)$, it will suffice to show that
$\phi' \circ \phi$ is an epimorphism and that $\psi \circ \psi'$ is a monomorphism.
By symmetry, it will suffice to prove the first assertion. Since $\phi'$ is evidently an epimorphism, we are reduced to showing that $\phi$ is an epimorphism. 

Let $K$ denote the kernel of the composite map
$$ \pi_{p+q} X(p-r,p) \rightarrow E^{p,q}_{r} \stackrel{d_r}{\rightarrow} E^{p-r, q+r-1}_{r}
\rightarrow \pi_{p+q-1} X(p-r-1,p-1),$$
so that the canonical map $K \rightarrow Z$ is an epimorphism. Choose a diagram
$$ \xymatrix{ \widetilde{K} \ar[r]^{f} \ar[d]^{g} & K \ar[d] \ar[r] & 0 \ar[d] \\
\pi_{p+q} X(p-r, p-1) \ar[r] & \pi_{p+q-1} X(p-r-1, p-r) \ar[r] & \pi_{p+q-1} X(p-r-1, p-1) }$$
where the square on the left is a pullback. The exactness of the bottom row implies that
$f$ is an epimorphism. Let $f'$ denote the composition
$$ \widetilde{K} \stackrel{g}{\rightarrow} \pi_{p+q} X(p-r,p-1) \rightarrow
\pi_{p+q} X(p-r,p).$$
The composition
$$ \widetilde{K} \stackrel{f'}{\rightarrow} \pi_{p+q} X(p-r,p) \rightarrow \pi_{p+q} X(p-1, p+r)$$
factors through $\pi_{p+q} X(p-1,p-1) \simeq 0$. Since
$E^{p,q}_{r} \rightarrow \pi_{p+q} X(p-1, p+r)$ is a monomorphism, we conclude that the composition $\widetilde{K} \stackrel{f'}{\rightarrow} \pi_{p+q} X(p-r, p) \rightarrow E^{p,q}_{r}$ is the zero map.
It follows that the composition
$$ \widetilde{K} \stackrel{f-f'}{\rightarrow} \pi_{p+q} X(p-r,p) \rightarrow Z$$
coincides with the composition $\widetilde{K} \stackrel{f}{\rightarrow} K \rightarrow Z$, and is therefore an epimorphism. 

Form a diagram
$$ \xymatrix{ \overline{K} \ar[r]^{f''} \ar[d] & \widetilde{K} \ar[d]^{f-f'} \ar[r] & 0 \ar[d] \\
\pi_{p+q} X(p-r-1, p) \ar[r] & \pi_{p+q} X(p-r,p) \ar[r] & \pi_{n-1} X(p-r-1, p-r) }$$
where the left square is a pullback. Since the bottom line is exact, we conclude that $f''$
is an epimorphism, so that the composition
$$ \overline{K} \stackrel{f''}{\rightarrow} \widetilde{K} \stackrel{f-f'}{\rightarrow}
\pi_{p+q} X(p-r, p) \rightarrow Z$$
is an epimorphism. This map coincides with the composition
$$ \overline{K} \rightarrow \pi_{p+q} X(p-r-1,p) \stackrel{\phi}{\rightarrow} Z,$$
so that $\phi$ is an epimorphism as well.
\end{proof}

\begin{definition}\label{sumtuo}\index{filtered object!of a stable $\infty$-category}
Let $\calC$ be a stable $\infty$-category. 
A {\it filtered object} of $\calC$ is a functor
$X: \Nerve(\Z) \rightarrow \calC$.

Suppose that $\calC$ is equipped with a t-structure, and let
$X: \Nerve(\Z) \rightarrow \calC$ be a filtered object of $\calC$. According to Lemma \ref{loij}, we
can extend $X$ to a complex in $\Gap( \Z \cup \{ - \infty \}, \calC)$. Let
$\overline{X}$ be the associated object of $\Gap( \Z, \calC)$, and let
$\{ E^{p,q}_{r}, d_r \}_{r \geq 1}$ be the spectral sequence described in
Construction \ref{sseq} and Proposition \ref{sest}. We will refer to
$\{ E^{p,q}_{r}, d_r \}_{r \geq 1}$ as the {\it spectral sequence associated to the filtered object $X$}.
\end{definition}

\begin{remark}
In the situation of Definition \ref{sumtuo}, Lemma \ref{loij} implies that $\overline{X}$ is determined up to contractible ambiguity by $X$. It follows that the spectral sequence
$\{ E^{p,q}_{r}, d_r \}_{r \geq 1}$ is independent of the choice of $\overline{X}$, up to canonical isomorphism.
\end{remark}

\begin{example}\index{filtered derived category}\index{derived category!filtered}
Let $\calA$ be a sufficiently nice abelian category, and let $\calC$ be the derived $\infty$-category of $\calA$ (see \S \stableref{stable10}). Let
$\Fun( \Nerve(\Z), \calC)$ be the $\infty$-category of filtered objects of $\calC$. Then the homotopy category $\h{ \Fun( \Nerve(\Z), \calC) }$ can be identified with the classical {\it filtered derived category} of $\calA$, obtained from the category of filtered complexes of objects of $\calA$ by inverting all filtered quasi-isomorphisms. In this case, Definition \ref{sumtuo} recovers the usual spectral sequence associated to a filtered complex. 
\end{example}

Our next goal is to establish the convergence of the spectral sequence of Definition \ref{sumtuo}. 
We will treat only the simplest case, which will be sufficient for our applications.

\begin{definition}\index{colimit!sequential}\index{sequential colimit}
Let $\calC$ be an $\infty$-category. We will say that $\calC$ {\it admits sequential colimits}
if every diagram $\Nerve( \Z_{\geq 0}) \rightarrow \calC$ has a colimit in $\calC$.

If $\calC$ is stable and admits sequential colimits, we will say that a t-structure on $\calC$ is {\it
compatible with sequential colimits} if the full subcategory $\calC_{\leq 0}$ is
stable under the colimits of diagrams indexed by $\Nerve( \Z_{\geq 0})$.
\end{definition}

\begin{remark}\label{seqcolim}
Let $\calC$ be a stable $\infty$-category equipped with a t-structure, so that
the heart of $\calC$ is equivalent to (the nerve of) an abelian category $\calA$.
Suppose that $\calC$ admits sequential colimits. Then
$\calC_{\geq 0}$ admits sequential colimits, so that
$\Nerve(\calA)$, being a localization of $\calC_{\geq 0}$, also admits sequential colimits.
If the t-structure on $\calC$ is compatible with sequential colimits, then the inclusion
$\Nerve(\calA) \subseteq \calC$ and the homological functors
$\{ \pi_{n}: \calC \rightarrow \Nerve(\calA) \}_{n \in \Z}$ preserve sequential colimits.
It follows that 
sequential colimits in the abelian category $\calA$ are exact: in other words, the direct limit of a sequence of monomorphisms in $\calA$ is again a monomorphism.
\end{remark}

\begin{proposition}\label{conseq}\index{spectral sequence!convergence of}
Let $\calC$ be a stable $\infty$-category equipped with a t-structure, and
let $X: \Nerve(\Z) \rightarrow \calC$ be a filtered object of $\calC$.
Assume that $\calC$ admits sequential colimits, and that the t-structure on
$\calC$ is compatible with sequential colimits. Suppose furthermore that $X(n) \simeq 0$ for $n \ll 0$. Then the associated spectral sequence $($Definition \ref{sumtuo}$)$ converges
$$ E^{p,q}_{r} \Rightarrow \pi_{p+q} \varinjlim(X).$$
\end{proposition}

\begin{proof}
Let $\calA$ be an abelian category such that the heart of $\calC$ is equivalent to (the nerve of)
$\calA$. 
The convergence assertion of the Proposition has the following meaning:
\begin{itemize}
\item[$(i)$] For fixed $p$ and $q$, the differentials
$d_r: E^{p,q}_{r} \rightarrow E^{p-r, q+r-1}_{r}$ vanish for $r \gg 0$. 
\end{itemize}

Consequently, for sufficiently large $r$ we obtain a sequence of epimorphisms
$$ E^{p,q}_{r} \rightarrow E^{p,q}_{r+1} \rightarrow E^{p,q}_{r+2} \rightarrow \ldots $$
Let $E^{p,q}_{\infty}$ denote the colimit of this sequence (in the abelian category $\calA$). 

\begin{itemize}
\item[$(ii)$] Let $n \in \Z$, and let $A_n = \pi_{n} \varinjlim(X)$. Then there exists a filtration
$$ \ldots \subseteq F^{-1} A_n \subseteq F^{0} A_n \subseteq F^1 A_n \subseteq \ldots $$
of $A_n$, with $F^{p} A_n \simeq 0$ for $p \ll 0$, and $\varinjlim( F^p A_n) \simeq A_n$.

\item[$(iii)$] For every $p,q \in \Z$, there exists an isomorphism
$E^{p,q}_{\infty} \simeq F^{p} A_{p+q} / F^{p-1} A_{p+q}$ in the abelian category $\calA$. 
\end{itemize}

To prove $(i)$, $(ii)$, and $(iii)$, we first extend $X$ to an object
$\overline{X} \in \Gap( \Z \cup \{ - \infty \}, \calC)$, so that for each $n \in \Z$ we have
$X(n) = \overline{X}(- \infty, n)$. Without loss of generality, we may suppose that
$X(n) \simeq \ast$ for $n < 0$. This implies that $\overline{X}(i,j) \simeq \ast$ for
$i, j < 0$. It follows that $E^{p-r,q+r-1}_{r}$, as a quotient
$\pi_{p+q} \overline{X}(p-2r,p-r)$, is zero for $r > p$. This proves $(i)$.

To satisfy $(ii)$, we set $F^{p} A_n = \im( \pi_{n} X(p) \rightarrow \pi_{n} \varinjlim(X) )$. 
It is clear that $F^{p} A_n \simeq \ast$ for $p < 0$, and the isomorphsim
$\varinjlim F^p A_n \simeq A_n$ follows from the compatibility of the homological functor
$\pi_{n}$ with sequential colimits (Remark \ref{seqcolim}). 

To prove $(iii)$, we note that for $r > p$, the object $E^{p,q}_{r}$ can be identified with
the image of the map $\pi_{p+q} X(p) \simeq \pi_{p+q} \overline{X}(p-r, p)
\rightarrow \pi_{p+q} \overline{X}(p-1, p+r)$. Let $Y = \varinjlim_{r} \overline{X}(p-1, p+r)$.
It follows that $E^{p,q}_{\infty}$ can be identified with the image of the map
$\pi_{p+q} X(p) \stackrel{f}{\rightarrow} \pi_{p+q} Y$. We have a distinguished triangle
$$ X(p-1) \rightarrow \varinjlim(X) \rightarrow Y \rightarrow X(p-1)[1],$$
which induces an exact sequence
$$ 0 \rightarrow F^{p-1} A_{p+q} \rightarrow A_{p+q} \stackrel{f'}{\rightarrow}
\pi_{p+q} Y.$$
We have a commutative triangle
$$ \xymatrix{ & A_{p+q} \ar[dr]^{f'} & \\
\pi_{p+q} X(p) \ar[ur]^{g} \ar[rr]^{f} & & \pi_{p+q} Y. }$$
Since the image of $g$ is $F^{p} A_{p+q}$, we obtain canonical isomorphisms
$$E^{p,q}_{\infty} \simeq \im( f ) \simeq \im( f' | F^p A_{p+q} ) \simeq
F^{p} A_{p+q} / \ker(f') \simeq F^{p} A_{p+q} / F^{p-1} A_{p+q}.$$
This completes the proof.
\end{proof}

\section{The $\infty$-Categorical Dold-Kan Correspondence}\label{doldkan}\index{Dold-Kan correspondence}

Let $\calA$ be an abelian category. Then the classical Dold-Kan correspondence (see \cite{weibel}) asserts that the category $\Fun( \cDelta^{op}, \calA)$ of simplicial objects of $\calA$ is equivalent to the category
$\Chain_{\geq 0}(\calA)$ of (homologically) nonnegatively graded chain complexes
$$ \ldots \stackrel{d}{\rightarrow} A_1 \stackrel{d}{\rightarrow} A_0 \rightarrow 0.$$
In this section, we will prove an analogue of this result when the abelian category $\calA$ is replaced by a stable $\infty$-category.\index{chain complex}

We begin by observing that if $X_{\bigdot}$ is a simplicial object in a stable $\infty$-category $\calC$, then $X_{\bigdot}$ determines a simplicial object of the homotopy category
$\h{\calC}$. The category $\h{\calC}$ is not abelian, but it is additive and has the following additional property (which follows easily from the fact that $\h{\calC}$ admits a triangulated structure):
\begin{itemize}
\item[$(\ast)$] If $i: X \rightarrow Y$ is a morphism in $\h{\calC}$ which admits a left inverse, then
there is an isomorphism $Y \simeq X \oplus X'$ such that $i$ is identified with the map $(\id, 0)$.
\end{itemize}
These conditions are sufficient to construct a Dold-Kan correspondence in $\h{\calC}$. Consequently, every simplicial object $X_{\bigdot}$ of $\calC$ determines a chain complex
$$ \ldots \rightarrow C_1 \rightarrow C_0 \rightarrow 0$$
in the homotopy category $\h{\calC}$. In \S \ref{filttt}, we described another construction which gives rise to the same type of data. More precisely, Lemma \ref{loij} and Remark \ref{makeplex}
show that every $\Z_{\geq 0}$-filtered object 
$$ Y(0) \stackrel{f_1}{\rightarrow} Y(1) \stackrel{f_2}{\rightarrow} \ldots $$
determines a chain complex $C_{\bigdot}$ with values in $\h{\calC}$, where
$C_{n} = \coker(f_i)[-n]$. This suggests a relationship between filtered objects of $\calC$ and simplicial objects of $\calC$. Our goal in this section is to describe this relationship in detail. Our main result, Theorem \ref{dkan}, asserts that the $\infty$-category of simplicial objects of $\calC$ is equivalent to a suitable $\infty$-category of (increasingly) filtered objects of $\calC$. The proof will require several preliminaries.

\begin{lemma}\label{surite}
Let $\calC$ be a stable $\infty$-category. A square
$$ \xymatrix{ X' \ar[r] \ar[d]^{f'} & X \ar[d]^{f} \\
Y' \ar[r] & Y }$$
in $\calC$ is a pullback if and only if the induced map
$\alpha: \coker(f') \rightarrow \coker(f)$ is an equivalence.
\end{lemma}

\begin{proof}
Form an expanded diagram
$$ \xymatrix{ X' \ar[r] \ar[d]^{f'} & X \ar[r] \ar[d]^{f} & 0 \ar[d] \\
Y' \ar[r] & Y \ar[r] & \coker(f) }$$
where the right square is a pushout. Since $\calC$ is stable, the right square is also a pullback.
Lemma \toposref{transplantt} implies that the left square is a pullback if and only if the outer square is a pullback, which is in turn equivalent to the assertion that $\alpha$ is an equivalence.
\end{proof}

\begin{lemma}\label{tuup}
Let $\calC$ be a stable $\infty$-category, let $K$ be a simplicial set, and suppose that $\calC$ admits $K$-indexed colimits. Let $\overline{\alpha}: K^{\triangleright} \times \Delta^1 \rightarrow \calC$
be a natural transformation between a pair of diagrams $\overline{p}, \overline{q}: K^{\triangleright} \rightarrow \calC$. Then $\overline{\alpha}$ is a colimit diagram if and only if
$\coker( \overline{\alpha} ): K^{\triangleright} \rightarrow \calC$ is a colimit diagram.
\end{lemma}

\begin{proof}
Let $p = \overline{p} | K$, $q = \overline{q} | K$, and $\alpha = \overline{\alpha} | K \times \Delta^1)$. Since $\calC$ admits $K$-indexed colimits, there exist colimit diagrams $\overline{p}', \overline{q}': K^{\triangleright} \rightarrow \calC$ extending $p$ and $q$, respectively. We obtain a square
$$ \xymatrix{ \overline{p}' \ar[d] \ar[r] & \overline{p} \ar[d] \\
\overline{q}' \ar[r] & \overline{q} }$$ 
in the $\infty$-category $\Fun( K^{\triangleright}, \calC)$. Let $\infty$ denote the cone point
of $K^{\triangleright}$. Using Corollary \toposref{util}, we deduce that $\alpha$ is a colimit diagram if and only if the induced square
$$ \xymatrix{ \overline{p}'(\infty) \ar[r] \ar[d]^{f'} & \overline{p}(\infty) \ar[d]^{f} \\
\overline{q}'(\infty) \ar[r] & \overline{q}(\infty) }$$
is a pushout. According to Lemma \ref{surite}, this is equivalent to the assertion that
the induced map $\beta: \coker(f') \rightarrow \coker(f)$ is an equivalence. We conclude by observing that $\beta$ can be identified with the natural map
$$ \varinjlim( \coker(\alpha) ) \rightarrow \coker( \overline{\alpha})(\infty),$$
which is an equivalence if and only if $\coker(\overline{\alpha})$ is a colimit diagram.
\end{proof}

Our next result is an analogue of Proposition \stableref{surose} which applies to cubical diagrams of higher dimension.

\begin{proposition}\label{cubate}
Let $\calC$ be a stable $\infty$-category, and let
$\sigma: (\Delta^1)^n \rightarrow \calC$ be a diagram. Then $\sigma$ is a colimit diagram
if and only if $\sigma$ is a limit diagram.
\end{proposition}

\begin{proof}
By symmetry, it will suffice to show that if $\sigma$ is a colimit diagram, then $\sigma$ is also a limit diagram. We work by induction on $n$. If $n=0$, then we must show that every initial object of $\calC$ is also final, which follows from the assumption that $\calC$ has a zero object.
If $n > 0$, then we may identify $\sigma$ with a natural transformation $\alpha: \sigma' \rightarrow \sigma''$ in the $\infty$-category $\Fun( ( \Delta^1)^{n-1}, \calC)$. Assume that $\sigma$ is a colimit diagram. Using Lemma \ref{tuup}, we deduce that $\coker(\alpha)$ is a colimit diagram. Since
$\coker(\alpha) \simeq \ker(\alpha)[1]$, we conclude that $\ker(\alpha)$ is a colimit diagram. Applying the inductive hypothesis, we deduce that $\ker(\alpha)$ is a limit diagram. The dual of Lemma \ref{tuup} now implies that $\sigma$ is a limit diagram, as desired.
\end{proof}

\begin{lemma}\label{subcont}
Fix $n \geq 0$, and let $S$ be a subset of the open interval $(0,1)$ of cardinality $\leq n$.
Let $Y$ be the set of all sequences of real numbers $0 \leq y_1 \leq \ldots \leq y_n \leq 1$ such
that $S \subseteq \{ y_1, \ldots, y_n \}$. Then $Y$ is a contractible topological space.
\end{lemma}

\begin{proof}
Let $S$ have cardinality $m \leq n$, and let $Z$ denote the set of sequences of real numbers
$0 \leq z_1 \leq \ldots \leq z_{n-m} \leq 1$. Then $Z$ is homeomorphic to a topological
$(n-m)$-simplex. Moreover, there is a homeomorphism $f: Z \rightarrow Y$, which carries
a sequence $\{ z_i \}$ to a suitable reordering of the sequence $\{ z_i \} \cup S$.
\end{proof}

\begin{lemma}\label{serpus}
Let $n \leq 0$, let $\cDelta_{\leq n}$
$\cDelta_{\leq n}$ denote the full subcategory of $\cDelta$ spanned
by the objects $\{ [m] \}_{ 0 \leq m \leq n}$, and let $\calI$ denote the full subcategory
of $( \cDelta_{\leq n})_{/[n]}$ spanned by the {\em injective} maps $[m] \rightarrow [n]$.
Then the induced map
$$ \Nerve( \calI)^{op} \rightarrow \Nerve( \cDelta_{\leq n})^{op}$$
is cofinal.
\end{lemma}

\begin{proof}
Fix $m \leq n$, and let $\calJ$ denote the category of diagrams
$$ [m] \leftarrow [k] \stackrel{i}{\rightarrow} [n]$$
where $i$ is injective. According to Theorem \toposref{hollowtt}, it will suffice to show that
the simplicial set $\Nerve(\calJ)$ is weakly contractible (for every $m \leq n$).

Let $X$ denote the simplicial subset of $\Delta^m \times \Delta^n$ spanned by those nondegenerate simplices whose projection to $\Delta^n$ is also nondegenerate. Then
$\Nerve(\calJ)$ can be identified with the barycentric subdivision of $X$. Consequently, it will suffice to show that the topological space $|X|$ is contractible. For this, we will show that the fibers of the map $\phi: |X| \rightarrow | \Delta^m |$ are contractible.

We will identify the topological $m$-simplex $| \Delta^m |$ with the set of all sequences
of real numbers $0 \leq x_1 \leq \ldots \leq x_m \leq 1$. Similarly, we may identify
points of $| \Delta^n |$ with sequences $0 \leq y_1 \leq \ldots \leq y_n \leq 1$. A pair of such sequences determines a point of $X$ if and only if each $x_i$ belongs to the set
$\{ 0, y_1, \ldots, y_n , 1\}$. Consequently, the fiber of $\phi$ over the point
$( 0 \leq x_1 \leq \ldots \leq x_m \leq 1)$ can be identified with the set
$$ Y = \{ 0 \leq y_1 \leq \ldots \leq y_n \leq 1 : \{ x_1, \ldots, x_m \} \subseteq \{ 0, y_1, \ldots, y_n , 1 \}  \} \subseteq | \Delta^n |,$$
which is contractible (Lemma \ref{subcont}).
\end{proof}

\begin{corollary}\label{shimmy}
Let $\calC$ be a stable $\infty$-category, and let $F: \Nerve( \cDelta_{\leq n})^{op} \rightarrow \calC$ be a functor such that $F([m]) \simeq 0$ for all $m < n$. Then there is a canonical isomorphism $\varinjlim(F) \simeq X[n]$ in the homotopy category $h{\calC}$, where $X = F([n])$.
\end{corollary}

\begin{proof}
Let $\calI$ be as in Lemma \ref{serpus}, let
$G''$ denote the composition
$\Nerve( \calI)^{op} \rightarrow \Nerve( \cDelta_{\leq n})^{op} \stackrel{F}{\rightarrow} \calC$, and let $G$ denote the constant map $\Nerve(\calI)^{op} \rightarrow \calC$ taking the value $X$. 
Let $\calI_0$ denote the full subcategory of $\calI$ obtained by deleting the initial object.
There is a canonical map $\alpha: G \rightarrow G''$, and $G' = \ker(\alpha)$ is
a left Kan extension of $G' | \Nerve(\calI_0)^{op}$. We obtain a distinguished triangle
$$ \varinjlim(G') \rightarrow \varinjlim(G) \rightarrow \varinjlim(G'') \rightarrow \varinjlim(G')[1]$$
in the homotopy category $\h{\calC}$. Lemma \ref{serpus} yields an equivalence
$\varinjlim(F) \simeq \varinjlim(G'')$, and Lemma \toposref{kan0} implies the existence of an equivalence $\varinjlim(G') \simeq \varinjlim(G'| \Nerve(\calI_0)^{op})$. 

We now observe that the simplicial set $\Nerve(\calI)^{op}$ can be identified with the barycentric subdivision of the standard $n$-simplex $\Delta^n$, and that $\Nerve(\calI_0)^{op}$ can be identified with the barycentric subdivision of its boundary $\bd \Delta^n$. It follows (see \S \toposref{quasilimit7}) that we may identify the map $\varinjlim(G') \rightarrow \varinjlim(G)$ with the map $\beta: X \otimes ( \bd \Delta^n) \rightarrow X \otimes \Delta^n$. The cokernel of $\beta$ is
canonically isomorphic (in $\h{\calC}$) to the $n$-fold suspension $X[n]$ of $X$.
\end{proof}

\begin{lemma}\label{spokenn}
Let $\calC$ be a stable $\infty$-category, let $n \geq 0$, and let
$F: \Nerve( \cDelta_{+, \leq n})^{op} \rightarrow \calC$ be a functor $($here
$\cDelta_{+, \leq n}$ denotes the full subcategory of $\cDelta_{+}$ spanned by
the objects $\{ [k] \}_{ -1 \leq k \leq n}$ $)$. The following conditions are
equivalent:
\begin{itemize}
\item[$(i)$] The functor $F$ is a left Kan extension of $F | \Nerve( \cDelta_{\leq n})^{op}$.
\item[$(ii)$] The functor $F$ is a right Kan extension of $F | \Nerve( \cDelta_{+, \leq n-1})^{op}$.
\end{itemize}
\end{lemma}

\begin{proof}
Condition $(ii)$ is equivalent to the assertion that the composition
$$ F': \Nerve( \cDelta_{+, \leq n-1}^{op})_{[n]/}^{\triangleleft} \rightarrow
\Nerve( \cDelta_{+, \leq n}^{op}) \stackrel{F}{\rightarrow} \calC$$
is a limit diagram. Since the source of $F$ is isomorphic to $(\Delta^1)^{n+1}$, Proposition
\ref{cubate} asserts that $F'$ is a limit diagram if and only if $F'$ is a colimit diagram.
In view of Lemma \ref{serpus}, $F'$ is a colimit diagram if and only if $F$ is a colimit diagram, which is equivalent to $(i)$.
\end{proof}

\begin{theorem}[$\infty$-Categorical Dold-Kan Correspondence]\label{dkan}
Let $\calC$ be a stable $\infty$-category. Then the $\infty$-categories
$\Fun( \Nerve( \Z_{\geq 0}), \calC)$ and $\Fun( \Nerve( \cDelta)^{op}, \calC)$ are
$($canonically$)$ equivalent to one another.
\end{theorem}\index{Dold-Kan correspondence!$\infty$-categorical version}

\begin{proof}
Our first step is to describe the desired equivalence in more precise terms.
Let $\calI_{+}$ denote the full subcategory of $\Nerve( \Z_{\geq 0} ) \times \Nerve( \cDelta_{+})^{op}$ spanned by those pairs $(n, [m])$, where $m \leq n$, and let
$\calI$ be the full subcategory of $\calI_{+}$ spanned by those pairs
$(n, [m])$ where $0 \leq m \leq n$. We observe that there is a natural projection $p: \calI \rightarrow \Nerve( \cDelta)^{op}$, and a natural embedding $i: \Nerve( \Z_{\geq 0} ) \rightarrow \calI_{+}$, which carries $n \geq 0$ to the object $(n, [-1])$. 

Let $\Fun^{0}( \calI, \calC)$ denote the full subcategory of $\Fun(\calI, \calC)$ spanned
by those functors $F: \calI \rightarrow \calC$ such that, for every $s \leq m \leq n$, 
the image under $F$ of the natural map $(m, [s]) \rightarrow (n, [s])$ is an equivalence in $\calC$.
Let $\Fun^{0}( \calI_{+}, \calC)$ denote the full subcategory of $\Fun( \calI_{+}, \calC)$
spanned by functors $F_{+}: \calI_{+} \rightarrow \calC$ such that $F = F_{+} | \calI$
belongs to $\Fun^{0}(\calI, \calC)$, and $F_{+}$ is a left Kan extension of $F$.
Composition with $p$, composition with $i$, and restriction from $\calI_{+}$ to $\calI$
yields a diagram of $\infty$-categories
$$ \Fun( \Nerve(\cDelta)^{op},\calC ) \stackrel{G}{\rightarrow}
\Fun^{0}( \calI, \calC) \stackrel{G'}{\leftarrow} \Fun^{0}(\calI_{+}, \calC) \stackrel{G''}{\rightarrow}
\Fun( \Nerve( \Z_{\geq 0}), \calC).$$
We will prove that $G$, $G'$, and $G''$ are equivalences of $\infty$-categories.

To show that $G$ is an equivalence of $\infty$-categories, we let
$\calI^{\leq k}$ denote the full subcategory of $\calI$ spanned by pairs
$(n, [m])$ where $m \leq n \leq k$, and let $\calI^{k}$ denote the full subcategory
of $\calI$ spanned by those pairs $(n,[m])$ where $m \leq n = k$.
Then the projection $p$ restricts to an equivalence $\calI^{k} \rightarrow \Nerve( \cDelta_{\leq n})^{op}$.
Let $\Fun^{0}( \calI^{\leq k}, \calC)$ denote the full subcategory of
$\Fun( \calI^{\leq k}, \calC)$ spanned by those functors $F: \calI^{\leq k} \rightarrow \calC$
such that, for every $s \leq m \leq n \leq k$, the image under $F$ of the natural map
$(m, [s]) \rightarrow (n,[s])$ is an equivalence in $\calC$. We observe that
this is equivalent to the condition that $F$ be a right Kan extension of
$F | \calI^{k}$. Using Proposition \toposref{lklk}, we deduce that the restriction map
$r: \Fun^{0}( \calI^{\leq k}, \calC) \rightarrow \Fun( \calI^{k}, \calC)$ is an equivalence
of $\infty$-categories. Composition with $p$ induces a functor
$G_{k}: \Fun( \Nerve(\cDelta_{\leq k})^{op} , \calC) \rightarrow \Fun^{0}( \calI^{\leq k}, \calC)$
which is a section of $r$. It follows that $G_{k}$ is an equivalence of $\infty$-categories.
We can identify $G$ with the homotopy inverse limit of the functors $\varprojlim( G_{k})$, so that $G$ is also an equivalence of $\infty$-categories.

The fact that $G'$ is an equivalence of $\infty$-categories follows immediately from Proposition \toposref{lklk}, since for each $n \geq 0$ the simplicial set $\calI_{/ (n, [-1])}$ is finite and $\calC$ admits finite colimits.

We now show that $G''$ is an equivalence of $\infty$-categories. Let
$\calI^{\leq k}_{+}$ denote the full subcategory of $\calI_{+}$ spanned by pairs
$(n, [m])$ where either $m \leq n \leq k$ or $m=-1$. We let $\calD(k)$ denote the full subcategory of
$\Fun( \calI^{\leq k}_{+}, \calC)$ spanned by those functors $F: \calI^{\leq k}_{+} \rightarrow \calC$
with the following pair of properties:
\begin{itemize}
\item[$(i)$] For every $s \leq m \leq n \leq k$, the image under $F$ of the natural map
$(m, [s]) \rightarrow (n,[s])$ is an equivalence in $\calC$.
\item[$(ii)$] For every $n \leq k$, $F$ is a left Kan extension of $F | \calI^{\leq k}$ at
$(n, [-1])$. 
\end{itemize}

Then $\Fun^{0}( \calI_{+}, \calC)$ is the inverse limit of the tower of restriction maps
$$ \ldots \rightarrow \calD(1) \rightarrow \calD(0) \rightarrow \calD(-1) = \Fun( \Nerve(\Z_{\geq 0}), \calC). $$
To complete the proof, we will show that for each $k \geq 0$, the restriction map
$\calD(k) \rightarrow \calD(k-1)$ is a trivial Kan fibration.

Let $\calI^{\leq k}_{0}$ be the full subcategory of $\calI^{\leq k}_{+}$ obtained by removing the object $(k,[k])$, and let $\calD'(k)$ be the full subcategory of $\Fun( \calI^{\leq k}_0, \calC)$
spanned by those functors $F$ which satisfy condition $(i)$ and satisfy $(ii)$ for $n < k$.
We have restriction maps
$$ \calD(k) \stackrel{\theta}{\rightarrow} \calD'(k) \stackrel{\theta'}{\rightarrow} \calD(k-1).$$
We observe
that a functor $F: \calI^{\leq k}_{0}$ belongs to $\calD'(k)$ if and only if
$F | \calI^{\leq k-1}_{+}$ belongs to $\calD(k-1)$ and $F$ is a left Kan extension of
$F | \calI^{\leq k-1}_{+}$. Using Proposition \toposref{lklk}, we conclude that $\theta'$ is a trivial Kan fibration. 

We will prove that $\theta$ is a trivial Kan fibration by a similar argument. According to Proposition \toposref{lklk}, it will suffice to show that a functor $F: \calI^{\leq k}_{+} \rightarrow \calC$
belongs to $\calD(k)$ if and only if $F | \calI^{\leq k}_0$ belongs to $\calD'(k)$ and
$F$ is a right Kan extension of $F| \calI^{\leq k}_0$. This follows immediately from
Lemma \ref{spokenn} and the observation that the inclusion $\calI^{k} \subseteq \calI^{\leq k}$ is cofinal.
\end{proof}

\begin{remark}
Let $\calC$ be a stable $\infty$-category.
We may informally describe the equivalence of Theorem \ref{dkan} as follows. To
a simplicial object $C_{\bigdot}$ of $\calC$, we assign the filtered object
$$ D(0) \rightarrow D(1) \rightarrow D(2) \rightarrow \ldots$$
where $D(k)$ is the colimit of the $k$-skeleton of $C_{\bigdot}$. In particular, we observe
that colimits $\varinjlim D(j)$ can be identified with geometric realizations of the simplicial object $C_{\bigdot}$.
\end{remark}

\begin{remark}\label{saltine}
Let $\calC$ be a stable $\infty$-category, and let $X$ be a simplicial object
of $\calC$. Using the Dold-Kan correspondence, we can associate to $X$ a chain complex
$$ \ldots \rightarrow C_2 \rightarrow C_1 \rightarrow C_0 \rightarrow 0$$
in the triangulated category $\h{\calC}$. More precisely, for each $n \geq 0$, let
$L_n \in \calC$ denote the $n$th {\it latching object} of $X$ (see \S \toposref{coreed}), so that $X$ determines a canonical map $\alpha: L_n \rightarrow X_n$. Then
$C_{n} \simeq \coker(\alpha)$, where the cokernel can be formed either in the $\infty$-category
$\calC$ or in its homotopy category $\h{\calC}$ (since $L_n$ is actually a direct summand of $X_n$). 

Using Theorem \ref{dkan}, we can also associate to $X$ a filtered object
$$ D(0) \rightarrow D(1) \rightarrow D(2) \rightarrow \ldots$$
of $\calC$. Using Lemma \ref{loij} and Remark \ref{makeplex}, we can associate to
this filtered abject another chain complex
$$ \ldots \rightarrow C'_1 \rightarrow C'_0 \rightarrow 0$$
with values in $\h{\calC}$. 
For each $n \geq 0$, let $X(n)$ denote the restriction of $X$ to $\Nerve( \cDelta_{\leq n}^{op})$, and let $X'(n)$ be a left Kan extension of $X(n-1)$ to $\Nerve( \cDelta_{\leq n}^{op})$. 
Then we have a canonical map $\beta = X'(n) \rightarrow X(n)$, which induces an equivalence
$X'(n)_{m} \rightarrow X(n)_m$ for $m < n$, while $X'(n)_n$ can be identified with the
latching object $L_n$. Let $X''(n) = \coker(\beta)$. Then $X''(n)_m = 0$ for $m < n$, while
$X''(n)_n \simeq C_n$. Corollary \ref{shimmy} determines a canonical isomorphism
$\varinjlim X''(n) \simeq C_{n}[n]$ in the homotopy category $\h{\calC}$
The map $D(n-1) \rightarrow D(n)$ can be identified with the composition
$$ D(n-1) \simeq \varinjlim X(n-1) \simeq \varinjlim X'(n) \rightarrow \varinjlim X(n) \simeq D(n).$$
It follows it follows that $C'_{n} \simeq \coker( D(n-1) \rightarrow D(n) )[-n] \simeq X''(n)_n[-n]$
is canonically isomorphic to $C_{n}$. It is not difficult to show that these isomorphisms are
compatible with the differentials, so that we obtain an isomorphism of chain complexes
$C_{\bigdot} \simeq C'_{\bigdot}$ with values in the triangulated category $\h{\calC}$.
\end{remark}

\begin{remark}
Let $\calC$ be a stable $\infty$-category, let $X_{\bigdot}$ be a simplicial object
of $\calC$, let $$D(0) \rightarrow D(1) \rightarrow \ldots $$
be the associated filtered object. Using the classical Dold-Kan correspondence and Remark
\ref{saltine}, we conclude that each $X_{n}$ is equivalent to a finite coproduct of objects of the form $\coker( D(m-1) \rightarrow D(m) )[-m]$, where $0 \leq m \leq n$ (here $D(-1) \simeq 0$ by convention). 
\end{remark}

\begin{remark}
Let $\calC$ be a stable $\infty$-category equipped with a t-structure, whose heart is equivalent to (the nerve of) an abelian category $\calA$. Let $X_{\bigdot}$ be a simplicial object of $\calC$, and let
$$ D(0) \rightarrow D(1) \rightarrow D(2) \rightarrow \ldots $$
be the associated filtered object (Theorem \ref{dkan}). Using Definition \ref{sumtuo} (and Lemma \ref{loij}), we can associate to this filtered object a spectral sequence
$\{ E_{r}^{p,q}, d_r \}_{r \geq 1}$ in the abelian category $\calA$. In view of Remarks \ref{saltine2} and \ref{saltine}, for each $q \in \Z$ we can identify the complex $(E_1^{\bigdot, q}, d_1)$ with
the normalized chain complex associated to the simplicial object $\pi_{q} X_{\bigdot}$ of $\calA$. 
Under the hypotheses of Proposition \ref{conseq}, this spectral sequence converges
to a filtration on the homotopy groups $\pi_{p+q} \varinjlim( D(n) ) \simeq \pi_{p+q} | X_{\bigdot} |$.

It possible to consider a slight variation on the spectral sequence described above. Namely, one can construct a new spectral sequence $\{ \overline{E}_{r}^{p,q}, d_r \}_{r \geq 1}$ which is isomorphic to $\{ E_{r}^{p,q}, d_r \}_{r \geq 1}$ from the $E_2$-page onward, but with
$\overline{E}_{1}^{\bigdot, q}$ given by the {\em unnormalized} chain complex
of $\pi_{q} X_{\bigdot}$. We can then write simply $\overline{E}^{p,q}_{1} \simeq \pi_{q} X_p$.
\end{remark}

\section{Homological Algebra}\label{stable10}

Let $\calA$ be an abelian category. In classical homological algebra, it is customary to associate to $\calA$ a certain triangulated category, called the {\it derived category} of $\calA$, the objects of which are chain complexes with values in $\calA$. In this section, we will review the theory of derived categories from the perspective of higher category theory. To simplify the discussion, we primarily consider only abelian categories $\calA$ which have enough projective objects (the dual case of abelian categories with enough injective objects can be understood by passing to the opposite category).\index{derived category!of an abelian category}

We begin by considering an arbitrary additive category $\calA$.\index{ZZZChainA@$\Chain(\calA)$}\index{ZZZChaingz@$\Chain_{\geq 0}(\calA)$}\index{chain complex}
Let $\Chain(\calA)$ denote the category whose objects are chain complexes
$$ \ldots \rightarrow A_1 \rightarrow A_0 \rightarrow A_{-1} \rightarrow \ldots$$
with values in $\calA$. The category $\Chain(\calA)$ is naturally enriched over simplicial sets. 
For $A_{\bigdot}, B_{\bigdot} \in \Chain(\calA)$, the simplicial set
$\bHom_{\Chain(\calA)}( A_{\bigdot}, B_{\bigdot} )$ is characterized by the property
that for every finite simplicial set $K$ there is a natural bijection
$$ \Hom_{\sSet}( K, \bHom_{\Chain(\calA)}( A_{\bigdot}, B_{\bigdot} ))
\simeq \Hom_{\Chain(\calA)}(A_{\bigdot} \otimes C_{\bigdot}(K), B_{\bigdot} )
.$$
Here $C_{\bigdot}(K)$ denotes the normalized chain complex for computing the homology of $K$, so that $C_{n}(K)$ is a free abelian group whose generators are in bijection with the nondegenerate $n$-simplices of $K$. Unwinding the definitions, we see that the vertices
of $\bHom_{\Chain(\calA)}(A_{\bigdot}, B_{\bigdot})$ are just the maps of chain complexes
from $A_{\bigdot}$ to $B_{\bigdot}$. An edge $e$ of $\bHom_{\Chain(\calA)}(A_{\bigdot}, B_{\bigdot})$ is determined by three pieces of data:
\begin{itemize}
\item[$(i)$] A vertex $d_0(e)$, corresponding to a chain map $f: A_{\bigdot} \rightarrow B_{\bigdot}$.
\item[$(ii)$] A vertex $d_1(e)$, corresponding to a chain map $g: A_{\bigdot} \rightarrow B_{\bigdot}$.
\item[$(iii)$] A map $h: A_{\bigdot} \rightarrow B_{\bigdot+1}$, which determines a chain homotopy from $f$ to $g$.
\end{itemize}

\begin{remark}\index{Dold-Kan correspondence}
Let $\Ab$ be the category of abelian groups, and let
$\Chain_{\geq 0}(\Ab)$ denote the full subcategory of $\Chain(\Ab)$ spanned by
those complexes $A_{\bigdot}$ such that $A_{n} \simeq 0$ for all $n < 0$. The classical Dold-Kan correspondence (see \cite{weibel}) asserts that $\Chain_{\geq 0}(\Ab)$ is equivalent to the category of {\em simplicial} abelian groups. In particular, there is a forgetful functor $\theta: \Chain_{\geq 0}(\Ab) \rightarrow \sSet$. 

Given a pair of complexes $A_{\bigdot}, B_{\bigdot} \in \Chain(\calA)$, the mapping space
$\bHom_{ \Chain(\calA)}(A_{\bigdot}, B_{\bigdot})$ can be defined as follows:
\begin{itemize}
\item[$(1)$] First, we extract the mapping complex $$[ A_{\bigdot}, B_{\bigdot} ] \in \Chain(\Ab),$$
where $[ A_{\bigdot}, B_{\bigdot}]_{n} = \prod \Hom_{\calA}( A_{m}, B_{n+m})$. 
\item[$(2)$] The inclusion $\Chain_{\geq 0}(\Ab) \subseteq \Chain(\Ab)$ has a right adjoint, which associates to an arbitrary chain complex $M_{\bigdot}$ the truncated complex
$$ \ldots \rightarrow M_1 \rightarrow \ker(M_0 \rightarrow M_{-1}) \rightarrow 0 \rightarrow \ldots$$
Applying this functor to $[A_{\bigdot}, B_{\bigdot}]$, we obtain a new complex $[A_{\bigdot}, B_{\bigdot}]_{\geq 0}$, whose degree zero term coincides with the set of chain maps from $A_{\bigdot}$ to $B_{\bigdot}$.
\item[$(3)$] Applying the Dold-Kan correspondence $\theta$, we can convert 
the chain complex $[A_{\bigdot}, B_{\bigdot}]_{\geq 0}$ into a simplicial set $\bHom_{ \Chain(\calA)}(A_{\bigdot}, B_{\bigdot})$. 
\end{itemize}
Because every simplicial abelian group is a Kan complex, the simplicial category $\Chain(\calA)$
is automatically {\em fibrant}. 
\end{remark}

\begin{remark}
Let $\calA$ be an additive category, and let $A_{\bigdot}, B_{\bigdot} \in \Chain(\calA)$. The homotopy group $$\pi_n \bHom_{\Chain(\calA)}(A_{\bigdot}, B_{\bigdot})$$ can be identified with the
group of chain-homotopy classes of maps from $A_{\bigdot}$ to
$B_{\bigdot+n}$.
\end{remark}

\begin{example}\label{urbus}
Let $\calA$ be an abelian category, and let $A_{\bigdot}, B_{\bigdot} \in \Chain(\calA)$. Suppose
that $A_n \simeq 0$ for $n < 0$, and that $B_{n} \simeq 0$ for $n > 0$. Then the simplicial
set $\bHom_{\Chain(\calA)}(A_{\bigdot}, B_{\bigdot})$ is constant, with value
$\Hom_{\calA}( \HH_0(A_{\bigdot}), \HH_0(B_{\bigdot}))$. 
\end{example}

\begin{lemma}\label{slurpus}
Let $\calA$ be an additive category. Then:
\begin{itemize}
\item[$(1)$] Let $$\xymatrix{ A_{\bigdot} \ar[r]^{f} \ar[d] & B_{\bigdot} \ar[d] \\
A'_{\bigdot} \ar[r] & B'_{\bigdot} }$$
be a pushout diagram in the (ordinary) category $\Chain(\calA)$, and suppose that
$f$ is {\em degreewise split} (so that each $B_{n} \simeq A_{n} \oplus C_{n}$, for some
$C_{n} \in \calA$). Then the above diagram determines a homotopy pushout square in the $\infty$-category
$\Nerve( \Chain(\calA) )$. 
\item[$(2)$] The $\infty$-category $\Nerve( \Chain(\calA) )$ is stable.
\end{itemize}
\end{lemma}

\begin{proof}
To prove $(1)$, it will suffice (Theorem \toposref{colimcomparee}) to show that
for every $D_{\bigdot} \in \Chain(\calA)$, the associated diagram of simplicial sets
$$ \xymatrix{ \bHom_{\Chain(\calA)}( B'_{\bigdot}, D_{\bigdot} ) \ar[r] \ar[d] & \bHom_{ \Chain(\calA)}( A'_{\bigdot}, D_{\bigdot} ) \ar[d] \\
\bHom_{\Chain(\calA)}( B_{\bigdot}, D_{\bigdot} ) \ar[r]^{f'} & \bHom_{ \Chain(\calA)}( A_{\bigdot}, D_{\bigdot}) }$$
is homotopy Cartesian. The above diagram is obviously a pullback, it will suffice to prove that $f'$ is a Kan fibration. This follows from the fact that $f'$ is the map of simplicial sets associated (under the Dold-Kan correspondence) to a map between complexes of abelian groups which is surjective in positive (homological) degrees.

It follows from $(1)$ that the $\infty$-category $\Nerve( \Chain(\calA) )$ admits pushouts:
it suffices to observe that any morphism $f: A_{\bigdot} \rightarrow B_{\bigdot}$ is chain homotopy-equivalent to a morphism which is degreewise split (replace $B_{\bigdot}$ by the mapping cylinder of $f$). It is obvious that $\Nerve( \Chain( \calA) )$ has a zero object (since $\Chain(\calA)$ has a zero object). Moreover, we can use $(1)$ to describe the suspension functor on $\Chain(\calA)$: for
each $A_{\bigdot} \in \Chain(\calA)$, let $C( A_{\bigdot} )$ denote the cone of $A_{\bigdot}$, so that 
$C(A_{\bigdot}) \simeq 0$ and there is a pushout diagram
$$ \xymatrix{ A_{\bigdot} \ar[r] \ar[d] & C( A_{\bigdot} ) \ar[d] \\
0 \ar[r] & A_{\bigdot-1}. }$$
It follows that the suspension functor $\Sigma$ can be identified with the shift functor
$$ A_{\bigdot} \mapsto A_{\bigdot-1}.$$
In particular, we conclude that $\Sigma$ is an equivalence of $\infty$-categories, so that $\Chain(\calA)$ is stable (Proposition \ref{charstut}). 
\end{proof}

\begin{remark}\label{stabte}
Let $\calA$ be an additive category, and let $\Chain'(\calA)$ be a full subcategory of
$\Chain(\calA)$. Suppose that $\Chain'(\calA)$ is stable under translations and the formation of mapping cones. Then the proof of Lemma \ref{slurpus} shows that $\Nerve( \Chain'(\calA) )$ is a stable subcategory of $\Nerve( \Chain(\calA) )$. In particular, if $\Chain^{\dplus}(\calA)$ denotes the full subcategory of $\Chain(\calA)$ spanned by those complexes $A_{\bigdot}$ such that
$A_{n} \simeq 0$ for $n \ll 0$, then $\Nerve( \Chain^{\dplus}(\calA) )$ is a stable subcategory of
$\Nerve( \Chain(\calA))$.\index{ZZZChainplus@$\Chain^{\dplus}(\calA)$}
\end{remark}

\begin{definition}\index{projective object}\index{ZZZDPlusA@$\calD^{\dplus}(\calA)$}\index{derived $\infty$-category!of an abelian category}\label{smucky}
Let $\calA$ be an abelian category with enough projective objects. We let
$\calD^{\dplus}(\calA)$ denote the nerve of the simplicial category $\Chain^{\dplus}(\calA_0)$, where
$\calA_0 \subseteq \calA$ is the full subcategory spanned by the projective objects of $\calA$. 
We will refer to $\calD^{\dplus}(\calA)$ as the {\it derived $\infty$-category of $\calA$}.
\end{definition}

\begin{remark}
The homotopy category $\h{ \calD^{\dplus}(\calA)}$ can be described as follows: objects are given by (bounded above) chain complexes of projective objects of $\calA$, and morphisms are given by homotopy classes of chain maps. Consequently, $\h{ \calD^{\dplus}(\calA)}$ can be identified
with the derived category of $\calA$ studied in classical homological algebra (with appropriate boundedness conditions imposed).
\end{remark}

\begin{lemma}\label{starit}
Let $\calA$ be an abelian category, and let $P_{\bigdot} \in \Chain(\calA)$ be a complex
of projective objects of $\calA$ such that $P_{n} \simeq 0$ for $n \ll 0$. Let 
$Q_{\bigdot} \rightarrow Q'_{\bigdot}$ be a quasi-isomorphism in $\Chain(\calA)$. Then
the induced map
$$ \bHom_{ \Chain(\calA)}( P_{\bigdot}, Q_{\bigdot}) \rightarrow \bHom_{ \Chain(\calA)}( P_{\bigdot}, Q'_{\bigdot})$$ is a homotopy equivalence.
\end{lemma}

\begin{proof}
We observe that $P_{\bigdot}$ is a homotopy colimit of its naive truncations
$$ \ldots \rightarrow 0 \rightarrow P_{n} \rightarrow P_{n-1} \rightarrow \ldots. $$
It therefore suffices to prove the result for each of these truncations, so we may assume that $P_{\bigdot}$ is concentrated in finitely many degrees. Working by induction, we can reduce to the case where $P_{\bigdot}$ is concentrated in a single degree. Shifting, we can reduce to the case where $P_{\bigdot}$ consists of a single projective object $P$ concentrated in degree zero. Since
$P$ is projective, we have isomorphisms
$$\Ext^i_{\Nerve(\Chain(\calA))}( P_{\bigdot}, Q_{\bigdot} )
\simeq \Hom_{\calA}( P, \HH_{-i}( Q_{\bigdot}) )
\simeq \Hom_{\calA}(P, \HH_{-i}( Q'_{\bigdot}) ) \simeq \Ext^i_{ \Nerve(\Chain(\calA))}( P_{\bigdot}, Q'_{\bigdot}).$$
\end{proof}

\begin{lemma}\label{surg}
Let $\calA$ be an abelian category. Suppose that $P_{\bigdot}, Q_{\bigdot} \in \Chain(\calA)$ have the following properties:
\begin{itemize}
\item[$(1)$] Each $P_{n}$ is projective, and $P_{n} \simeq 0$ for $n < 0$.
\item[$(2)$] The homologies $\HH_n(Q_{\bigdot})$ vanish for $n > 0$.
\end{itemize}
Then the space $\bHom_{\Chain(\calA)}(P_{\bigdot}, Q_{\bigdot})$ is discrete, and we have a canonical isomorphism of abelian groups $$\Ext^0( P_{\bigdot}, Q_{\bigdot} ) \simeq \Hom_{\calA}( \HH_0(P_{\bigdot}),
\HH_0(Q_{\bigdot}) ).$$
\end{lemma}

\begin{proof}
Let $Q'_{\bigdot}$ be the complex
$$ \ldots \rightarrow 0 \rightarrow \coker( Q_1 \rightarrow Q_0 ) \rightarrow Q_{-1} \rightarrow \ldots.$$
Condition $(2)$ implies that the canonical map $Q_{\bigdot} \rightarrow Q'_{\bigdot}$ is a quasi-isomorphism. In view of $(1)$ and Lemma \ref{starit}, it will suffice to prove the
result after replacing $Q_{\bigdot}$ by $Q'_{\bigdot}$. The result now follows from Example \ref{urbus}.
\end{proof}

\begin{proposition}\label{hurkin}\index{t-structure!on $\calD^{\dplus}(\calA)$}
Let $\calA$ be an abelian category with enough projective objects. Then:
\begin{itemize}
\item[$(1)$] The $\infty$-category $\calD^{\dplus}(\calA)$ is stable.
\item[$(2)$] Let $\calD^{\dplus}_{\geq 0}(\calA)$ be the full subcategory of
$\calD^{\dplus}(\calA)$ spanned by those complexes $A_{\bigdot}$ such that the
homology objects $\HH_n(A_{\bigdot}) \in \calA$ vanish for $n < 0$, and let
$\calD^{\dplus}_{\leq 0}(\calA)$ be defined similarly. Then $( \calD^{\dplus}_{\geq 0}(\calA),
\calD^{\dplus}_{\leq 0}(\calA))$ determines a t-structure on $\calD^{\dplus}(\calA)$.
\item[$(3)$] The heart of $\calD^{\dplus}(\calA)$ is equivalent to (the nerve of) the abelian category $\calA$.
\end{itemize}
\end{proposition}
\begin{proof}
Assertion $(1)$ follows from Remark \ref{stabte}. 

To prove $(2)$, we first make the following observation:
\begin{itemize}
\item[$(\ast)$] For any object $A_{\bigdot} \in \Chain(\calA)$, there exists a map
$f: P_{\bigdot} \rightarrow A_{\bigdot}$ where each $P_n$ is projective, $P_n \simeq 0$ for $n < 0$, and the induced map $\HH_k(P_{\bigdot}) \rightarrow \HH_{k}(A_{\bigdot})$ is an isomorphism for $k \geq 0$.
\end{itemize}
This is proven by a standard argument in homological algebra, using the assumption that
$\calA$ has enough projectives. We also note that if $A_{\bigdot} \in \calD^{\dplus}(\calA)$ and the homologies $\HH_n(A_{\bigdot})$ vanish for $n < 0$, then $f$ is a quasi-isomorphism between projective complexes and therefore a chain homotopy equivalence.

It is obvious that $\calD^{\dplus}_{\leq 0}(\calA)[-1] \subseteq \calD^{\dplus}_{\leq 0}(\calA)$ and
$\calD^{\dplus}_{\geq 0}(\calA)[1] \subseteq \calD^{\dplus}_{\geq 0}(\calA)$. Suppose now that
$A_{\bigdot} \in \calD^{\dplus}_{\geq 0}(\calA)$ and $B_{\bigdot} \in \calD^{\dplus}_{\leq -1}(\calA)$; we wish to show that $\Ext^0_{\calD^{\dplus}(\calA)}(A_{\bigdot}, B_{\bigdot} ) \simeq 0$. Using $(\ast)$, we may reduce to the case where $A_{n} \simeq 0$ for $n < 0$. The desired result now follows immediately from Lemma \ref{surg}. Finally, choose an arbitrary object $A_{\bigdot} \in \calD^{\dplus}(\calA)$, and let
$f: P_{\bigdot} \rightarrow A_{\bigdot}$ be as in $(\ast)$. It is easy to see that
$\coker(f) \in \calD^{\dplus}_{\leq -1}(\calA)$. This completes the proof of $(2)$.

To prove $(3)$, we begin by observing that the functor
$A_{\bigdot} \mapsto \HH_0( A_{\bigdot} )$ determines a functor
$\theta: \Nerve( \Chain(\calA) ) \rightarrow \Nerve(\calA)$. Let $\calC \subseteq \Nerve(\Chain(\calA))$
be the full subcategory spanned by complexes $P_{\bigdot}$ such that each $P_{n}$ is projective, $P_{n} \simeq 0$ for $n < 0$, and $\HH_{n}( P_{\bigdot} ) \simeq 0$ for $n \neq 0$. Assertion
$(\ast)$ implies that the inclusion $\calC \subseteq \heart{ \calD^{\dplus}(\calA) }$ is an equivalence of $\infty$-categories. Lemma \ref{surg} implies that $\theta | \calC$ is fully faithful. Finally, we can apply $(\ast)$ in the case where $A_{\bigdot}$ is concentrated in degree zero to deduce that
$\theta | \calC$ is essentially surjective. This proves $(3)$.
\end{proof}

\begin{remark}
Let $\calA$ be an abelian category with enough projective objects. Then
$\calD^{\dplus}(\calA)$ is a colocalization of $\Nerve( \Chain^{\dplus}(\calA) )$. To prove this, 
it will suffice to show that for every $A_{\bigdot} \in \Chain^{\dplus}(\calA)$, there exists a map of chain complexes $f: P_{\bigdot} \rightarrow A_{\bigdot}$ where $P_{\bigdot} \in \calD^{\dplus}(\calA)$, and
such that $f$ induces a homotopy equivalence 
$$\bHom_{\Chain(\calA)}( Q_{\bigdot}, P_{\bigdot}) \rightarrow \bHom_{ \Chain(\calA)}( Q_{\bigdot}, A_{\bigdot} )$$ for every $Q_{\bigdot} \in \calD^{\dplus}(\calA)$ (Proposition \toposref{testreflect}). 
According Lemma \ref{starit}, it will suffice to choose $f$ to be a quasi-isomorphism; the existence now follows from $(\ast)$ in the proof of Proposition \ref{hurkin}.

Let $L: \Nerve( \Chain^{\dplus}(\calA) ) \rightarrow \calD^{\dplus}(\calA)$ be a right adjoint to the inclusion.
Roughly speaking, the functor $L$ associates to each complex $A_{\bigdot}$ a projective
resolution $P_{\bigdot}$ as above. We observe that, if $f: A_{\bigdot} \rightarrow B_{\bigdot}$ is a map of complexes, then $Lf$ is a chain homotopy equivalence if and only if $f$ is a quasi-isomorphism. Consequently, we may regard $\calD^{\dplus}(\calA)$ as the $\infty$-category
obtained from $\Nerve( \Chain^{\dplus}(\calA) )$ by inverting all quasi-isomorphism. 
\end{remark}
 
\section{The Universal Property of $\calD^{\dplus}(\calA)$}\label{stable14}
 
In this section, we will apply the results of \S \toposref{stable11} and \S \toposref{stable12} to characterize the derived $\infty$-category $\calD^{\dplus}(\calA)$ by a universal mapping property. Here $\calA$ denotes an abelian category with enough projective objects; to simplify the discussion, we will assume that $\calA$ is small.\index{derived $\infty$-category!universal property of}

Let $\calA_0 \subseteq \calA$ be the full subcategory of $\calA$ spanned by the projective objects, and let $\bfA$ denote the category of product-preserving functors from $\calA_{0}^{op}$ to the category of simplicial sets, as in \S \toposref{stable12}. Let $\calA^{\vee}$ denote the category of product-preserving functors from $\calA_{0}^{op}$ to sets, so that we can identify $\bfA$ with the category of simplicial objects of $\calA^{\vee}$. Our first goal is to understand the category $\calA^{\vee}$.\index{ZZZcalAvee@$\calA^{\vee}$}
 
\begin{lemma}\label{tooce}
Let $\calA$ be an abelian category with enough projective objects, and let $\calB$ be an arbitrary category which admits finite colimits. Let $\calC$ be the category of right exact functors from
$\calA$ to $\calB$, and let $\calC'$ be the category of coproduct-preserving functors
from $\calA_0$ to $\calB$. Then the restriction functor
$\theta: \calC \rightarrow \calC'$ is an equivalence of categories.
\end{lemma}

\begin{proof}
We will describe an explicit construction of an inverse functor. Let $f: \calA_0 \rightarrow \calB$
be a functor which preserves finite coproducts. Let $A \in \calA$ be an arbitrary object. Since $\calA$ has enough projectives, there exists a projective resolution
$$\ldots \rightarrow P_1 \stackrel{u}{\rightarrow} P_0 \rightarrow A \rightarrow 0.$$
We now define $F(A)$ to be the coequalizer of the map
$$\xymatrix{ f(P_1) \ar@<.4ex>[r]^{f(0)} \ar@<-.4ex>[r]_{f(u)} & f(P_0)}.$$
Of course, this definition appears to depend not only on $A$ but on a choice of projective resolution. However, because any two projective resolutions of $A$ are chain homotopy equivalent to one another, $F(A)$ is well-defined up to canonical isomorphism.
It is easy to see that $F: \calA \rightarrow \calB$ is a right exact functor which extends $f$, and that $F$ is uniquely determined (up to unique isomorphism) by these properties.
\end{proof}
 
\begin{proposition}
Let $\calA$ be an abelian category with enough projective objects. Then:
\begin{itemize}
\item[$(1)$] The category $\calA^{\vee}$ can be identified with the category of $\Ind$-objects of $\calA$.
\item[$(2)$] The category $\calA^{\vee}$ is abelian.
\item[$(3)$] The abelian category $\calA^{\vee}$ has enough projective objects.
\end{itemize}
\end{proposition}

\begin{proof}
Assertion $(1)$ follows immediately from Lemma \ref{tooce} (taking $\calB$ to be the opposite of the category of sets). Part $(2)$ follows formally from $(1)$ and the assumption that $\calA$ is an abelian category (see, for example, \cite{artinmazur}). We may identify $\calA$ with a full subcategory of $\calA^{\vee}$ via the Yoneda embedding. Moreover, if $P$ is a projective object of $\calA$, then $P$ is also projective when viewed as an object of $\calA^{\vee}$. An arbitrary object
of $\calA^{\vee}$ can be written as a filtered colimit $A = \varinjlim \{ A_{\alpha} \}$, where each $A_{\alpha} \in \calA$.  Using the assumption that $\calA$ has enough projective objects, we can choose
epimorphisms $P_{\alpha} \rightarrow A_{\alpha}$, where each $P_{\alpha}$ is projective. We then have an epimorphism $\oplus P_{\alpha} \rightarrow A$. Since $\oplus P_{\alpha}$ is projective, we conclude that $\calA^{\vee}$ has enough projectives.
\end{proof}

\begin{warning}
Let $\calA$ be an abelian category with enough projective objects, and let
$\bfA$ be the category of product-preserving functors $\calA^{op}_{0} \rightarrow \sSet$. The Dold-Kan correspondence determines an equivalence of categories $\theta: \bfA \simeq \Chain_{\geq 0}( \calA^{\vee} )$. However, this is {\em not} an equivalence of simplicial categories. Let $K$ be a simplicial set, and let $\Z K$ denote the free simplicial abelian group generated by $K$ (so that
the group of $n$-simplices of $\Z K$ is the free abelian group generated by the set of $n$-simplices of $K$, for each $n \geq 0$). Then $\bfA$ is tensored over the category of simplicial sets in two different ways:
\begin{itemize}
\item[$(i)$] Given a simplicial set $K$ and an object $A_{\bigdot} \in \bfA$ viewed as a simplicial object of $\calA^{\vee}$, we can form the tensor product $A_{\bigdot} \otimes K$ given by 
the formula $(A_{\bigdot} \otimes K)_{n} = A_{n} \otimes (\Z K)_n$.
\item[$(ii)$] Given a simplicial set $K$ and an object $A_{\bigdot} \in \bfA$, we can construct a new object $A_{\bigdot} \boxtimes K$, which is characterized by the existence of an isomorphism
$$ \theta( A_{\bigdot} \boxtimes K) \simeq \theta(A_{\bigdot}) \otimes \theta'( \Z K)$$
in the category $\Chain( \calA^{\vee} )$. Here $\theta'( \Z K)$ denotes the object of $\Chain( \Ab)$ determined by $\Z K$.
\end{itemize}
However, it is easy to see that both of these simplicial structures on $\bfA$ are compatible with the model structure of Proposition \toposref{sutcoat}. Moreover, the classical {\it Alexander-Whitney map} determines a natural transformation $A_{\bigdot} \otimes K \rightarrow A_{\bigdot} \boxtimes K$, which endows $\theta^{-1}: \Chain_{\geq 0}( \calA^{\vee}) \rightarrow \bfA$ with the structure of a simplicial functor.
\end{warning}

We observe that every object of $\bfA$ is fibrant, and that an object of $\bfA$ is cofibrant if and only if it corresponds (under $\theta$) to a complex of projective objects of $\calA^{\vee}$. Applying 
Corollary \toposref{urchug}, we obtain an equivalence of $\infty$-categories
$\calD^{\dplus}_{\geq 0}( \calA^{\vee} ) \rightarrow \Nerve( \bfA^{\degree} )$. Here
$\calD^{\dplus}_{\geq 0}( \calA^{\vee} )$ denotes the full subcategory of $\calD^{\dplus}(\calA^{\vee})$ spanned by those complexes $P_{\bigdot}$ such that $P_n \simeq 0$ for $n < 0$. Composing with the equivalence of Corollary \toposref{smokerr}, we obtain the following result:

\begin{proposition}\label{urchite}
Let $\calA$ be an abelian category with enough projective objects. Then there exists
an equivalence of $\infty$-categories
$$ \psi: \calD^{\dplus}_{\geq 0}( \calA^{\vee} ) \rightarrow \calP_{\Sigma}( \Nerve(\calA_0) )$$
whose composition with the inclusion $\Nerve( \calA_0 ) \subseteq \calD^{\dplus}_{\geq }( \calA^{\vee})$
is equivalent to the Yoneda embedding $\Nerve(\calA_0) \rightarrow \calP_{\Sigma}( \Nerve(\calA_0))$. 
\end{proposition}

\begin{remark}\label{sandboat}
We can identify $\calD^{\dplus}(\calA)$ with a full subcategory of $\calD^{\dplus}(\calA^{\vee})$. Moreover, an object $P_{\bigdot} \in \calD^{\dplus}( \calA^{\vee} )$ belongs to the essential image of $\calD^{\dplus}(\calA)$ if and only if the homologies $\HH_n(P_{\bigdot})$ belong to $\calA$, for all $n \in \Z$. 
\end{remark}

 \begin{proposition}
Let $\calA$ be an abelian category with enough projective objects. Then
the t-structure on $\calD^{\dplus}(\calA)$ is right bounded and left complete.
\end{proposition}

\begin{proof}
The right boundedness of $\calD^{\dplus}(\calA)$ is obvious. To prove the left completeness,
we must show that $\calD^{\dplus}(\calA)$ is a homotopy inverse limit of the tower of $\infty$-categories
$$ \ldots \rightarrow \calD^{\dplus}(\calA)_{\leq 1} \rightarrow \calD^{\dplus}(\calA)_{\leq 0} \rightarrow \ldots$$
Invoking the right boundedness of $\calD^{\dplus}(\calA)$, we may reduce to proving that for each $n \in \Z$, $\calD^{\dplus}(\calA)_{\geq n}$ is a homotopy inverse limit of the tower
$$ \ldots \rightarrow \calD^{\dplus}( \calA)_{ \leq 1, \geq n} \rightarrow \calD^{\dplus}(\calA)_{\leq 0, \geq n} \rightarrow \ldots $$
Shifting if necessary, we may suppose that $n = 0$. Using Remark \ref{sandboat}, we can replace $\calA$ by $\calA^{\vee}$. For each $k \geq 0$, we let $\calP_{\Sigma}^{\leq k}( \Nerve( \calA_0) )$
denote\index{ZZZcalPSigmak@$\calP_{\Sigma}^{\leq k}(\calC)$} the $\infty$-category of product-preserving functors from $\Nerve( \calA_0 )^{op}$ to 
$\tau_{\leq k} \SSet$; equivalently, we can define $\calP_{\Sigma}^{\leq k}( \Nerve( \calA_0 ))$ to be the $\infty$-category of $k$-truncated objects of $\calP_{\Sigma}( \Nerve(\calA_0) )$. We observe that the equivalence $\psi$ of Proposition \ref{urchite} restricts to an equivalence
$$ \psi(k): \calD^{\dplus}_{\geq 0}(\calA^{\vee})_{\leq k} \rightarrow \calP_{\Sigma}^{\leq k}( \Nerve(\calA_0) ).$$
Consequently, it will suffice to show that $\calP_{\Sigma}( \Nerve(\calA_0 ))$ is a homotopy inverse limit for the tower
$$ \ldots \rightarrow \calP_{\Sigma}^{\leq 1}( \Nerve(\calA_0) ) \rightarrow
\calP_{\Sigma}^{\leq 0}( \Nerve(\calA_0 )).$$
Since the truncation functors on $\SSet$ commute with finite products (Lemma \toposref{slurpy} ), we may reduce to the problem of showing that $\SSet$ is a homotopy inverse limit of the tower
$$ \ldots \rightarrow \tau_{\leq 1} \SSet \rightarrow \tau_{\leq 0} \SSet. $$
This amounts to the classical fact that every space $X$ can be recovered as the limit of its Postnikov tower (see for example \S \toposref{homdim}). 

\end{proof}
 
Our goal is to characterize the derived $\infty$-category $\calD^{\dplus}(\calA)$ by a universal mapping property. Propositions \ref{urchite} and \toposref{surottt} give a characterization of $\calD^{\dplus}_{\geq 0}(\calA^{\vee})$ of the right flavor. The next step is to understand the embedding of $\calD^{\dplus}_{\geq 0}(\calA)$ into $\calD^{\dplus}_{\geq 0}( \calA^{\vee} )$.

\begin{definition}\label{hurtman}\index{functor!right t-exact}\index{right t-exact}\index{functor!left t-exact} Let $\calC$ and $\calC'$ be stable $\infty$-categories equipped with t-structures. We will say
that a functor $f: \calC \rightarrow \calC'$ is {\it right t-exact} if it is exact, and carries
$\calC_{\geq 0}$ into $\calC'_{\geq 0}$.
\end{definition}

\begin{lemma}\label{simpenough}
\begin{itemize}
\item[$(1)$] Let $\calC$ be an $\infty$-category which admits finite coproducts and geometric realizations. Then $\calC$ admits all finite colimits. Conversely, if
$\calC$ is an $n$-category which admits finite colimits, then $\calC$ admits geometric realizations.
\item[$(2)$] Let $F: \calC \rightarrow \calD$ be a functor between $\infty$-categories which admit finite coproducts and geometric realizations. If $F$ preserves finite coproducts and geometric realizations, then $F$ is right exact. The converse holds if $\calC$ and $\calD$ are $n$-categories.
\end{itemize}
\end{lemma}

\begin{proof}
We will prove $(1)$; the proof of $(2)$ follows by the same argument. 
Now suppose that $\calC$ admits finite coproducts and geometric realizations of simplicial objects. We wish to show that $\calC$ admits all finite colimits. According to Proposition \toposref{appendixdiagram}, it will suffice to show that $\calC$ admits coequalizers. Let $\cDelta_{s}^{ \leq 1}$ be the full subcategory of $\cDelta$ spanned by the objects $[0]$ and $[1]$, and {\em injective} maps between them, so that a coequalizer diagram in $\calC$ can be identified with a functor $\Nerve( \cDelta_{s}^{ \leq 1})^{op} \rightarrow \calC$. Let
$j: \Nerve( \cDelta_{s}^{\leq 1})^{op} \rightarrow \Nerve( \cDelta)^{op}$ be the inclusion functor. We claim that every diagram $f: \Nerve( \cDelta_{s}^{ \leq 1})^{op} \rightarrow \calC$ has a left Kan extension along $j$. To prove this, it suffices to show that for each $n \geq 0$, the associated diagram
$$ \Nerve( \cDelta_{s}^{ \leq 1})^{op} \times_{ \Nerve(\cDelta)^{op} } (\Nerve(\cDelta)^{op})_{[n]/}
\rightarrow \calC$$
has a colimit. We now observe that this last colimit is equivalent to a coproduct: more precisely, we have $(j_{!} f)([n]) \simeq f([0]) \coprod f([1]) \coprod \ldots \coprod f([1])$, where there are precisely
$n$ summands equivalent to $f([1])$. Since $\calC$ admits finite coproducts, the desired Kan extension $j_{!} f$ exists. We now observe that $\varinjlim(f)$ can be identified with
$\varinjlim( j_{!} f)$, and the latter exists in virtue of our assumption that $\calC$ admits geometric realizations for simplicial objects.

Now suppose that $\calC$ is an $n$-category which admits finite colimits; we wish to show that $\calC$ admits geometric realizations. Passing to a larger universe if necessary, we may suppose that $\calC$ is small. Let $\calD = \Ind(\calC)$, and let $\calC' \subseteq \calD$ denote the essential image of the Yoneda embedding $j: \calC \rightarrow \calD$. Then $\calD$ admits small colimits (Theorem \toposref{pretop}) and $j$ is fully faithful (Proposition \toposref{fulfaith}); it will therefore suffice to show that $\calC'$ is stable under geometric realization of simplicial objects in $\calD$. 

Fix a simplicial object $U_{\bigdot}: \Nerve(\cDelta)^{op} \rightarrow \calC' \subseteq \calD$.
Let $V_{\bigdot}: \Nerve(\cDelta)^{op} \rightarrow \calD$ be a left Kan extension of
$U_{\bigdot} | \Nerve(\cDelta^{\leq n})^{op}$, and $\alpha_{\bigdot}: V_{\bigdot} \rightarrow U_{\bigdot}$ the induced map. The geometric realization of $V_{\bigdot}$ can be identified with the colimit of $U_{\bigdot} | \Nerve( \cDelta^{\leq n})^{op}$, and therefore belongs to $\calC'$ since
$\calC'$ is stable under finite colimits in $\calD$ (Proposition \toposref{turnke}). Consequently, it will suffice to prove that $\alpha_{\bigdot}$ induces an equivalence from the geometric realization of $V_{\bigdot}$ to the geometric realization of $U_{\bigdot}$. 

Let $L: \calP(\calC) \rightarrow \calD$ be a left adjoint to the inclusion. 
Let $|U_{\bigdot}|$ and $|V_{\bigdot}|$ be colimits of $U_{\bigdot}$ and $V_{\bigdot}$ in the
$\infty$-category $\calP(\calC)$, and let $| \alpha_{\bigdot}|: |V_{\bigdot} | \rightarrow |U_{\bigdot} |$ be the induced map. We wish to show that $L | \alpha_{\bigdot} |$ is an equivalence in $\calD$. 
Since $\calC$ is an $n$-category, we have inclusions
$\Ind(\calC) \subseteq \Fun( \calC^{op}, \tau_{\leq n-1} \SSet) \subseteq \calP(\calC)$. 
It follows that $L$ factors through the truncation functor $\tau_{\leq n-1}: \calP(\calC)
\rightarrow \calP(\calC)$. Consequently, it will suffice to prove that $\tau_{\leq n-1} |\alpha_{\bigdot}|$ is an equivalence in $\calP(\calC)$. For this, it will suffice to show
that the morphism $| \alpha_{\bigdot} |$ is $n$-connective (in the sense of Definition \toposref{stooog}). This follows from Lemma \toposref{bball2}, since $\alpha_k: V_k \rightarrow U_{k}$ is an equivalence for $k \leq n$.
\end{proof}

\begin{lemma}\label{surput}
Let $\calC$ and $\calC'$ be stable $\infty$-categories equipped with t-structures. Let $\theta: \Fun(\calC, \calC')
\rightarrow \Fun( \calC_{\geq 0}, \calC')$ be the restriction map. Then:
\begin{itemize}
\item[$(1)$] If $\calC$ is right-bounded, then $\theta$ induces an equivalence from
the full subcategory of $\Fun(\calC, \calC')$ spanned by the right t-exact functors to the full subcategory of $\Fun(\calC_{\geq 0}, \calC'_{\geq 0})$ spanned by the right exact functors.
\item[$(2)$] Let $\calC$ and $\calC'$ be left complete. Then the $\infty$-categories
$\calC_{\geq 0}$ and $\calC'_{\geq 0}$ admit geometric realizations of simplicial objects.
Furthermore, a functor $F: \calC_{\geq 0} \rightarrow \calC'_{\geq 0}$ is right exact if and only if
if it preserves finite coproducts and geometric realizations of simplicial objects.
\end{itemize}
\end{lemma}

\begin{proof}
We first prove $(1)$. If $\calC$ is right bounded, then $\Fun(\calC, \calC')$ is the (homotopy) inverse limit of the tower
$$ \ldots \rightarrow \Fun( \calC_{\geq -1}, \calC' ) \rightarrow \Fun( \calC_{\geq 0}, \calC' ),$$
where the functors are given by restriction. The full subcategory of right t-exact functors is then given by the homotopy inverse limit
$$ \ldots \rightarrow \Fun'( \calC_{\geq -1}, \calC'_{\geq -1}) \stackrel{\theta(0)}{\rightarrow} \Fun( \calC_{\geq 0}, \calC'_{\geq 0} )$$
where $\Fun'(\calC, \calD)$ denotes the full subcategory of $\Fun(\calC, \calD)$ spanned by the right exact functors. To complete the proof, it will suffice to show that this tower is essentially constant; in other words, that each $\theta(n)$ is an equivalence of $\infty$-categories. Without loss of generality, we may suppose that $n=0$. Choose shift functors on the $\infty$-categories
$\calC$ and $\calC'$, and define
$$ \psi: \Fun'( \calC_{\geq 0}, \calC'_{\geq 0}) \rightarrow \Fun'( \calC_{ \geq -1}, \calC'_{\geq -1})$$
by the formula $\psi(F)= \Sigma^{-1} \circ F \circ \Sigma$. We claim that $\psi$ is a homotopy inverse to $\theta(0)$. To prove this, we observe that the right exactness of $F \in \Fun'(\calC_{\geq 0}, \calC'_{\geq 0})$, $G \in \Fun'( \calC_{\geq -1}, \calC'_{\geq -1})$ determines canonical equivalences
$$ (\theta(0) \circ \psi)(F) \simeq F$$
$$ ( \psi \circ \theta(0))(G) \simeq G.$$

We now prove $(2)$. Since $\calC$ is left complete, $\calC_{\geq 0}$ is the (homotopy) inverse limit of the tower of $\infty$-categories $\{ (\calC_{\geq 0})_{\leq n} \}$ with transition maps given by right exact truncation functors. Lemma \ref{simpenough} implies that each $(\calC_{\geq 0})_{\leq n}$
admits geometric realizations of simplicial objects, and that each of the truncation functors preserves geometric realizations of simplicial objects. It follows that $\calC_{\geq 0}$ admits geometric realizations for simplicial objects. Similarly, $\calC'_{\geq 0}$ admits geometric realizations for simplicial objects.

If $F$ preserves finite coproducts and geometric realizations of simplicial objects, then $F$ is right exact (Lemma \ref{simpenough}). Conversely, suppose that $F$ is right exact; we wish to prove that $F$ preserves geometric realizations of simplicial objects. It will suffice to show that each composition
$$ \calC_{\geq 0} \stackrel{F}{\rightarrow} \calC'_{\geq 0} \stackrel{\tau_{\leq n}}{\rightarrow}
( \calC'_{\geq 0})_{\leq n}$$
preserves geometric realizations of simplicial objects. We observe that, in virtue of the right exactness of $F$, this functor is equivalent to the composition
$$ \calC_{\leq 0} \stackrel{\tau_{\leq n}}{\rightarrow} (\calC_{\geq 0})_{\leq n}
\stackrel{ \tau_{\leq n} \circ F }{\rightarrow} ( \calC'_{\geq 0})_{\leq n}.$$
It will therefore suffice to prove that $\tau_{\leq n} \circ F$ preserves geometric realizations
of simplicial objects, which follows from Lemma \ref{simpenough} since both the source and target are equivalent to $n$-categories.
\end{proof}

\begin{lemma}\label{geomista}
Let $\calA$ be a small abelian category with enough projective objects, and let
$\calC \subseteq \calP_{\Sigma}( \Nerve( \calA_0 ))$ be the essential image of
$\calD^{\dplus}_{\geq 0}(\calA) \subseteq \calD^{\dplus}_{\geq 0}( \calA^{\vee} )$ under the equivalence
$ \psi: \calD^{\dplus}_{\geq 0}( \calA^{\vee} ) \rightarrow \calP_{\Sigma}( \Nerve(\calA_0) )$ of Proposition \ref{urchite}. Then $\calC$ is the smallest full subcategory of $\calP( \Nerve(\calA_0) )$ which is closed under geometric realization and contains the essential image of the Yoneda embedding.
\end{lemma}

\begin{proof}
It is clear that $\calC$ contains the essential image of the Yoneda embedding. Lemma \ref{surput} implies that $\calD^{\dplus}_{\geq 0}( \calA)$ admits geometric realizations and that the inclusion
$\calD^{\dplus}_{\geq 0}(\calA) \subseteq \calD^{\dplus}_{\geq 0}(\calA^{\vee})$ preserves geometric realizations. It follows that $\calC$ is closed under geometric realizations in $\calP( \Nerve(\calA_0))$. 

To complete the proof, we will show that every object of $X \in \calD^{\dplus}_{\geq 0}(\calA)$ can be obtained as the geometric realization, in $\calD^{\dplus}_{\geq 0}(\calA^{\vee})$, of a simplicial
object $P_{\bigdot}$ such that each $P_{n} \in \calD^{\dplus}_{\geq 0}( \calA^{\vee} )$ consists of a projective object of $\calA$, concentrated in degree zero. In fact, we can take $P_{\bigdot}$ to be the simplicial object of $\calA_0$ which corresponds to $X \in \Chain_{\geq 0}( \calA_0 )$ under the Dold-Kan correspondence. It follows from Theorem \toposref{colimcomparee} and Proposition \toposref{eggers} that $X$ can be identified with the geometric realization of $P_{\bigdot}$.
\end{proof}

We are now ready to establish our characterization of $\calD^{\dplus}_{\geq 0}(\calA)$. 

\begin{theorem}\label{uternal}\index{derived $\infty$-category!universal property of}
Let $\calA$ be an abelian category with enough projective objects, $\calA_0 \subseteq \calA$ the full subcategory spanned by the projective objects, and $\calC$
an arbitrary $\infty$-category which admits geometric realizations. Let
$\Fun'( \calD^{\dplus}_{\geq 0}(\calA), \calC)$ denote the full subcategory of
$\Fun(\calD^{\dplus}_{\geq 0}(\calA), \calC)$ spanned by those functors which preserve geometric realizations. Then:
\begin{itemize}
\item[$(1)$] The restriction map
$$ \Fun'( \calD^{\dplus}_{\geq 0}(\calA), \calC) \rightarrow \Fun( \Nerve(\calA_0), \calC)$$ is an equivalence of $\infty$-categories. 
\item[$(2)$] A functor $F \in \Fun'( \calD^{\dplus}_{\geq 0}(\calA), \calC)$ preserves preserves finite coproducts if and only if the restriction $F | \Nerve(\calA_0)$ preserves finite coproducts.
\end{itemize}
\end{theorem}

\begin{proof}
Part $(1)$ follows from Lemma \ref{geomista}, Remark \toposref{poweryoga}, and Proposition \toposref{lklk}. The ``only if'' direction of $(2)$ is obvious. To prove the ``if'' direction, let us suppose that $F| \Nerve(\calA_0)$ preserves finite coproducts. We may assume without loss of generality that $\calC$ admits filtered colimits (Lemma \toposref{diverti}), so that $F$ extends to a functor $F': \calD^{\dplus}_{\geq 0}( \calA^{\vee} )$ which preserves filtered colimits and geometric realizations (Propositions \ref{urchite} and \toposref{surottt}). It follows from Proposition \toposref{surottt} that $F'$ preserves finite coproducts, so that $F = F' | \calD^{\dplus}_{\geq 0}(\calA)$ also preserves finite coproducts.
\end{proof}

\begin{corollary}\label{truceborn}
Let $\calA$ be an abelian category with enough projective objects, and let
$\calC$ be a stable $\infty$-category equipped with a left complete t-structure. 
Then the restriction functor 
$$ \Fun( \calD^{\dplus}(\calA), \calC) \rightarrow \Fun( \Nerve(\calA_0), \calC )$$
induces an equivalence from the full subcategory of $\Fun(\calD^{\dplus}(\calA), \calC)$ spanned by the right t-exact functors to the full subcategory of $\Fun( \Nerve(\calA_0), \calC_{\geq 0} )$ spanned by functors which preserve finite coproducts $($here $\calA_0$ denotes the full subcategory of $\calA$ spanned by the projective objects$)$. 
\end{corollary}

\begin{proof}
Let $\Fun'( \calD^{\dplus}(\calA), \calC)$ be the full subcategory of $\Fun(\calD^{\dplus}(\calA), \calC)$ spanned by the right t-exact functors. Lemma \ref{surput} implies that $\Fun'(\calD^{\dplus}(\calA), \calC)$ is equivalent (via restriction) to the full subcategory $$\Fun'( \calD^{\dplus}_{\geq 0}(\calA), \calC_{\geq 0}) \subseteq \Fun( \calD^{\dplus}_{\geq 0}(\calA), \calC_{\geq 0})$$ spanned by those functors which preserve finite coproducts and geometric realizations of simplicial objects.
Theorem \ref{uternal} and Proposition \toposref{surottt} allow us to identify $\Fun'( \calD^{\dplus}_{\geq 0}(\calA), \calC_{\geq 0})$ with the $\infty$-category of finite-coproduct preserving functors
from $\Nerve(\calA_0)$ into $\calC_{\geq 0}$.
\end{proof}

\begin{corollary}\label{sutty}
Let $\calA$ be an abelian category with enough projective objects, let $\calC$ be a stable $\infty$-category equipped with a left complete t-structure, and let $\calE \subseteq \Fun( \calD^{\dplus}(\calA), \calC)$ be the full subcategory spanned by those right t-exact functors which carry projective objects of $\calA$ into the heart of $\calC$. Then $\calE$ is equivalent to 
(the nerve of) the ordinary category of right exact functors from $\calA$ to the heart of $\calC$. 
\end{corollary}

\begin{proof}
Corollary \ref{truceborn} implies that the restriction map
$$ \calE \rightarrow \Fun( \Nerve(\calA_0), \heart{\calC} )$$
is fully faithful, and that the essential image of $\theta$ consists of the collection of coproduct-preserving functors from $\Nerve(\calA_0)$ to $\heart{\calC}$. Lemma \ref{tooce} allows us to identify the latter $\infty$-category with the nerve of the category of right exact functors from $\calA$ to the heart of $\calC$. 
\end{proof}

If $\calA$ and $\calC$ are as in Proposition \ref{truceborn}, then any
right exact functor from $\calA$ to $\heart{\calC}$ can be extended (in an essentially unique way) to a functor $\calD^{\dplus}(\calA) \rightarrow \calC$. In particular, if the abelian category 
$\heart{\calC}$ has enough projective objects, then we obtain an induced map
$\calD^{\dplus}( \heart{\calC} ) \rightarrow \calC$. 

\begin{example}\index{left derived functor}\index{functor!left derived}
Let $\calA$ and $\calB$ be abelian categories equipped with enough projective objects. Then
any right-exact functor $f: \calA \rightarrow \calB$ extends to a
right t-exact functor $F: \calD^{\dplus}(\calA) \rightarrow \calD^{\dplus}(\calB)$. One typically refers to
$F$ as the {\it left derived functor} of $f$. 
\end{example}

\begin{example}
Let $\Spectra$ be the stable $\infty$-category of spectra (see \S \ref{stable8}),
with its natural t-structure. Then the heart of $\Spectra$
is equivalent to the category $\calA$ of abelian groups. We therefore obtain
a functor $\calD^{\dplus}(\calA) \rightarrow \Spectra$, which carries a complex of abelian groups to the corresponding {\it generalized Eilenberg-MacLane spectrum}.
\end{example}

\section{Presentable Stable $\infty$-Categories}\label{stable15}

In this section, we will study the class of {\em presentable} stable $\infty$-categories: that is, stable $\infty$-categories which admit small colimits and are generated (under colimits) by a set of small objects. In the stable setting, the condition of presentability can be formulated in reasonably simple terms.

\begin{proposition}\label{stablecolimits}
\begin{itemize}
\item[$(1)$] A stable $\infty$-category $\calC$ admits small colimits if and only if $\calC$ admits small coproducts.

\item[$(2)$] Let $F: \calC \rightarrow \calD$ be an exact functor between stable $\infty$-categories satisfying $(1)$. Then $F$ preserves small colimits if and only if $F$ preserves small coproducts.

\item[$(3)$] Let $\calC$ be a stable $\infty$-category satisfying $(1)$, and let $X$ be an object of $\calC$. Then $X$ is compact if and only if the following condition is satisfied:
\begin{itemize}
\item[$(\ast)$] For every map $f: X \rightarrow \coprod_{ \alpha \in A } Y_{\alpha}$ in $\calC$,
there exists a finite subset $A_0 \subseteq A$ such that $f$ factors (up to homotopy)
through $\coprod_{ \alpha \in A_0} Y_{\alpha}$.
\end{itemize}
\end{itemize}
\end{proposition}

\begin{proof}
The ``only if'' direction of $(1)$ is obvious, and the converse follows from Proposition \toposref{appendixdiagram}. Assertion $(2)$ can be proven in the same way.

The ``only if'' direction of $(3)$ follows from the fact that an arbitrary coproduct
$\coprod_{\alpha \in A} Y_{\alpha}$ can be obtained as a filtered colimit of 
finite coproducts $\coprod_{ \alpha \in A_0} Y_{\alpha}$ (see \S \toposref{quasilimit1}). 
Conversely, suppose that an object $X \in \calC$ satisfies $(\ast)$; we wish to show that
$X$ is compact. Let $f: \calC \rightarrow \widehat{\SSet}$ be the functor corepresented by $C$ (recall that $\hat{\SSet}$ denotes the $\infty$-category of spaces which are not necessarily small).
Proposition \toposref{yonedaprop} implies that $f$ is left exact. According to Proposition \ref{urtusk21}, we can assume that $f = \Omega^{\infty} \circ F$, where
$F: \calC \rightarrow \widehat{\Spectra}$ is an exact functor; here $\widehat{\Spectra}$ denotes the $\infty$-category of spectra which are not necessarily small. We wish to prove that $f$ preserves
filtered colimits. Since $\Omega^{\infty}$ preserves filtered colimits, it will suffice to show that
$F$ preserves {\em all} colimits. In view of $(2)$, it will suffice to show that $F$ preserves coproducts. In virtue of Remark \ref{sunrise}, we are reduced to showing that each of the induced functors
$$ \calC \stackrel{F}{\rightarrow} \widehat{\Spectra} \stackrel{\pi_n}{\rightarrow} \Nerve(\Ab)$$
preserves coproducts, where $\Ab$ denotes the category of (not necessarily small) abelian groups. Shifting if necessary, we may suppose $n = 0$. In other words, we must show that for any collection of objects $\{ Y_{\alpha} \}_{\alpha \in A}$, the natural map
$$ \theta: \bigoplus \Ext^0_{\calC}(X,Y_{\alpha}) \rightarrow \Ext^0_{\calC}(X, \coprod Y_{\alpha})$$
is an isomorphism of abelian groups. The surjectivity of $\theta$ amounts to the assumption $(\ast)$, while the injectivity follows from the observations that each
$Y_{\alpha}$ is a retract of the coproduct $\coprod Y_{\alpha}$ and that the natural map
$\bigoplus \Ext^0_{\calC}(X, Y_{\alpha}) \rightarrow \prod \Ext^0_{\calC}(X, Y_{\alpha})$ 
is injective. 
\end{proof}

If $\calC$ is a stable $\infty$-category, then we will say that an object $X \in \calC$
{\it generates} $\calC$ if the condition $\pi_0 \bHom_{\calC}(X,Y) \simeq \ast$
implies that $Y$ is a zero object of $\calC$.

\begin{corollary}\label{stablecolimcor}
Let $\calC$ be a stable $\infty$-category. Then $\calC$ is
presentable if and only if the following conditions are satisfied:

\begin{itemize}
\item[$(1)$] The $\infty$-category $\calC$ admits small coproducts.
\item[$(2)$] The homotopy category $\h{\calC}$ is locally small.
\item[$(3)$] There exists regular cardinal $\kappa$ and a $\kappa$-compact
generator $X \in \calC$.
\end{itemize}
\end{corollary}

\begin{proof}
Suppose first that $\calC$ is presentable. Conditions $(1)$ and $(2)$ are obvious. To establish $(3)$, we may assume without loss of generality that $\calC$ is an accessible localization of
$\calP(\calD)$, for some small $\infty$-category $\calD$. Let $F: \calP(\calD) \rightarrow \calC$
be the localization functor and $G$ its right adjoint. Let $j: \calD \rightarrow \calP(\calD)$ be the Yoneda embedding, and let
$X$ be a coproduct of all suspensions (see \S \ref{stable3}) of objects of the form
$F( j(D) )$, where $D \in \calD$. Since $\calC$ is presentable, $X$ is $\kappa$-compact provided that $\kappa$ is sufficiently large. We claim that $X$ generates $\calC$. To prove this, we consider an arbitrary $Y \in \calC$ such that $\pi_0 \bHom_{\calC}( X,Y) \simeq \ast$. It follows that
the space $$\bHom_{\calC}( F(j(D)), Y) \simeq \bHom_{\calP(\calD)}( j(D), G(Y)) \simeq
G(Y)(D)$$
is contractible for all $D \in \calD$, so that $G(Y)$ is a final object of $\calP(\calD)$. Since
$G$ is fully faithful, we conclude that $Y$ is a final object of $\calC$, as desired.

Conversely, suppose that $(1)$, $(2)$, and $(3)$ are satisfied. We first claim that $\calC$ is itself locally small. It will suffice to show that for every morphism $f: X \rightarrow Y$ in $\calC$ and
every $n \geq 0$, the homotopy group $\pi_n( \Hom^{\rght}_{\calC}(X,Y), f)$ is small. We note that
$\Hom^{\rght}_{\calC}( X,Y)$ is equivalent to the loop space of $\Hom^{\rght}_{\calC}(X,Y[1])$; the question is therefore independent of base point, so we may assume that $f$ is the zero map. We conclude that the relevant homotopy group can identified with $\Hom_{\h{\calC}}( X[n], Y)$, which
is small in virtue of assumption $(2)$.

Fix a regular cardinal $\kappa$ and a $\kappa$-compact object $X$ which generates $\calC$.
We now define a transfinite sequence of full subcategories
$$ \calC(0) \subseteq \calC(1) \subseteq \ldots $$
as follows. Let $\calC(0)$ be the full subcategory of $\calC$ spanned by the objects
$\{ X[n] \}_{ n \in \Z}$. If $\lambda$ is a limit ordinal, let $\calC(\lambda) = \bigcup_{\beta < \lambda} \calC(\beta)$.
Finally, let $\calC(\alpha+1)$ be the full subcategory of $\calC$ spanned by all objects
which can be obtained as the colimit of $\kappa$-small diagrams in $\calC(\alpha)$. Since
$\calC$ is locally small, it follows that each $\calC(\alpha)$ is essentially small. It follows by
induction that each $\calC(\alpha)$ consists of $\kappa$-compact objects of $\calC$ and is
stable under translation. Finally, we observe that $\calC(\kappa)$ is stable under $\kappa$-small colimits. It follows from Lemma \ref{surefoot} that $\calC(\kappa)$ is a stable subcategory
of $\calC$. Choose a small $\infty$-category $\calD$ and an equivalence $f: \calD \rightarrow \calC(\kappa)$. According to Proposition \toposref{uterr}, we may suppose that
$f$ factors as a composition
$$ \calD \stackrel{j}{\rightarrow} \Ind_{\kappa}(\calD) \stackrel{F}{\rightarrow} \calC$$
where $j$ is the Yoneda embedding and $F$ is a $\kappa$-continuous, fully faithful functor.
We will complete the proof by showing that $F$ is an equivalence.

Proposition \toposref{sumatch} implies that $F$ preserves small colimits. It follows that
$F$ admits a right adjoint $G: \calC \rightarrow \Ind_{\kappa}(\calD)$ (Remark \toposref{afi}). We wish to show that the counit map $u: F \circ G \rightarrow \id_{\calC}$ is an equivalence of functors. Choose an object $Z \in \calC$, and
let $Y$ be a cokernel for the induced map $u_Z: (F \circ G)(Z) \rightarrow Z$. 
Since $F$ is fully faithful, 
$G(u_Z)$ is an equivalence. Because $G$ is an exact functor, we deduce that
$G(Y) = 0$. It follows that $\bHom_{\calC}( F(D), Y) \simeq \bHom_{ \Ind_{\kappa}(\calD)}(D, G(Y) ) \simeq \ast$ for all $D \in \Ind_{\kappa}(\calD)$. In particular, we conclude that
$\pi_0 \bHom_{\calC}(X,Y) \simeq \ast$. Since $X$ generates $\calC$, we deduce that
$Y \simeq 0$. Thus $u_Z$ is an equivalence as desired.
\end{proof}

\begin{remark}
In view of Proposition \ref{stablecolimits} and Corollary \ref{stablecolimcor}, the hypothesis
that a stable $\infty$-category $\calC$ be compactly generated can be formulated entirely in terms of the homotopy category $\h{\calC}$. Consequently, one can study this condition entirely in the setting of triangulated categories, without making reference to (or assuming the existence of) an underlying stable $\infty$-category. We refer to reader to \cite{neeman} for further discussion.
\end{remark}

The following result gives a good class of examples of presentable $\infty$-categories.




\begin{proposition}\label{urtusk22}
Let $\calC$ and $\calD$ be presentable $\infty$-categories, and suppose that $\calD$ is stable. 

\begin{itemize}
\item[$(1)$] The $\infty$-category $\Stab(\calC)$ is presentable.

\item[$(2)$] The functor $\Omega^{\infty}: \Stab(\calC) \rightarrow \calC$ admits a left adjoint
$\Sigma^{\infty}: \calC \rightarrow \Stab(\calC)$. 

\item[$(3)$] An exact functor $G: \calD \rightarrow \Stab(\calC)$ admits a left adjoint if and only if
$\Omega^{\infty} \circ G: \calD \rightarrow \calC$ admits a left adjoint.
\end{itemize}
\end{proposition}

\begin{proof}
We first prove $(1)$. Assume that $\calC$ is presentable, and let $1$ be a final object of $\calC$. Then $\calC_{\ast}$ is equivalent to $\calC_{1/}$, and therefore presentable (Proposition \toposref{stabslic}). The loop functor $\Omega: \calC_{\ast} \rightarrow \calC_{\ast}$ admits a left adjoint $\Sigma: \calC_{\ast} \rightarrow \calC_{\ast}$. Consequently, we may view the tower
$$ \ldots \stackrel{\Omega}{\rightarrow} \calC_{\ast} \stackrel{\Omega}{\rightarrow} \calC_{\ast}$$
as a diagram in the $\infty$-category $\RPres$. Invoking Theorem \toposref{surbus}, we deduce $(1)$ and the following modified versions of $(2)$ and $(3)$:

\begin{itemize}
\item[$(2')$] The functor $\Omega^{\infty}_{\ast}: \Stab(\calC) \rightarrow \calC_{\ast}$ admits a left adjoint
$\Sigma^{\infty}_{\ast}: \calC_{\ast} \rightarrow \Stab(\calC)$. 

\item[$(3')$] An exact functor $G: \calD \rightarrow \Stab(\calC)$ admits a left adjoint if and only if
$\Omega^{\infty}_{\ast} \circ G: \calD \rightarrow \calC_{\ast}$ admits a left adjoint.
\end{itemize}

To complete the proof, it will suffice to verify the following:
\begin{itemize}
\item[$(2'')$] The forgetful functor $\calC_{\ast} \rightarrow \calC$ admits a left adjoint
$\calC \rightarrow \calC_{\ast}$. 

\item[$(3'')$] A functor $G: \calD \rightarrow \calC_{\ast}$ admits a left adjoint if and only if
the composition $\calD \stackrel{G}{\rightarrow} \calC_{\ast} \rightarrow \calC$ admits a left adjoint.
\end{itemize}

To prove $(2'')$ and $(3'')$, we recall that a functor $G$ between presentable $\infty$-categories admits a left adjoint if and only if $G$ preserves small limits and small, $\kappa$-filtered colimits, for some regular cardinal $\kappa$ (Corollary \toposref{adjointfunctor}). The desired results now follow from Propositions \toposref{goeselse} and \toposref{needed17}.
\end{proof}

\begin{corollary}\label{mapprop}
Let $\calC$ and $\calD$ be presentable $\infty$-categories, and suppose that $\calD$ is stable.
Then composition with $\Sigma^{\infty}: \calC \rightarrow \Stab(\calC)$ induces an equivalence
$$ \LPres( \Stab(\calC), \calD) \rightarrow \LPres(\calC, \calD).$$
\end{corollary}

\begin{proof}
This is equivalent to the assertion that composition with $\Omega^{\infty}$ induces an equivalence
$$ \RPres( \calD, \Stab(\calC) ) \rightarrow \RPres( \calD, \calC ),$$ which follows from Propositions \ref{urtusk21} and \ref{urtusk22}.
\end{proof}

Using Corollary \ref{mapprop}, we obtain another characterization of the $\infty$-category of
spectra. Let $\Sphere \in \Spectra$ denote the image under $\Sigma^{\infty}: \SSet \rightarrow \Spectra$ of the final object $\ast \in \SSet$. We will refer to $\Sphere$ as the {\it sphere spectrum}. 

\begin{corollary}\label{choccrok}
Let $\calD$ be a stable, presentable $\infty$-category. Then evaluation on the sphere spectrum induces an equivalence of $\infty$-categories
$$ \theta: \LPres( \Spectra, \calD) \rightarrow \calD.$$
\end{corollary}

In other words, we may regard the $\infty$-category $\Spectra$ as the stable $\infty$-category which is freely generated, under colimits, by a single object.

\begin{proof}
We can factor the evaluation map $\theta$ as a composition
$$ \LPres( \Stab(\SSet), \calD ) \stackrel{\theta'}{\rightarrow} \LPres( \SSet, \calD) \stackrel{\theta''}{\rightarrow} \calD$$
where $\theta'$ is given by composition with $\Sigma^{\infty}$ and $\theta''$ by evaluation
at the final object of $\SSet$. We now observe that $\theta'$ and $\theta''$ are both equivalences of $\infty$-categories (Corollary \ref{mapprop} and Theorem \toposref{charpresheaf}).
\end{proof}

We conclude this section by establishing a characterization of the class of stable, presentable $\infty$-categories.

\begin{lemma}\label{lexloc}
Let $\calC$ be a stable $\infty$-category, and let $\calC' \subseteq \calC$ be a localization of $\calC$. Let $L: \calC \rightarrow \calC'$ be a left adjoint to the inclusion. Then $L$ is left exact if and only if $\calC'$ is stable.
\end{lemma}

\begin{proof}
The ``if'' direction follows from Proposition \ref{funrose}, since $L$ is right exact. Conversely, suppose that $L$ is left exact. Since $\calC'$ is a localization of $\calC$, it is closed under finite limits. In particular, it is closed under the formation of kernels and contains a zero object of $\calC$. To complete the proof, it will suffice to show that $\calC'$ is stable under the formation of pushouts in $\calC$. Choose a pushout diagram $\sigma: \Delta^1 \times \Delta^1 \rightarrow \calC$
$$ \xymatrix{ X \ar[r] \ar[d] & X' \ar[d] \\
Y \ar[r] & Y' }$$
in $\calC$, where $X, X', Y \in \calC'$. Proposition \ref{surose} implies that $\sigma$ is also a pullback square. Let $L: \calC \rightarrow \calC'$ be a left adjoint to the inclusion. Since $L$ is left exact, 
we obtain a pullback square $L(\sigma)$:
$$ \xymatrix{ LX \ar[r] \ar[d] & LX' \ar[d] \\
LY \ar[r] & LY'. }$$
Applying Proposition \ref{surose} again, we deduce that $L(\sigma)$ is a pushout square in $\calC$. The natural transformation $\sigma \rightarrow L(\sigma)$ is an equivalence when restricted to $\Lambda^2_0$, and therefore induces an equivalence $Y' \rightarrow LY'$. It follows that $Y'$ belongs to the essential image of $\calC'$, as desired.
\end{proof}

\begin{lemma}\label{funkyfaith}
Let $\calC$ be a stable $\infty$-category, $\calD$ an $\infty$-category which admits
finite limits, and $G: \calC \rightarrow \Stab(\calD)$ an exact functor. Suppose that
$g = \Omega^{\infty} \circ G: \calC \rightarrow \calD$ is fully faithful. Then $G$ is fully faithful.
\end{lemma}

\begin{proof}
It will suffice to show that each of the composite maps
$$ g_{n}: \calC \rightarrow \Stab(\calD) \stackrel{ \Omega^{\infty - n}_{\ast} }{\rightarrow}
\calD_{\ast} $$
is fully faithful. Since $g_{n}$ can be identified with $g_{n+1} \circ \Omega$, where
$\Omega: \calC \rightarrow \calC$ denotes the loop functor, we can reduce to the case
$n = 0$. Fix objects $C, C' \in \calC$; we will show that the map
$\bHom_{\calC}( C,C' ) \rightarrow \bHom_{\calD_{\ast}}( g_0(C), g_0(C') )$ is a homotopy equivalence. We have a homotopy fiber sequence
$$ \bHom_{\calD_{\ast}}( g_0(C), g_0(C') ) \stackrel{\theta}{\rightarrow}
\bHom_{\calD}( g(C), g(C')) \rightarrow \bHom_{\calD}( \ast, g(C') ).$$ 
Here $\ast$ denotes a final object of $\calD$. 
Since $g$ is fully faithful, it will suffice to prove that $\theta$ is a homotopy equivalence. 
For this, it suffices to show that $\bHom_{\calD}( \ast, g(C') )$ is contractible. Since $g$ is left exact, this space can be identified with $\bHom_{\calD}( g(\ast), g(C') )$, where $\ast$ is the final object of $\calC$. Invoking once again our assumption that $g$ is fully faithful, we are reduced to proving that $\bHom_{\calC}( \ast, C')$ is contractible. This follows from the assumption that $\calC$ is pointed (since $\ast$ is also an initial object of $\calC$). 
\end{proof}

\begin{proposition}\label{giraudstable}
Let $\calC$ be an $\infty$-category. The following conditions are equivalent:
\begin{itemize}
\item[$(1)$] The $\infty$-category $\calC$ is presentable and stable.
\item[$(2)$] There exists a presentable, stable $\infty$-category $\calD$ and an accessible left-exact
localization $L: \calD \rightarrow \calC$.
\item[$(3)$] There exists a small $\infty$-category $\calE$ such that $\calC$ is equivalent to
an accessible left-exact localization of $\Fun(\calE, \Spectra)$.
\end{itemize}
\end{proposition}

\begin{proof}
The $\infty$-category $\Spectra$ is stable and presentable, so for every small $\infty$-category
$\calE$, the functor $\infty$-category $\Fun(\calE, \Spectra)$ is also stable (Proposition \ref{expfun}) and presentable (Proposition \toposref{presexp}). This proves $(3) \Rightarrow (2)$. The implication $(2) \Rightarrow (1)$ follows from Lemma \ref{lexloc}. We will complete the proof by showing that $(1) \Rightarrow (3)$.

Since $\calC$ is presentable, there exists a small $\infty$-category $\calE$ and a fully faithful embedding $g: \calC \rightarrow \calP(\calE)$, which admits a left adjoint (Theorem \toposref{pretop} ). Propositions \ref{urtusk21} and \ref{urtusk22} implies that $g$ is equivalent to a composition $$ \calC \stackrel{G}{\rightarrow} \Stab(\calP(\calE)) \stackrel{\Omega^{\infty}}{\rightarrow} \calP(\calE),$$
where the functor $G$ admits a left adjoint. Lemma \ref{funkyfaith} implies that $G$ is fully faithful. It follows that $\calC$ is an (accessible) left exact localization of $\Stab( \calP(\calE) )$. We now invoke Example \ref{expire} to identify $\Stab( \calP(\calE) )$ with $\Fun( \calE, \Spectra)$. 
\end{proof}

\begin{remark}
Proposition \ref{giraudstable} can be regarded as an analogue of Giraud's characterization of
topoi as left exact localizations of presheaf categories (\cite{giraud}). Other variations on this theme include the $\infty$-categorical version of Giraud's theorem (Theorem \toposref{mainchar}) and the Gabriel-Popesco theorem for abelian categories (see \cite{gabriel}). 
\end{remark}
 
\section{Accessible t-Structures}\label{stable16}
 
Let $\calC$ be a stable $\infty$-category. If $\calC$ is presentable, then it is reasonably easy to construct t-structures on $\calC$: for any small collection of objects $\{ X_{\alpha} \}$ of $\calC$, there exists a t-structure {\em generated} by the objects $X_{\alpha}$. More precisely, we have the following result:\index{t-structure!generators for}

\begin{proposition}\label{kura}
Let $\calC$ be a presentable stable $\infty$-category.
\begin{itemize}
\item[$(1)$] If $\calC' \subseteq \calC$ is a full subcategory which is presentable, closed under small colimits, and closed under extensions, then there exists a t-structure on $\calC$ such that
$\calC' = \calC_{\geq 0}$.

\item[$(2)$] Let $\{ X_{\alpha} \}$ be a small collection of objects of $\calC$, and let $\calC'$ be the smallest full subcategory of $\calC$ which contains
each $X_{\alpha}$ and is closed under extensions and small colimits. Then $\calC'$ is presentable.
\end{itemize}
\end{proposition}

\begin{proof}
We will give the proof of $(1)$ and defer the (somewhat technical) proof of $(2)$ until the end of this section. Fix $X \in \calC$, and let $\calC'_{/X}$ denote the fiber product $\calC_{/X} \times_{\calC} \calC'$.
Using Proposition \toposref{horse22}, we deduce that $\calC'_{/X}$ is presentable, so that
it admits a final object $f: Y \rightarrow X$. It follows that composition with $f$ induces
a homotopy equivalence
$$ \bHom_{\calC}(Z,Y) \rightarrow \bHom_{\calC}(Z,X)$$
for each $Z \in \calC'$. Proposition \toposref{testreflect} implies that
$\calC'$ is a colocalization of $\calC$. Since $\calC'$ is stable under extensions, Proposition
\ref{condit} implies the existence of a (uniquely determined) t-structure such that
$\calC' = \calC_{\geq 0}$. 
\end{proof}

\begin{definition}\label{accesstt}
Let $\calC$ be a presentable stable $\infty$-category. We will say that a t-structure on $\calC$ is {\it accessible} if the subcategory $\calC_{\geq 0} \subseteq \calC$ is presentable.
\end{definition}\index{t-structure!accessible}

Proposition \ref{kura} can be summarized as follows: any small collection of objects
$\{ X_{\alpha} \}$ of a presentable stable $\infty$-category $\calC$ determines an accessible t-structure on $\calC$, which is minimal among t-structures such that each $X_{\alpha}$ belongs to $\calC_{\geq 0}$. 

Definition \ref{accesstt} has a number of reformulations:

\begin{proposition}\label{accesst}
Let $\calC$ be a presentable stable $\infty$-category equipped with a t-structure. The following conditions are equivalent:
\begin{itemize}
\item[$(1)$] The $\infty$-category $\calC_{\geq 0}$ is presentable (equivalently: the t-structure on
$\calC$ is accessible).
\item[$(2)$] The $\infty$-category $\calC_{\geq 0}$ is accessible.
\item[$(3)$] The $\infty$-category $\calC_{\leq 0}$ is presentable.
\item[$(4)$] The $\infty$-category $\calC_{\leq 0}$ is accessible.

\item[$(5)$] The truncation functor $\tau_{\leq 0}: \calC \rightarrow \calC$
is accessible.

\item[$(6)$] The truncation functors $\tau_{ \geq 0}: \calC \rightarrow \calC$ is accessible.
\end{itemize}
\end{proposition}

\begin{proof}
We observe that $\calC_{\geq 0}$ is stable under all colimits which exist in $\calC$, and that
$\calC_{\leq 0}$ is a localization of $\calC$. It follows that $\calC_{\geq 0}$ and $\calC_{\leq 0}$ admit small colimits, so that $(1) \Leftrightarrow (2)$ and $(3) \Leftrightarrow (4)$. 
We have a distinguished triangle of functors
$$ \tau_{\geq 0} \stackrel{\alpha}{\rightarrow} \id_{\calC} \stackrel{\beta}{\rightarrow} \tau_{\leq -1} \rightarrow \tau_{\geq 0}[1]$$
in the homotopy category $\h{\Fun(\calC, \calC)}$. The collection of accessible functors from
$\calC$ to itself is stable under shifts and under small colimits. Since $\tau_{\leq 0} \simeq
\coker(\alpha)[1]$ and $\tau_{\geq 0} \simeq \coker(\beta)[-1]$, we conclude that
$(5) \Leftrightarrow (6)$. The equivalence $(1) \Leftrightarrow (5)$ follows from
Proposition \toposref{accloc}. We will complete the proof by showing that
$(1) \Leftrightarrow (3)$.

Suppose first that $(1)$ is satisfied. Then $\calC_{ \geq 1} = \calC_{\geq 0}[1]$ is generated under colimits by a set of objects $\{ X_{\alpha} \}$. Let $S$ be the collection of all morphisms
$f$ in $\calC$ such that $\tau_{\leq 0}(f)$ is an equivalence. Using Proposition \ref{condit},
we conclude that $S$ is generated by $\{ 0 \rightarrow X_{\alpha} \}$ as a quasisaturated class of morphisms, and therefore also as a strongly saturated class of morphisms (Definition \toposref{saturated2}). We now apply Proposition \toposref{local} to conclude that $\calC_{\leq 0} = S^{-1} \calC$ is presentable; this proves $(3)$.

We now complete the proof by showing that $(3) \Rightarrow (1)$. If $\calC_{\leq -1} = \calC_{\leq 0}[-1]$ is presentable, then Proposition \toposref{postbluse}
implies that $S$ is of small generation (as a strongly saturated class of morphisms). Proposition \ref{condit} implies that $S$ is generated (as a strongly saturated class) by the morphisms $\{ 0 \rightarrow X_{\alpha} \}_{\alpha \in A}$, where
$X_{\alpha}$ ranges over the collection of all objects of $\calC_{\geq 0}$. It follows that there
is a small subcollection $A_0 \subseteq A$ such that $S$ is generated by the morphisms
$\{ 0 \rightarrow X_{\alpha} \}_{\alpha \in A_0}$. Let
$\calD$ be the smallest full subcategory of $\calC$ which contains the objects $\{ X_{\alpha} \}_{\alpha \in A_0}$ and is closed under colimits and extensions. Since
$\calC_{\geq 0}$ is closed under colimits and extensions, we have $\calD \subseteq \calC_{ \geq 0}$. Consequently, $\calC_{ \leq -1}$ can be characterized as full subcategory 
of $\calC$ spanned by those objects $Y \in \calC$ such that $\Ext^{k}_{\calC}(X,Y)$
for all $k \leq 0$ and $X \in \calD$. Propositions \ref{kura} implies that $\calD$ is the collection of nonnegative objects for some accessible t-structure on $\calC$. Since the negative objects of this new t-structure coincide with the negative objects of the original t-structure, we conclude that
$\calD = \calC_{\geq 0}$, which proves $(1)$.
\end{proof}

The following result provides a good source of examples of accessible t-structures:

\begin{proposition}\label{tmonster}
Let $\calC$ be a presentable $\infty$-category, and let $\Stab(\calC)_{\leq -1}$ be the full subcategory of $\Stab(\calC)$ spanned by those objects $X$ such that 
$\Omega^{\infty}(X)$ is a final object of $\calC$. Then $\Stab(\calC)_{\leq -1}$ determines
an accessible t-structure on $\Stab(\calC)$.
\end{proposition}

\begin{proof}
Choose a small collection of objects $\{ C_{\alpha} \}$ which generate $\calC$ under colimits. 
We observe that an object $X \in \Stab(\calC)$ belongs to $\Stab(\calC)_{\leq -1}$ if and only if
each of the spaces
$$ \bHom_{\calC}( C_{\alpha}, \Omega^{\infty}(X) ) \simeq \bHom_{ \Stab(\calC)}( \Sigma^{\infty}(C_{\alpha}), X) $$
is contractible. Let $\Stab(\calC)_{\geq 0}$ be the smallest full subcategory of
$\Stab(\calC)$ which is stable under colimits and extensions, and contains each $\Sigma^{\infty}(C_{\alpha})$. Proposition \ref{kura} implies that $\Stab(\calC)_{\geq 0}$ is the collection of nonnegative objects of the desired t-structure on $\Stab(\calC)$.
\end{proof}

\begin{remark}
The proof of Proposition \ref{tmonster} gives another characterization of the t-structure on
$\Stab(\calC)$: the full subcategory $\Stab(\calC)_{\geq 0}$ is generated, under extensions and colimits, by the essential image of the functor $\Sigma^{\infty}: \calC \rightarrow \Stab(\calC)$.
\end{remark}

We conclude this section by completing the proof of Proposition \ref{kura}.

\begin{proof}[Proof of part $(2)$ of Proposition \ref{kura}]
Choose a regular cardinal $\kappa$ such that every object of $X_{\alpha}$ is $\kappa$-compact, and let $\calC^{\kappa}$ denote the full subcategory of $\calC$ spanned by the $\kappa$-compact objects. Let ${\calC'}^{\kappa} = \calC' \cap \calC^{\kappa}$, and let $\calC''$ be the smallest
full subcategory of $\calC'$ which contains ${\calC'}^{\kappa}$ and is closed under small colimits.
The $\infty$-category $\calC''$ is $\kappa$-accessible, and therefore presentable.
To complete the proof, we will show that $\calC' \subseteq \calC''$. For this, it will suffice to show that $\calC''$ is stable under extensions.

Let $\calD$ be the full subcategory of $\Fun(\Delta^1, \calC)$ spanned by those morphisms
$f: X \rightarrow Y$ where $Y \in \calC''$, $X \in \calC''[-1]$. We wish to prove that the cokernel functor $\coker: \calD \rightarrow \calC$ factors through $\calC''$. Let $\calD^{\kappa}$ be the full subcategory of $\calD$ spanned by those morphisms $f: X \rightarrow Y$ where both $X$ and $Y$ are $\kappa$-compact objects of $\calC$. By construction, $\coker | \calD^{\kappa}$ factors through $\calC''$. Since $\coker: \calD \rightarrow \calC$ preserves small colimits, it will suffice to show that
$\calD$ is generated (under small colimits) by $\calD^{\kappa}$.

Fix an object $f: X \rightarrow Y$ in $\calD$. To complete the proof, it will suffice to show that
the canonical map $(\calD^{\kappa}_{/f})^{\triangleright} \rightarrow \calD$ is a colimit diagram.
Since $\calD$ is stable under colimits in $\Fun(\Delta^1, \calC)$ and colimits in $\Fun(\Delta^1,\calC)$ are computed pointwise (Proposition \toposref{limiteval}), it will suffice to show that composition with the evaluation maps give colimit diagrams
$(\calD^{\kappa}_{/f})^{\triangleright} \rightarrow \calC$. Lemma \toposref{waitlong1} implies that
the maps $({\calC'}^{\kappa}[-1])_{/X}^{\triangleright} \rightarrow \calC$, $({\calC'}^{\kappa})_{/Y}^{\triangleright} \rightarrow \calC$ are colimit diagrams. It will therefore suffice to show that
the evaluation maps
$$ ({\calC'}^{\kappa}[-1])_{/X} \stackrel{\theta}{\leftarrow} (\calD^{\kappa}_{/f}) \stackrel{\theta'}{\rightarrow} ({\calC'}^{\kappa})_{/Y}$$
are cofinal.

We first show that $\theta$ is cofinal. According to Theorem \toposref{hollowtt}, it will suffice to show that for every morphism $\alpha: X' \rightarrow X$ in $\calC'[-1]$, where $X'$ is $\kappa$-compact, the $\infty$-category
$$ \calE_{\theta}: \calD^{\kappa}_{/f} \times_{ {\calC'}^{\kappa}[-1]_{/X} } ( {\calC'}^{\kappa}[-1]_{/X} )_{X'/} $$
is weakly contractible. For this, it is sufficient to show that $\calE_{\theta}$ is filtered (Lemma \toposref{stull2}). 

We will show that $\calE_{\theta}$ is $\kappa$-filtered. Let $K$ be a $\kappa$-small simplicial set, and $p: K \rightarrow \calE_{\theta}$ a diagram; we will extend $p$ to a diagram
$\overline{p}: K^{\triangleright} \rightarrow \calE_{\theta}$. We can identify $p$ with two pieces of data:
\begin{itemize}
\item[$(i)$] A map $p': K^{\triangleleft} \rightarrow {\calC'}^{\kappa}[-1]_{/X}$.
\item[$(ii)$] A map $p'': ( K \star \{ \infty \} ) \times \Delta^1 \rightarrow \calC$, with the properties
that $p'' | (K \star \{ \infty \}) \times \{0\}$ can be identified with $p'$, 
$p'' | \{ \infty \} \times \Delta^1$ can be identified with $f$, and
$p'' | K \times \{1\}$ factors through ${ \calC' }^{\kappa}$. 
\end{itemize}
Let $\overline{p}': (K^{\triangleleft})^{\triangleright} \rightarrow {\calC'}^{\kappa}[-1]_{/X}$ be a colimit of $p'$. To complete the proof that $\calE_{\theta}$ is $\kappa$-filtered, it will suffice to show that
we can find a compatible extension $\overline{p}'': (K^{\triangleright} \star \{ \infty \}) \times \Delta^1 \rightarrow \calC$ with the appropriate properties. Let $L$ denote the full simplicial subset
of $(K^{\triangleright} \star \{\infty \}) \times \Delta^1)$ spanned by every vertex except 
$(v, 1)$, where $v$ denotes the cone point of $K^{\triangleright}$. We first choose a map
$q: L \rightarrow \calC$ compatible with $p''$ and $\overline{p}'$. This is equivalent
to solving the lifting problem
$$ \xymatrix{ & \calC_{/f} \ar[d] \\
K^{\triangleright} \ar[r] \ar@{-->}[ur] & \calC_{/X}, }$$
which is possible since the vertical arrow is a trivial fibration. Let $L' = L \cap (K^{\triangleright} \times \Delta^1)$. Then $q$ determines a map $q_0: L' \rightarrow \calC_{/Y}$. Finding the
desired extension $\overline{p}''$ is equivalent to finding a map
$\overline{q}_0: {L'}^{\triangleright} \rightarrow \calC_{/Y}$, which carries the cone point
into ${\calC'}^{\kappa}$. 

Let $g: Z \rightarrow Y$ be a colimit of $q_0$
(in the $\infty$-category $\calC_{/Y}$). We observe that $Z$ is a $\kappa$-small colimit of $\kappa$-compact objects of $\calC$, and therefore $\kappa$-compact. Since $Y \in \calC''$, $Y$ can be written as the colimit of a $\kappa$-filtered diagram $\{ Y_{\alpha} \}$, taking values in ${\calC'}^{\kappa}$. Since $Z$ is $\kappa$-compact, the map $g$ factors through some $Y_{\alpha} $; it follows that there exists an extension $\overline{q}_0$ as above, which carries the cone point to $Y_{\alpha}$. This completes the proof that $\calE_{\theta}$ is $\kappa$-filtered, and also the proof that $\theta$ is cofinal.

The proof that $\theta'$ is cofinal is similar but slightly easier: it suffices to show that for every
map $Y' \rightarrow Y$ in $\calC'$, where $Y'$ is $\kappa$-compact, the fiber product 
$$\calE_{\theta'} =  \calD^{\kappa}_{/f} \times_{ {\calC'}^{\kappa}_{/Y} } ( {\calC'}^{\kappa}_{/Y} )_{Y'/} $$
is filtered. For this, we can either argue as above, or simply observe that $\calE_{\theta'}$ admits $\kappa$-small colimits.
\end{proof}